\documentstyle[amscd,amssymb,xypic,11pt]{amsart}
\xyoption{all}

\newcommand{\nc}{\newcommand}

\nc{\exto}[1]{\stackrel{#1}{\longrightarrow}}
\nc{\dlim}{{\mathop{\lim\limits_{\longrightarrow}}}}
\nc{\lan}{\big\langle}
\nc{\ran}{\big\rangle}

\nc{\kk}{{\mathsf{k}}}
\nc{\ix}{{\mathsf{i}}}
\nc{\jx}{{\mathsf{j}}}

\nc{\C}{{\mathbb{C}}}
\nc{\HH}{{\mathbb{H}}}
\nc{\PP}{{\mathbb{P}}}
\nc{\QQ}{{\mathbb{Q}}}
\nc{\ZZ}{{\mathbb{Z}}}

\nc{\CA}{{\mathcal{A}}}
\nc{\CB}{{\mathcal{B}}}
\nc{\CC}{{\mathcal{C}}}
\nc{\D}{{\mathcal{D}}}
\nc{\CE}{{\mathcal{E}}}
\nc{\CF}{{\mathcal{F}}}
\nc{\CG}{{\mathcal{G}}}
\nc{\CH}{{\mathcal{H}}}
\nc{\CL}{{\mathcal{L}}}
\nc{\CM}{{\mathcal{M}}}
\nc{\CN}{{\mathcal{N}}}
\nc{\CO}{{\mathcal{O}}}
\nc{\CQ}{{\mathcal{Q}}}
\nc{\CR}{{\mathcal{R}}}
\nc{\CS}{{\mathcal{S}}}
\nc{\CT}{{\mathcal{T}}}
\nc{\CU}{{\mathcal{U}}}
\nc{\CV}{{\mathcal{V}}}
\nc{\CW}{{\mathcal{W}}}
\nc{\CX}{{\mathcal{X}}}
\nc{\CY}{{\mathcal{Y}}}
\nc{\CMo}{{\mathcal{M}^\circ}}
\nc{\Co}{{{C}^\circ}}

\nc{\BY}{{\overline{Y}}}
\nc{\BYD}{{\overline{Y}{}^{|D|}}}
\nc{\OZ}{{\overline{Z}}}
\nc{\bg}{{\bar{g}}}

\nc{\bq}{{\mathbf{q}}}
\nc{\BD}{{\mathbf{D}}}
\nc{\BG}{{\mathbf{G}}}
\nc{\BM}{{\mathbf{M}}}
\nc{\BP}{{\mathbf{P}}}
\nc{\BZ}{{\mathbf{Z}}}
\nc{\BPr}{{\mathsf{P}}}
\nc{\BR}{{\mathbf{R}}}
\nc{\BRO}[1]{{{\mathbf{R}}^{\circ}_{#1}}}
\nc{\BRD}[1]{{{\mathbf{R}}^{|D|}_{#1}}}
\nc{\BRP}[1]{{{\mathbf{R}}^{1}_{#1}}}
\nc{\BRTP}[1]{{{\mathbf{\tilde{R}}}{}^{1}_{#1}}}
\nc{\BS}{{\mathbf{S}}}
\nc{\BMS}{{{\mathbf{M}}^{{s}}}}
\nc{\BMSS}{{{\mathbf{M}}^{{ss}}}}
\nc{\BMZ}{{\mathbf{M}^{\circ}}}
\nc{\BCL}{{\mathbf{L}}}

\nc{\PCC}{{{}^\perp\CC}}

\nc{\Cl}{{\mathsf{Cliff}}}
\nc{\Clev}{{\mathop{\mathsf{Cliff}}^{\circ}}}

\nc{\FA}{{\mathfrak{A}}}
\nc{\FB}{{\mathfrak{B}}}

\nc{\fa}{{\mathfrak{a}}}
\nc{\fb}{{\mathfrak{b}}}
\nc{\fg}{{\mathfrak{g}}}
\nc{\fn}{{\mathfrak{n}}}
\nc{\fp}{{\mathfrak{p}}}
\nc{\FD}{{\mathfrak{D}}}
\nc{\FE}{{\mathfrak{E}}}
\nc{\FL}{{\mathfrak{L}}}
\nc{\FM}{{\mathfrak{M}}}
\nc{\FS}{{\mathsf{S}}}

\nc{\sfc}{{\mathsf{c}}}
\nc{\sfch}{{\mathsf{ch}}}
\nc{\sfh}{{\mathsf{h}}}

\nc{\SK}{{\mathsf{K}}}
\nc{\SO}{{\mathsf{O}}}
\nc{\SQ}{{\mathsf{Q}}}
\nc{\SPV}{{\mathsf{S}^+\mathsf{V}}}
\nc{\SMV}{{\mathsf{S}^-\mathsf{V}}}
\nc{\SPMV}{{\mathsf{S}^\pm\mathsf{V}}}
\nc{\SX}{{S_X}}
\nc{\SY}{{S_Y}}
\nc{\phipsi}{{q}}
\nc{\eps}{\varepsilon}

\nc{\pim}{{\pi_-}}
\nc{\pip}{{\pi_+}}

\nc{\BE}{{\overline{\CE}}}
\nc{\TE}{{\tilde{\CE}}}
\nc{\TQ}{{\tilde{Q}}}
\nc{\TCF}{{\tilde{\CF}}}
\nc{\TCG}{{\tilde{\CG}}}
\nc{\TCL}{{\tilde{\CL}}}
\nc{\TF}{{\tilde{F}}}
\nc{\TW}{{\tilde{W}}}
\nc{\TCC}{{\tilde{\CC}}}
\nc{\TCX}{{\tilde{\CX}}}
\nc{\TCY}{{\tilde{\CY}}}
\nc{\TPhi}{{\tilde{\Phi}}}
\nc{\OPhi}{{\bar{\Phi}}}
\nc{\txi}{{\tilde{\xi}}}
\nc{\tp}{{\tilde{p}}}
\nc{\tq}{{\tilde{q}}}
\nc{\tzeta}{{\tilde{\zeta}}}
\nc{\tpi}{{\tilde{\pi}}}

\nc{\HE}{{\widehat{\CE}}}
\nc{\HX}{{\hat{X}}}
\nc{\hxi}{{\hat{\xi}}}

\nc{\UH}{{\mathcal{H}}}

\nc{\TM}{{\widetilde{M}}}
\nc{\TCM}{{\widetilde{\CM}}}
\nc{\TU}{{\widetilde{U}}}
\nc{\TX}{{\widetilde{X}}}
\nc{\TY}{{\widetilde{Y}}}
\nc{\TYO}{{{\widetilde{Y}}^\circ}}
\nc{\barf}{{\bar{f}}}
\nc{\te}{{\tilde{e}}{}}
\nc{\tf}{{\tilde{f}}}
\nc{\tg}{{\tilde{g}}}
\nc{\ti}{{\tilde{\imath}}}
\nc{\tj}{{\tilde{\jmath}}}
\nc{\ty}{{\tilde{y}}}
\nc{\tphi}{{\tilde{\phi}}}

\nc{\urho}{{\underline{\rho}}}

\nc{\LRA}{\Leftrightarrow}
\nc{\RA}{\Rightarrow}
\nc{\lotimes}{\mathbin{\mathop{\otimes}\limits^{\mathbb{L}}}}
\nc{\CEnd}{\mathop{\mathcal{E}\mathit{nd}}\nolimits}
\nc{\CExt}{\mathop{\mathcal{E}\mathit{xt}}\nolimits}
\nc{\CHom}{\mathop{\mathcal{H}\mathit{om}}\nolimits}
\nc{\RH}{\mathop{{\mathsf{R}}\Gamma}\nolimits}
\nc{\RGamma}{\mathop{{\mathsf{R}}\Gamma}\nolimits}
\nc{\RHom}{\mathop{\mathsf{RHom}}\nolimits}
\nc{\RCHom}{\mathop{\mathsf{R}\mathcal{H}\mathit{om}}\nolimits}
\nc{\RG}{\mathop{\mathsf{R\Gamma}}\nolimits}
\nc{\Hom}{\mathop{\mathsf{Hom}}\nolimits}
\nc{\Ext}{\mathop{\mathsf{Ext}}\nolimits}
\nc{\End}{\mathop{\mathsf{End}}\nolimits}
\nc{\Tor}{\mathop{\mathsf{Tor}}\nolimits}
\nc{\Tordim}{\mathop{\mathsf{Tor}\text{\rm-}\mathsf{dim}}\nolimits}
\nc{\Hilb}{\mathop{\mathsf{Hilb}}\nolimits}
\nc{\Spec}{\mathop{\mathsf{Spec}}\nolimits}
\nc{\Pic}{\mathop{\mathsf{Pic}}\nolimits}
\renewcommand{\Im}{\mathop{\mathsf{Im}}\nolimits}
\nc{\Tw}{\mathop{\mathsf{Tw}}\nolimits}
\nc{\Ker}{\mathop{\mathsf{Ker}}\nolimits}
\nc{\Coker}{\mathop{\mathsf{Coker}}\nolimits}
\nc{\codim}{\mathop{\mathsf{codim}}\nolimits}
\nc{\sing}{{\mathsf{sing}}}
\nc{\supp}{\mathop{\mathsf{supp}}}
\nc{\perf}{{\mathsf{perf}}}
\nc{\rank}{\mathop{\mathsf{rank}}}
\nc{\Pf}{{\mathsf{Pf}}}
\nc{\Gr}{{\mathsf{Gr}}}
\nc{\OGr}{{\mathsf{OGr}}}
\nc{\Flag}{{\mathsf{Fl}}}
\nc{\Kosz}{{\mathsf{Kosz}}}
\nc{\LGr}{{\mathsf{LGr}}}
\nc{\GTGr}{{\mathsf{G_2Gr}}}
\nc{\GTF}{{\mathsf{G_2F}}}
\nc{\OF}{{\mathsf{OF}}}
\nc{\Fl}{{\mathsf{Fl}}}
\nc{\Bl}{{\mathsf{Bl}}}
\nc{\GL}{{\mathsf{GL}}}
\nc{\PGL}{{\mathsf{PGL}}}
\nc{\SL}{{\mathsf{SL}}}
\nc{\SP}{{\mathsf{Sp}}}
\nc{\Spin}{{\mathsf{Spin}}}
\nc{\Tot}{{\mathsf{Tot}}}
\nc{\ev}{{\mathsf{ev}}}
\nc{\od}{{\mathsf{odd}}}
\nc{\coev}{{\mathsf{coev}}}
\nc{\id}{{\mathsf{id}}}
\nc{\opp}{{\mathsf{opp}}}
\nc{\PS}{{{\PP^3}}}
\nc{\Qu}{{{Q^3}}}
\nc{\tdim}{\mathop{\Tor\dim}}
\nc{\ecart}{{\fbox{$\scriptstyle\mathsf{EC}$}}}
\nc{\ad}{{\mathop{\mathsf ad}}}
\nc{\sg}{{\mathop{\mathsf sg}}}
\nc{\hf}{{\mathop{\mathsf hf}}}
\nc{\gr}{{\mathop{\mathsf gr}}}
\nc{\qgr}{{\mathop{\mathsf qgr}}}
\nc{\Coh}{{\mathop{\mathsf Coh}}}
\nc{\Ab}{{\mathop{\mathcal{A}\mathit{b}}}}
\nc{\Qcoh}{{\mathop{\mathsf Qcoh}}}

\nc{\AAV}{{\mathcal{AAV}}}

\nc{\Rep}{{\mathsf{Rep}}}

\nc{\Cubics}{{{\mathcal{S}}_3}}
\nc{\VFT}{{{\mathcal{S}}_{14}}}
\nc{\VFTE}{{{\mathcal{N}}_{\mathrm{reg,sm}}}}
\nc{\MX}{{\CM_X}}
\nc{\MY}{{\CM_Y}}
\nc{\MYE}{{\CM_{Y,\CE}}}
\nc{\Yd}{{Y_d}}
\nc{\Yfive}{{Y_5}}
\nc{\Xg}{{X_{2g-2}}}
\nc{\Xtt}{{X_{22}}}
\nc{\Xst}{{X_{16}}}
\nc{\Xtw}{{X_{12}}}
\nc{\Xe}{{X_{8}}}
\nc{\Xf}{{X_{4}}}

\nc{\git}{{/\!\!/\!{}_\chi}}

\theoremstyle{plain}

\newtheorem{theorem}{Theorem}[section]
\newtheorem{conjecture}[theorem]{Conjecture}
\newtheorem{lemma}[theorem]{Lemma}
\newtheorem{proposition}[theorem]{Proposition}
\newtheorem{corollary}[theorem]{Corollary}

\theoremstyle{definition}

\newtheorem{definition}[theorem]{Definition}

\theoremstyle{remark}

\newtheorem{remark}[theorem]{Remark}

\newenvironment{proof}{\noindent{\sf Proof:}}{\qed\medskip}
\newenvironment{proofof}[1]{\noindent{\sf Proof #1:}}{\qed\medskip}

\title{Homological Projective Duality}
\author{Alexander Kuznetsov}
\address{
Algebra Section, Steklov Mathematical Institute,
8 Gubkin str., Moscow 119991 Russia
}
\email{akuznet@@mi.ras.ru}
\date{}
\thanks{I was partially supported by RFFI grants 05-01-01034 and 02-01-01041,
Russian Presidential grant for young scientists No. MK-3926.2004.1,
CRDF Award No. RUM1-2661-MO-05, and the Russian Science Support Foundation.}

\begin{document}

\begin{abstract}
We introduce a notion of Homological Projective Duality for smooth algebraic
varieties in dual projective spaces, a homological extension of the classical
projective duality. If algebraic varieties $X$ and $Y$ in dual projective
spaces are Homologically Projectively Dual, then we prove that the orthogonal
linear sections of $X$ and $Y$ admit semiorthogonal decompositions with
an equivalent nontrivial component. In particular, it follows that
triangulated categories of singularities of these sections are equivalent.
We also investigate Homological Projective Duality for projectivizations
of vector bundles.
\end{abstract}

\maketitle


\section{Introduction}

Investigation of derived categories of coherent sheaves on algebraic varieties
became one of the most important topics in the modern algebraic geometry.
Besides other reasons this is caused by the Homological Mirror Symmetry
conjecture of Maxim Kontsevich \cite{Ko} predicting that there is
an equivalence of categories between the derived category of coherent
sheaves on a Calabi--Yau variety and the derived Fukai category
of its mirror. There is an extension of the Mirror Symmetry
to the non Calabi-Yau case \cite{HV}. According to this,
the mirror of an arbitrary variety is a so-called
Landau--Ginzburg model, that is an algebraic variety with a 2-form
and a holomorphic function (superpotential) such that the restriction
of the 2-form to smooth fibers of the superpotential is symplectic.
It is expected that singular fibers of the superpotential of the mirror
Landau--Ginzburg model give a decomposition of the derived category
of coherent sheaves on the initial algebraic variety into semiorthogonal
pieces, a semiorthogonal decomposition.


Thus from the point of view of mirror symmetry it is important
to investigate when the derived category of coherent sheaves on a variety
admits a semiorthogonal decomposition.
The goal of the present paper is to answer the following more precise question:
$$
\parbox{0.7\textwidth}{\em
Assume that $X$ is a smooth projective variety and denote by $\D^b(X)$
the bounded derived category of coherent sheaves on $X$.
Supposing that we are given a semiorthogonal decomposition of $\D^b(X)$,
is it possible to construct a semiorthogonal decomposition of $\D^b(X_H)$,
where $X_H$ is a hyperplane section of $X$?
}\leqno{(\dagger)}
$$
Certainly this question is closely related to the question
{\em what does the operation of taking a hyperplane section of a projective
algebraic variety mean on the other side of the mirror?}

In general one cannot expect an affirmative answer on ($\dagger$).
However, there is an important particular case, when something can be said.
Explicitly, assume that $X \subset \PP(V)$ is a smooth projective variety,
$\CO_X(1)$ is the corresponding very ample line bundle, and assume that
there is a semiorthogonal decomposition of its derived category
of the following type
$$
\D^b(X) = \langle \CA_0,\CA_1(1),\dots,\CA_{\ix-1}(\ix-1) \rangle,
\qquad
0 \subset \CA_{\ix-1} \subset \CA_{\ix-2} \subset \dots \subset \CA_1 \subset \CA_0,
$$
where $(k)$ stands for the twist by $\CO_X(k)$.
A decomposition of this type will be called {\sf Lefschetz decomposition}
because as we will see its behavior with respect to hyperplane sections
is similar to that of the Lefschetz decomposition of the cohomology groups.
An easy calculation shows that for any hyperplane section $X_H$ of $X$
with respect to $\CO_X(1)$ the composition of the embedding and the restriction
functors $\CA_k(k) \to \D^b(X) \to \D^b(X_H)$ is fully faithful and
$\lan \CA_1(1),\dots,\CA_{\ix-1}(\ix-1) \ran$ is a semiorthogonal collection
in $\D^b(X_H)$. In other words, dropping the first (the biggest) component
of the Lefschetz decomposition of $\D^b(X)$ we obtain a semiorthogonal
collection in $\D^b(X_H)$. Denoting by $\CC_H$ the orthogonal in $\D^b(X_H)$
to the subcategory of $\D^b(X_H)$ generated by this collection we consider
$\{\CC_H\}_{H\in\PP(V^*)}$ as a family of triangulated categories
over the projective space $\PP(V^*)$. Assuming geometricity of this family,
i.e. roughly speaking that there exists an algebraic variety $Y$ with a map
$Y \to \PP(V^*)$ such that  for all $H$ we have $\CC_H \cong \D^b(Y_H)$,
where $Y_H$ is the fiber of $Y$ over $H\in\PP(V^*)$,
we prove the main result of the paper

\begin{theorem}\label{int_th}
The derived category of $Y$ admits a dual Lefschetz decomposition
$$
\D^b(Y) = \lan \CB_{\jx-1}(1-\jx),\CB_{\jx-2}(2-\jx),\dots,\CB_{1}(-1),\CB_0\ran,
\qquad
0 \subset \CB_{\jx-1} \subset \CB_{\jx-2} \subset \dots \subset \CB_1 \subset \CB_0.
$$
Moreover, if $L\subset V^*$ is a linear subspace and
$L^\perp \subset V$ is the orthogonal subspace
such that the linear sections
$X_L = X\times_{\PP(V)}\PP(L^\perp)$ and
$Y_L = Y\times_{\PP(V^*)}\PP(L)$,
are of expected dimension
$\dim X_L = \dim X - \dim L$,
and
$\dim Y_L = \dim Y - \dim L^\perp$,
then there exists a triangulated category $\CC_L$ and
semiorthogonal decompositions
$$
\begin{array}{lll}
\D^b(X_L) &=& \langle \CC_L,\CA_{\dim L}(1),\dots,\CA_{\ix-1}(\ix-\dim L)\rangle\\
\D^b(Y_L) &=& \langle \CB_{\jx-1}(\dim L^\perp - \jx),\dots,\CB_{\dim L^\perp}(-1),\CC_L\rangle.
\end{array}
$$
\end{theorem}

In other words, the derived categories of $X_L$ and $Y_L$ have semiorthogonal
decompositions with several ``trivial'' components coming from the Lefschetz
decompositions of the ambient varieties $X$ and $Y$ respectively,
and with {\em equivalent}\/ nontrivial components.

We would like to emphasize a similarity in behavior of derived categories
and cohomology groups with respect to the hyperplane section operation.
Thus, theorem~\ref{int_th} can be considered as a homological generalization
of the Lefschetz Theorem about hyperplane sections.

\medskip

A simple corollary of theorem~\ref{int_th} is an equivalence of the derived
categories of singularities (see \cite{O4}) of $X_L$ and $Y_L$. In particular,
it easily follows that $Y_L$ is singular if and only if $X_L$ is singular.
This means that we have an equality of the following two closed subsets
of the dual projective space $\PP(V^*)$:
$$
\{ H \in \PP(V^*)\ |\ \text{$X_H$ is singular} \} =
\{ \text{critical values of the projection $Y \to \PP(V^*)$}\}
$$
Note that the first of these subsets is the classical projectively dual variety of $X$.
Thus $Y$ can be considered as a homological generalization of the projectively dual.
In accordance with this we say that $Y$ is a {\sf Homologically Projectively Dual}\/ variety of $X$.

\medskip

The simplest example of a Lefschetz decomposition is given by
the standard exceptional collection $(\CO,\CO(1),\dots,\CO(\ix-1))$
on a projective space $X = \PP^{\ix-1}$
(we take $\CA_0 = \CA_1 = \dots = \CA_{\ix-1} = \lan \CO \ran$).
It is easy to see that the corresponding Homological Projectively
Dual variety is an empty set, and we obtain nothing interesting.
However, considering a {\em relative}\/ projective space
we already obtain some interesting results.
More precisely, consider a projectivization of a vector bundle $X = \PP_S(E)$
over a base scheme $S$, embedded into the projectivization
of the vector space $H^0(S,E^*)^* = H^0(X,\CO_{X/S}(1))^*$
with the following Lefschetz decomposition
$$
\D^b(X) = \lan\D^b(S),\D^b(S)\otimes\CO_{X/S}(1),\dots,\D^b(S)\otimes\CO_{X/S}(\ix-1)\ran.
$$
We prove that $Y = \PP_S(E^\perp)$, where $E^\perp = \Ker(V^*\otimes\CO_S \to E^*)$
is a Homologically Projectively Dual variety of $X$.
As a consequence we get a bunch of semiorthogonal decompositions and
equivalences between derived categories of linear sections of $\PP_S(E)$
and $\PP_S(E^\perp)$. For example, applying a relative version of theorem~\ref{int_th}
we can deduce that there is an equivalence of derived categories between the following
two varieties related by a special birational transformation called a {\em flop}\/
(it is conjectured in \cite{BO2} that the derived categories of any pair
of algebraic varieties related by a flop are equivalent).
Consider a morphism of vector bundles $F \exto{\phi} E^*$ of equal ranks on $S$
and consider
$$
X_F = \{(s,e) \in \PP_S(E)\ |\ \phi^*_s(e) = 0\},
\quad
Y_F = \{(s,f) \in \PP_S(F)\ |\ \phi_s(f) = 0\},
\ \text{and}\
Z_F = \{s\in S\ |\ \det\phi_s = 0\}.
$$
If $\dim X_F = \dim Y_F = \dim S - 1$ then the natural projections
$X_F \to Z_F$ and $Y_F \to Z_F$ are birational and the corresponding
birational transformation $\xymatrix@1{X_F \ar@{..>}[r] & Y_F}$ is a flop.
We prove an equivalence of categories $\D^b(X_F) \cong \D^b(Y_F)$ if
additionally $\dim X_F\times_S Y_F = \dim S - 1$.

The next example of a Lefschetz decomposition is a decomposition
of $\D^b(X)$ for $X = \PP(W)$ with respect to $\CO_X(2)$ given
by $\CA_0 = \CA_1 = \dots = \CA_{\ix-2} = \lan\CO_X,\CO_X(1)\ran$
and either $\CA_{\ix-1} = \lan\CO_X,\CO_X(-1)\ran$ for $\dim W = 2\ix$,
or $\CA_{\ix-1} = \lan\CO_X\ran$ for $\dim W = 2\ix-1$.
In a companion paper~\cite{K4} we show that the universal sheaf
of even parts of Clifford algebras on $\PP(S^2W^*)$ is
a Homologically Projectively Dual variety to $X$ with respect
to the double Veronese embedding $\PP(W) \subset \PP(S^2W)$.
This gives immediately a proof of the theorem of Bondal and
Orlov \cite{BO2,BO3} about derived categories of intersections
of quadrics.

Finally, let us mention that the Homological Projective Duality
for Lefschetz decompositions with $\CA_0$ generated by exceptional pair
and $\CA_0 = \CA_1 = \dots = \CA_{\ix-1}$ was considered in \cite{K}.
There were constructed such decompositions for
$X = \Gr(2,5)$,
$X = \OGr_+(5,10)$, a connected component of the Grassmannian
of 5-dimensional subspaces in $\kk^{10}$ isotropic with respect to
a nondegenerate quadratic form,
$X = \LGr(3,6)$, the Lagrangian Grassmannian of 3-dimensional
subspaces in $\kk^6$ with respect to a symplectic form, and
$X = \GTGr(2,7)$, the Grassmannian of the Lie group ${\mathsf G}_2$,
and it was shown that Homologically Projectively Dual varieties for them are
$Y = \Gr(2,5)$, $Y = \OGr_-(5,10)$, a quartic hypersurface in $\PP^{13}$,
and a double covering of $\PP^{13}$ ramified in a sextic hypersurface
(in the last two cases one must consider the derived category of sheaves
of modules over a suitable sheaf of Azumaya algebras on~$Y$ instead of the usual
derived category).


\medskip

Now we describe the structure of the paper.
In section~2 we recall the necessary material concerning admissible
subcategories, semiorthogonal decompositions, remind an important
technical result, the faithful base change theorem proved in~\cite{K},
and check that a property of a linear over a base functor to be fully faithful
is local over the base.
In section~3 we define {\em splitting functors}\/ and give a criterion
for a functor to be splitting.
In section~4 we define {\em Lefschetz decompositions}\/ of triangulated categories.
In section~5 we consider derived category of the {\em universal hyperplane section}\/
of a variety admitting a Lefschetz decomposition of the derived category.
In section~6 we define {\em Homological Projective Duality}\/ and
prove theorem~\ref{int_th} and its relative versions.
In section~7 we discuss relation of the Homological Projective Duality
to the classical projective duality.
In section~8 we consider the Homological Projective Duality
for a projectivization of a vector bundle.
Finally, in section~9 we consider some explicit examples
of Homological Projective Duality.

{\bf Acknowledgements.}
I am grateful to A.~Bondal, D.~Kaledin and D.~Orlov for many useful
discussions.
Also I would like to mention that an important example of Homological
Projective Duality (the case of $X = \Gr(2,6)$ which is not discussed
in this paper) first appeared in a conversation with A.~Samokhin.

\section{Preliminaries}

\subsection{Notation}

The base field $\kk$ is assumed to be algebraically closed of zero characteristic.
All algebraic varieties are assumed to be embeddable (i.e. admitting a finite morphism
onto a smooth algebraic variety) and of finite type over $\kk$.

Given an algebraic variety $X$ we denote by $\D^b(X)$ the bounded derived category of
coherent sheaves on $X$. Similarly, $\D^-(X)$, $\D^+(X)$ and $\D(X)$ stand for
the bounded above, the bounded below and the unbounded derived categories.
Further, $\D^b_{qc}(X)$, $\D^-_{qc}(X)$, $\D^+_{qc}(X)$, and $\D_{qc}(X)$,
stand for the corresponding derived categories of quasicoherent sheaves,
and $\D^\perf(X)$ denotes the category of perfect complexes on $X$,
i.e. the full subcategory of $\D(X)$ consisting of all objects
locally quasiisomorphic to bounded complexes of locally free sheaves
of finite rank.

Given a morphism $f:X \to Y$ we denote by $f_*$ and $f^*$
the {\em total}\/ derived pushforward and
the {\em total}\/ derived pullback functors.
The twisted pullback functor \cite{H} is denoted by $f^!$
(it is right adjoint to $f_*$).
Similarly, $\otimes$ stands for the derived tensor product,
and $\RHom$, $\RCHom$ stand for the global and local $\RHom$ functors.

Given an object $F \in \D(X)$ we denote by $\CH^k(F)$ the $k$-th
cohomology sheaf of $F$.

\subsection{Semiorthogonal decompositions}


If $\CA$ is a full subcategory of $\CT$ then the {\sf right orthogonal}\/ to $\CA$ in $\CT$
(resp.\ the {\sf left orthogonal}\/ to $\CA$ in $\CT$) is the full subcategory
$\CA^\perp$ (resp.\ ${}^\perp\CA$) consisting of all objects $T\in\CT$
such that $\Hom_\CT(A,T) = 0$ (resp.\ $\Hom_\CT(T,A) = 0$) for all $A \in \CA$.

\begin{definition}[\cite{BO1}]
A semiorthogonal decomposition of $\CT$ is a sequence of full subcategories
$\CA_1,\dots,\CA_n$ in $\CT$ such that $\Hom_{\CT}(\CA_i,\CA_j) = 0$ for $i > j$
and for every object $T \in \CT$ there exists a chain of morphisms
$0 = T_n \to T_{n-1} \to \dots \to T_1 \to T_0 = T$ such that
the cone of the morphism $T_k \to T_{k-1}$ is contained in $\CA_k$
for each $k=1,2,\dots,n$.
%
%
\end{definition}

In other words, every object $T$ admits a ``filtration''
with factors in $\CA_1$, \dots, $\CA_n$ respectively.
Semiortho\-gonality implies that this filtration is unique
and functorial.

For any sequence of subcategories $\CA_1,\dots,\CA_n$ in $\CT$ we denote
by $\lan\CA_1,\dots,\CA_n\ran$ the minimal triangulated subcategory of $\CT$
containing $\CA_1$, \dots, $\CA_n$.

If $\CT = \lan\CA_1,\dots,\CA_n\ran$ is a semiorthogonal decomposition then
$\CA_i = \lan\CA_{i+1},\dots,\CA_n\ran{}^\perp \cap {}^\perp\lan\CA_1,\dots,\CA_{i-1}\ran$.

\begin{definition}[\cite{BK,B}]
A full triangulated subcategory $\CA$ of a triangulated category $\CT$ is called
{\sf right admissible}\/ if for the inclusion functor $i:\CA \to \CT$ there is
a right adjoint $i^!:\CT \to \CA$, and
{\sf left admissible}\/ if there is a left adjoint $i^*:\CT \to \CA$.
Subcategory $\CA$ is called {\sf admissible}\/ if it is both right and left admissible.
\end{definition}




\begin{lemma}[\cite{B}]\label{sod_adm}
If $\CT = \lan\CA,\CB\ran$ is a semiorthogonal decomposition then
$\CA$ is left amissible and $\CB$ is right admissible.
\end{lemma}

\begin{lemma}[\cite{B}]\label{sos_sod}
If\/ $\CA_1,\dots,\CA_n$ is a semiorthogonal sequence in $\CT$
such that $\CA_1,\dots,\CA_k$ are left admissible and
$\CA_{k+1},\dots,\CA_n$ are right admissible then
$$
\lan\CA_1,\dots,\CA_k,
{}^\perp\lan\CA_1,\dots,\CA_k\ran \cap \lan\CA_{k+1},\dots,\CA_n\ran{}^\perp,
\CA_{k+1},\dots,\CA_n\ran
$$
is a semiorthogonal decomposition.
\end{lemma}

Assume that $\CA\subset\CT$ is an admissible subcategory. Then
$\CT = \lan\CA,{}^\perp\CA\ran$ and $\CT = \lan\CA^\perp,\CA\ran$ are
semiorthogonal decompositions, hence ${}^\perp\CA$ is right admissible
and $\CA^\perp$ is left admissible.
Let $i_{{}^\perp\CA}:{}^\perp\CA \to \CT$
and $i_{\CA{}^\perp}:\CA{}^\perp \to \CT$
be the inclusion functors.

\begin{definition}\cite{B}
The functor $R_\CA = i_{{}^\perp\CA}i_{{}^\perp\CA}^!$ is called
the {\sf right mutation through $\CA$}.
The functor $L_\CA = i_{\CA{}^\perp}i_{\CA{}^\perp}^*$ is called
the {\sf left mutation through $\CA$}.
\end{definition}

\begin{lemma}\cite{B}\label{rl_mut}
We have $R_\CA(\CA) = 0$ and the restriction of $R_\CA$ to $\CA^\perp$
is an equivalence $\CA^\perp \to {}^\perp\CA$.
Similarly, we have $L_\CA(\CA) = 0$ and the restriction of $L_\CA$ to ${}^\perp\CA$
is an equivalence ${}^\perp\CA \to \CA^\perp$.
\end{lemma}

\begin{lemma}\cite{B}
If $\CA_1,\dots,\CA_n$ is a semiorthogonal sequence of admissible subcategories in $\CT$
then $R_{\langle\CA_1,\dots,\CA_n\rangle} = R_{\CA_n}\circ\dots\circ R_{\CA_1}$
and $L_{\langle\CA_1,\dots,\CA_n\rangle} = L_{\CA_1}\circ\dots\circ L_{\CA_n}$.
\end{lemma}

\begin{lemma}[\cite{O2}]\label{sod_proj}
If $E$ is a vector bundle of rank $r$ on $S$,
$\PP_S(E)$ is its projectivization,
$\CO(1)$ is the corresponding Grothendieck ample line bundle,
and $p:\PP_S(E) \to S$ is the projection then the pullback
$p^*:\D^b(S) \to \D^b(\PP_S(E))$ is fully faithful and
$$
\D^b(\PP_S(E)) = \lan
p^*(\D^b(S))\otimes\CO(k),
p^*(\D^b(S))\otimes\CO(k+1),\dots,
p^*(\D^b(S))\otimes\CO(k+r-1)
\ran
$$
is a semiorthogonal decomposition for any $k\in\ZZ$.
\end{lemma}

\subsection{Saturatedness and Serre functors}

\begin{definition}[\cite{B}]
A triangulated category $\CT$ is called
{\sf left saturated}\/ if every exact covariant functor
$\CT \to \D^b(\kk)$ is representable, and
{\sf right saturated}\/ if every exact contravariant functor
$\CT \to \D^b(\kk)$ is representable.
A triangulated category $\CT$ is called {\sf saturated}\/ if it is both left
and right saturated.
\end{definition}

\begin{lemma}[\cite{B}]\label{adm_sat}
A left {\rm(}resp.\ right{\rm)} admissible subcategory of a saturated category is
saturated.
\end{lemma}
\begin{proof}
Assume that $\CA$ is a left admissible subcategory in a saturated triangulated
category~$\CT$, $i:\CA \to \CT$ is the inclusion functor and
$i^*:\CT \to \CA$ is its left adjoint functor.
Let $\phi:\CA \to \D^b(\kk)$ be an exact covariant functor.
Then $\phi\circ i^*:\CT \to \D^b(\kk)$ is representable since $\CT$
is saturated. Therefore there exists $T\in\CT$ such that
$\phi\circ i^* \cong \Hom_\CT(T,-)$. Then
$$
\phi \cong \phi\circ i^*\circ i \cong \Hom_\CT(T,i(-)) \cong \Hom_\CA(i^*T,-),
$$
therefore $i^*T$ represents $\phi$.

Let $\psi:\CA \to \D^b(\kk)$ be an exact contravariant functor.
Then $\psi\circ i^*:\CT \to \D^b(\kk)$ is representable since $\CT$
is saturated. Therefore there exists $T\in\CT$ such that
$\psi\circ i^* \cong \Hom_\CT(-,T)$.
Note that $i^*({}^\perp\CA) = 0$, hence
$\Hom_\CT({}^\perp\CA,T) = 0$ which means that
$T \in ({}^\perp\CA)^\perp = \CA$, thus $T \cong i(A)$ with $A\in\CA$.
Finally
$$
\psi \cong \psi\circ i^*\circ i \cong \Hom_\CT(i(-),i(A)) \cong \Hom_\CA(-,A),
$$
since $i$ is fully faithful, therefore $A$ represents $\psi$.

A similar argument works for right admissible subcategories.
\end{proof}

\begin{lemma}[\cite{B}]\label{sat_adm}
If $\CA$ is saturated then $\CA$ is admissible.
\end{lemma}
\begin{proof}
For any object $T\in\CT$ consider the functor
$\Hom_\CT(T,i(-)):\CA \to \D^b(\kk)$. Since $\CA$ is saturated there exists
$A_T \in \CA$, such that this functor is isomorphic to $\Hom_\CA(A_T,-)$.
Since $\Hom_\CT(T,i(A_T)) \cong \Hom_\CA(A_T,A_T)$ we have a canonical morphism
$T \to i(A_T)$ and since $i$ is fully faithful it is easy to see that its cone
is contained in ${}^\perp\CA$. It follows that any morphism $T\to S$
composed with $S \to i(A_S)$ factors in a unique way as
$T \to i(A_T) \to i(A_S)$. Since
$\Hom_\CT(i(A_T),i(A_S)) \cong \Hom_\CA(A_T,A_S)$ the correspondence
$T\mapsto A_T$ is a functor $\CT \to \CA$, left adjoint to $i:\CA \to \CT$.
Similarly one can construct a right adjoint functor.
\end{proof}

\begin{lemma}[\cite{BV}]
If $X$ is a smooth projective variety then $\D^b(X)$ is saturated.
\end{lemma}

\begin{corollary}
If $X$ is a smooth projective variety and $\CA$ is a left {\rm(}resp.\ right{\rm)}
admissible subcategory in $\D^b(X)$ then $\CA$ is saturated.
\end{corollary}

\begin{definition}[\cite{BK},\cite{BO4}]
Let $\CT$ be a triangulated category.
A covariant additive functor $\FS:\CT\to\CT$ is a {\sf Serre functor}\/
if it is a category equivalence and for all objects $F,G\in\CT$ there are
given bi-functorial isomorphisms $\Hom(F,G)\to\Hom(G,\FS(F))^*$.
\end{definition}

\begin{lemma}[\cite{BK}]\label{serre_bk}
If a Serre functor exists then it is unique up to a canonical
functorial isomorphism. If $X$ is a smooth projective variety
then $\FS(F):=F\otimes\omega_X[\dim X]$ is a Serre functor
in $\D^b(X)$.
\end{lemma}

\begin{definition}\cite{BV}
A triangulated category $\CT$ is called $\Ext$-finite if for any objects
$F,G\in\CT$ the vector space $\oplus_{n\in\ZZ} \Hom_\CT(F,G[n])$ is
finite dimensional.
\end{definition}

\begin{lemma}[\cite{BK}]
If $\CT$ is an $\Ext$-finite saturated category then $\CT$ admits a Serre functor.
\end{lemma}

\begin{lemma}[\cite{BK}]\label{ser_loro}
If $\FS$ is a Serre functor for $\CT$ and $\CA$ is a subcategory of $\CT$
then $\FS({}^\perp\CA) = \CA^\perp$ and $\FS^{-1}(\CA^\perp) = {}^\perp\CA$.
In particular, if $\CT = \lan\CA_1,\CA_2\ran$ is a semiorthogonal decomposition
then $\CT = \lan\FS(\CA_2),\CA_1\ran$ and $\CT = \lan\CA_2,\FS^{-1}(\CA_1)\ran$
are semiorthogonal decompositions.
\end{lemma}

\begin{lemma}[\cite{B}]\label{mut_funct}
If $\CT$ admits a Serre functor $\FS$ and $\CA\subset\CT$ is right admissible
then $\CA$ admits a Serre functor $\FS_\CA = i^!\circ\FS\circ i$,
where $i:\CA \to \CT$ is the inclusion functor.
\end{lemma}
\begin{proof}
If $A,A'\in\CA$ then
$\Hom_\CA(A,i^!\FS i A') \cong
\Hom_\CT(i A,\FS i A') \cong
\Hom_\CT(i A',i A)^* \cong
\Hom_\CA(A',A)^*$.
\end{proof}

\subsection{$\Tor$ and $\Ext$-amplitude}

Let $f:X \to Y$ be a morphism of algebraic varieties.
For any subset $I \subset\ZZ$ we denote by $\D^I(X)$ the full subcategory
of $\D(X)$ consisting of all objects $F\in\D(X)$ with $\CH^k(F) = 0$
for $k\not\in I$.

\begin{definition}[\cite{K}]
An object $F \in \D(X)$ has {\sf finite $\Tor$-amplitude over $Y$}
{\rm(}resp.\ {\sf finite $\Ext$-amplitude over $Y$}{\rm)},
if there exist integers $p,q$ such that for any object $G \in \D^{[s,t]}(Y)$
we have $F\otimes f^*G \in \D^{[p+s,q+t]}(X)$
{\rm(}resp.\ $\RCHom(F,f^!G) \in \D^{[p+s,q+t]}(X)${\rm)}.
Morphism $f$ has {\sf finite $\Tor$-dimension},
{\rm(}resp.\ {\sf finite $\Ext$-dimension}{\rm)},
if the sheaf $\CO_X$ has finite $\Tor$-amplitude over $Y$
{\rm(}resp.\ finite $\Ext$-amplitude over $Y${\rm)}.
\end{definition}

The full subcategory of $\D(X)$ consisting of objects of finite $\Tor$-amplitude
(resp.\ of finite $\Ext$-amplitude) over $Y$ is denoted by $\D_{fTd/Y}(X)$
(resp.\ $\D_{fEd/Y}(X)$).
Both are triangulated subcategories of $\D^b(X)$.

\begin{lemma}[\cite{K}]\label{isfted}
If $i:X \to X'$ is a finite morphism over $Y$ then
$F\in\D_{fTd/Y}(X)$ $\LRA$ $i_*F\in\D_{fTd/Y}(X')$ and
$F\in\D_{fEd/Y}(X)$ $\LRA$ $i_*F\in\D_{fEd/Y}(X')$.
\end{lemma}

\begin{lemma}[\cite{K}]\label{ftedperf}
If morphism $f:X \to Y$ has finite $\Tor$-dimension
{\rm(}resp.\ $\Ext$-dimension{\rm)} then any perfect complex on $X$
has finite $\Tor$-amplitude {\rm(}resp.\ $\Ext$-amplitude{\rm)} over $Y$.
\end{lemma}

\begin{lemma}[\cite{K}]\label{snav_sm1}
If $f:X \to Y$ is a smooth morphism then $\D_{fTd/Y}(X) = \D^\perf(X) = \D_{fEd/Y}(X)$.
\end{lemma}

\begin{definition}[\cite{K}]
A triangulated category $\CT$ is {\sf $\Ext$-bounded}, if
for any objects $F,G\in\CT$ the set $\{n\in\ZZ\ |\ \Hom(F,G[n])\ne 0\}$
is finite.
\end{definition}

\begin{lemma}[\cite{K}]\label{snav_sm}
The following conditions for an algebraic variety $X$ are equivalent:

\noindent$(i)$
$X$ is smooth;

\noindent$(ii)$
$\D^b(X) = \D^\perf(X)$.

\noindent$(iii)$
the bounded derived category $\D^b(X)$ is $\Ext$-bounded.
\end{lemma}

\begin{lemma}\label{abt}
Assume that $\CT = \lan\CA,\CB\ran$ is a semiorthogonal decomposition.
If both $\CA$ and $\CB$ are $\Ext$-bounded and either $\CA$ or $\CB$
is admissible then $\CT$ is $\Ext$-bounded.
\end{lemma}
\begin{proof}
Let $F,G\in\CT$. Then there exist exact triangles
$$
\beta\beta^! F \to F \to \alpha\alpha^*F,\qquad
\beta\beta^! G \to G \to \alpha\alpha^*G.
$$
Computing $\Hom(F,G[n])$ and using semiorthogonality of $\CA$ and $\CB$
we obtain a long exact sequence
$$
\dots \to
\Hom(\alpha\alpha^*F,\beta\beta^!G[n]) \to
\Hom(F,G[n]) \to
\Hom(\beta^!F,\beta^!G[n]) \oplus \Hom(\alpha^*F,\alpha^*G[n]) \to
\dots
$$
Since $\CA$ and $\CB$ are $\Ext$-bounded, the third term vanishes for $|n|\gg 0$.
On the other hand, if $\CA$ is admissible then
$\Hom(\alpha\alpha^*F,\beta\beta^!G[n]) \cong \Hom(\alpha^*F,\alpha^!\beta\beta^!G[n])$,
hence the first term also vanishes for $|n|\gg 0$.
Similarly, if $\CB$ is admissible then
$\Hom(\alpha\alpha^*F,\beta\beta^!G[n]) \cong \Hom(\beta^*\alpha\alpha^*F,\beta^!G[n])$,
hence the first term also vanishes for $|n|\gg 0$.
In both cases we deduce that $\Hom(F,G[n])$ vanishes for $|n|\gg0$, hence
$\CT$ is $\Ext$-bounded.
\end{proof}

\subsection{Kernel functors}

Let $X_1$, $X_2$ be algebraic varieties and
let $p_i:X_1\times X_2 \to X_i$ denote the projections.
Take any $K\in\D_{qc}^-(X_1\times X_2)$ and define functors
$$
\Phi_K(F_1) := {p_2}_*(p_{1\circ}^*F_1\otimes K),\qquad
\Phi_K^!(F_2) := {p_1}_*\RCHom(K,p_{2\circ}^!F_2).
$$
Then $\Phi_K$ is an exact functor $\D_{qc}^-(X_1)\to\D_{qc}^-(X_2)$ and
$\Phi_K^!$ is an exact functor $\D_{qc}^+(X_2)\to\D_{qc}^+(X_1)$.
We call $\Phi_K$ the {\sf kernel functor}\/ with kernel $K$,
and $\Phi_K^!$ the {\sf kernel functor of the second type}\/ with kernel~$K$ (cf.\ \cite{K}).
In smooth case any kernel functor of the second type is isomorphic to
a usual kernel functor:
$\Phi_K^! \cong \Phi_{\RCHom(K,\omega_{X_1}[\dim X_1])}$.

\begin{lemma}[\cite{K}]\label{phi_bounded}
$(i)$
If $K$ has coherent cohomologies, finite $\Tor$-amplitude over $X_1$ and
$\supp(K)$ is projective over $X_2$ then $\Phi_K$ takes
$\D^b(X_1)$ to $\D^b(X_2)$.

\noindent $(ii)$
If $K$ has coherent cohomologies, finite $\Ext$-amplitude over $X_2$ and
$\supp(K)$ is projective over $X_1$ then $\Phi_K^!$
takes $\D^b(X_2)$ to $\D^b(X_1)$.

\noindent $(iii)$
If both $(i)$ and $(ii)$ hold then $\Phi_K^!$ is right adjoint to~$\Phi_K$.
Moreover, $\Phi_K$ takes $\D^\perf(X_1)$ to $\D^\perf(X_2)$.
\end{lemma}

\begin{lemma}[\cite{K}]\label{ladj}
If $K$ is a perfect complex, $X_2$ is smooth and $\supp(K)$
is projective both over $X_1$ and over $X_2$, then the functor $\Phi_K$
admits a left adjoint functor $\Phi_K^*$ which is isomorphic
to a kernel functor $\Phi_{K^\#}$ with the kernel
$$
K^\#:= \RCHom(K,\omega_{X_2}[\dim X_2]).
$$
\end{lemma}

Consider kernels
$K_{12} \in \D^-(X_1\times X_2)$,
$K_{23} \in \D^-(X_2\times X_3)$.
Denote by $p_{ij}:X_1\times X_2\times X_3 \to X_i\times X_j$
the projections.
We define the {\sf convolution}\/ of kernels as follows
$$
K_{23}\circ K_{12} :=
{p_{13}}_*(p_{12\circ}^*K_{12}\otimes p_{23\circ}^*K_{23}),
$$

\begin{lemma}\label{kerconv}
For $K_{12} \in \D^-(X_1\times X_2)$,
$K_{23} \in \D^-(X_2\times X_3)$
we have
$\Phi_{K_{23}} \circ \Phi_{K_{12}} = \Phi_{K_{23}\circ K_{12}}$.
\end{lemma}

Assume that $\Phi_1,\Phi_2,\Phi_3:\D \to \D'$ are exact functors
between triangulated categories, and $\alpha:\Phi_1\to\Phi_2$,
$\beta:\Phi_2\to\Phi_3$, $\gamma:\Phi_3\to\Phi_1[1]$ are morphisms
of functors. We say that
$$
\Phi_1 \exto{\alpha} \Phi_2 \exto{\beta} \Phi_3 \exto{\gamma} \Phi_1[1]
$$
is an {\sf exact triangle of functors}, if for any object $F\in\D$
the triangle
$$
\Phi_1(F) \exto{\alpha(F)} \Phi_2(F) \exto{\beta(F)}
\Phi_3(F) \exto{\gamma(F)} \Phi_1(F)[1]
$$
is exact in $\D'$.

\begin{lemma}\label{etf}
If  $K_1 \exto{\alpha} K_2 \exto{\beta} K_3 \exto{\gamma} K_1[1]$ is
an exact triangle in $\D^-(X\times Y)$ then we have the following
exact triangles of functors
$$
\Phi_{K_1} \exto{\alpha_*} \Phi_{K_2} \exto{\beta_*}
\Phi_{K_3} \exto{\gamma_*} \Phi_{K_1}[1]
$$
$$
\Phi^!_{K_3} \exto{\beta^!} \Phi^!_{K_2} \exto{\alpha^!}
\Phi^!_{K_1} \exto{\gamma^!} \Phi^!_{K_3}[1]
$$
If additionally kernels $K_1$, $K_2$ and $K_3$ satisfy the conditions
of lemma~$\ref{ladj}$ then we have also the following exact triangle of functors
$$
\Phi^*_{K_3} \exto{\beta^*} \Phi^*_{K_2} \exto{\alpha^*}
\Phi^*_{K_1} \exto{\gamma^*} \Phi^*_{K_3}[1]
$$
\end{lemma}
\begin{proof}
Evident.
\end{proof}

\subsection{Exact cartesian squares}

Consider a cartesian square
$$
\vcenter{\xymatrix{
X\times_S Y  \ar[r]^-q \ar[d]_p         &       Y \ar[d]_g \\
X \ar[r]^f                              &       S
}}
$$
Consider the functors $q_*p^*$ and $g^*f_*:\D^b(X) \to \D^b(Y)$.
It is easy to see that both are kernel functors. Explicitly, the first
is given by the structure sheaf of the fiber product $\CO_{X\times_S Y}$
and the second is given by the convolution of the structure sheaves
of graphs of $f$ and $g$ respectively. It is easy to see that the latter
kernel is a complex supported on the fiber product, the top cohomology
of which is isomorphic to $\CO_{X\times_S Y}$.
The natural map from this complex to its top cohomology induces a morphism
of functors $g^*f_* \to q_*p^*$.
A cartesian square is called {\sf exact cartesian}\/ \cite{K} if this morphism
of functors is an isomorphism. As explained above a square is exact cartesian
if and only if the convolution of the structure sheaves of graphs of $f$ and $g$
is isomorphic to its top cohomology.

\begin{lemma}[\cite{K}]\label{ec}
Consider a cartesian square as above.

\noindent$(i)$
If either $f$ or $g$ is flat then the square is exact cartesian.

\noindent$(ii)$
A square is exact cartesian, if and only if the transposed square is exact cartesian.

\noindent$(iii)$
If $g$ is a closed embedding, $Y \subset S$ is a locally complete intersection,
both $S$ and $X$ are Cohen--Macaulay, and $\codim_X (X\times_S Y) = \codim_S Y$,
then the square is exact cartesian.
\end{lemma}


\subsection{Derived categories over a base}

Consider a pair of algebraic varieties $X$ and $Y$
over the same smooth algebraic variety $S$. In other words, we have
a pair of morphisms $f:X \to S$ and $g:Y\to S$.

A functor $\Phi:\D(X) \to \D(Y)$ is called {\sf $S$-linear} \cite{K}
if for all $F\in\D(X)$, $G\in\D^b(S)$ there are given bifunctorial isomorphisms
$$
\Phi(f^*G\otimes F) \cong g^*G\otimes \Phi(F).
$$
Note that since $S$ is smooth any object $G\in\D^b(S)$ is a perfect complex.

\begin{lemma}[\cite{K}]\label{adjslin}
If $\Phi$ is $S$-linear and admits a right adjoint functor $\Phi^!$
then $\Phi^!$ is also $S$-linear.
If $K\in\D^-(X\times_S Y)$ then the kernel functors
$\Phi_{i_*K}$ and $\Phi^!_{i_*K}$ are $S$-linear.
\end{lemma}

A strictly full subcategory $\CC \subset \D(X)$
is called {\sf $S$-linear} if for all $F\in\CC$,
$G\in\D^b(S)$ we have $f^*G\otimes F\in \CC$.

\begin{lemma}[\cite{K}]\label{slinperp}
If $\CC\subset\D^b(X)$ is a strictly full $S$-linear left
{\rm(}resp.\ right{\rm)} admissible triangulated subcategory
then its left {\rm(}resp.\ right{\rm)} orthogonal is also $S$-linear.
\end{lemma}

\subsection{Faithful base change theorem}

Consider morphisms $f:X \to S$ and $g:Y \to S$ with smooth $S$.
For any base change $\phi:T \to S$ we consider the fiber products
$$
X_T := X\times_S T,\qquad
Y_T := Y\times_S T,\qquad
X_T\times_T Y_T = (X\times_S Y)\times_S T
$$
and denote the projections $X_T \to X$, $Y_T \to Y$, and
$X_T\times_T Y_T \to X\times_S Y$ also by $\phi$.
For any kernel $K \in \D^-(X\times_S Y)$ we denote
$K_T = \phi^* K \in \D^-(X_T\times_T Y_T)$.

\begin{definition}[\cite{K}]
A change of base $\phi:T \to S$ is called {\sf faithful} with respect
to a morphism $f:X \to S$ if the cartesian square
$$
\xymatrix@=20pt{X_T \ar[r] \ar[d] &  X \ar[d]^f \\ T \ar[r]^\phi & S}
$$
is exact cartesian.
A change of base $\phi:T \to S$ is called {\sf faithful} for a pair $(X,Y)$
if $\phi$ is faithful with respect to morphisms $f:X\to S$, $g:Y\to S$,
and $f\times_{\scriptscriptstyle S} g:X\times_S Y \to S$.
\end{definition}

Using the criterions of lemma~\ref{ec} it is easy to deduce the following

\begin{lemma}[\cite{K}]\label{fisf}
Let $f:X \to S$ be a morphism and $\phi:T \to S$ a base change.

\noindent$(i)$
If\/ $\phi$ is flat then it is faithful.

\noindent$(ii)$
If\/ $T$ and $X$ are smooth and $\dim X_T = \dim X + \dim T - \dim S$ then
the base change $\phi:T \to S$ is faithful with respect to the morphism $f:X \to S$.
\end{lemma}

\begin{lemma}[\cite{K}]\label{fted_bc}
If $\phi:T \to S$ is a faithful base change for a morphism $f:X \to S$
then we have $\phi^*(\D_{fTd/S}(X)) \subset \D_{fTd/T}(X_T)$,
and $\phi^*(\D_{fEd/S}(X)) \subset \D_{fEd/T}(X_T)$.
\end{lemma}

\begin{lemma}[\cite{K}]\label{phit}
If $\phi:T \to S$ is a base change faithful for a pair $(X,Y)$,
and $f$ is projective then we have
$\Phi_{K_T}\phi^*   = \phi^*\Phi_{K},\quad
\Phi_{K}\phi_*     = \phi_*\Phi_{K_T},\quad
\Phi_{K_T}^!\phi^* = \phi^*\Phi_{K}^!,\quad\text{and}\quad
\Phi_{K}^!\phi_*   = \phi_*\Phi_{K_T}^!$.
\end{lemma}

\begin{proposition}[\cite{K}]\label{fbc_sod}
If $\phi$ is faithful for a pair $(X,Y)$, varieties $X$ and $Y$ are
projective over $S$ and smooth, and $K \in \D^b(X\times_S Y)$
is a kernel such that
$\Phi_K:\D^b(X) \to \D^b(Y)$ is fully faithful then
$\Phi_{K_T}:\D^b(X_T) \to \D^b(Y_T)$ is fully faithful.
\end{proposition}

\begin{theorem}[\cite{K}]\label{phitsod}
If $\D^b(Y) = \langle\Phi_{K_1}(\D^b(X_1)),\dots,\Phi_{K_n}(\D^b(X_n))\rangle$
is a semiorthogonal decomposition, with $K_i \in \D^b(X_i\times_S Y)$,
a base change $\phi$ is faithful for all pairs $(X_1,Y)$, \dots, $(X_n,Y)$,
and all varieties $X_1,\dots,X_n,Y$ are projective over $S$ and smooth
then
$\D^b(Y_T) = \langle\Phi_{K_{1T}}(\D^b(X_{1T})),\dots,\Phi_{K_{nT}}(\D^b(X_{nT}))\rangle$
is a semiorthogonal decomposition.
\end{theorem}

Note that though $X_1,\dots,X_n,Y$ are smooth in the assumptions
of the theorem, their pullbacks $X_{1T},\dots,X_{nT},Y_T$
under the base change $\phi$ are singular in general.

We will need also the following theorem.

\begin{theorem}\label{loc_sod}
If $S$ and $Y$ are smooth and for any point $s\in S$ there exists an open neighborhood $U \subset S$ such that
$\D^b(Y_U) = \langle\Phi_{K_{1U}}(\D^b(X_{1U})),\dots,\Phi_{K_{nU}}(\D^b(X_{nU}))\rangle$
is a semiorthogonal decomposition then also
$\D^b(Y) = \langle\Phi_{K_1}(\D^b(X_1)),\dots,\Phi_{K_n}(\D^b(X_n))\rangle$
is a semiorthogonal decomposition.
\end{theorem}
\begin{proof}
We must check that for every $i = 1, \dots, n$ the functor
$\Phi_{K_i}:\D^b(X_i) \to \D^b(Y)$ is fully faithful.
Equivalently, we must show that the morphism of functors
$\id_{\D^b(X_i)} \to \Phi_{K_i}^!\Phi_{K_i}$ is an isomorphism.
Note that $\D^b(X_{iU})$ being a semiorthogonal summand of an $\Ext$-bounded
category $\D^b(Y_U)$ is $\Ext$-bounded, hence $X_{iU}$ is smooth,
hence $X_i$ is smooth for any $i$. Therefore the functors $\Phi_{K_i}^!$
are kernel functors. Note also that the morphism of functors
$\id_{\D^b(X_i)} \to \Phi_{K_i}^!\Phi_{K_i}$ is induced by morphism
of kernels. Moreover, restricting this morphism of kernels
from $S$ to $U$ we obtain precisely the morphism of kernels
corresponding to the canonical morphism of functors
$\id_{\D^b(X_{iU})} \to \Phi_{K_{iU}}^!\Phi_{K_{iU}}$.
Since the latter morphism is an isomorphism by assumptions
for suitable $U$, it follows that the corresponding morphism
of kernels is an isomorphism over $U$. Since this is true
for a suitable neighborhood of every point $s \in S$,
we deduce that the morphism of kernels is an isomorphism
over the whole $S$, hence $\Phi_{K_i}$ is fully faithful.

Further, we must check the semiorthogonality.
Equivalently, we must show that the functor
$\Phi_{K_j}^!\Phi_{K_i}$ is zero for all $1\le i < j \le n$.
As above we note that this functor is a kernel functor.
Restricting its kernel from $S$ to $U$ we obtain precisely the kernel
of the functor $\Phi_{K_{jU}}^!\Phi_{K_{iU}}$. Since the latter functor
is zero by assumptions for suitable $U$, it follows that the corresponding
kernel is zero over $U$. Since this is true for a suitable neighborhood
of every point $s \in S$, we deduce that the kernel is zero over the whole~$S$,
hence $\Phi_{K_j}^!\Phi_{K_i} = 0$.

Finally, we must check that our semiorthogonal collection generates $\D^b(Y)$.
Assume that there is an object in the orthogonal to
$\langle\Phi_{K_1}(\D^b(X_1)),\dots,\Phi_{K_n}(\D^b(X_n))\rangle$.
Then it is easy to see that its restriction from $S$ to $U$ is in the orthogonal to
$\langle\Phi_{K_{1U}}(\D^b(X_{1U})),\dots,\Phi_{K_{nU}}(\D^b(X_{nU}))\rangle$.
By assumptions we deduce that this object is zero over $U$. Since this is true
for a suitable neighborhood of every point $s \in S$, we deduce that the object
is zero over the whole $S$.
\end{proof}

\section{Splitting functors}

Assume that $\CA$ and $\CB$ are triangulated categories and
$\Phi:\CB \to \CA$ is an exact functor.
Consider the following full subcategories of $\CA$ and $\CB$:
$$
\Ker\Phi = \{B\in\CB\  |\ \Phi(B) = 0 \} \subset \CB,
\qquad
\Im\Phi  = \{A \cong \Phi(B)\ |\ B\in\CB\} \subset \CA.
$$
Note that $\Ker\Phi$ is a triangulated subcategory of $\CB$,
and if $\Phi$ is fully faithful then $\Im\Phi$ is also triangulated.
However, if $\Phi$ is not fully faithful, in general $\Im\Phi$ is not
triangulated. If $\Phi$ admits an adjoint functor then we have
$$
\begin{array}{ll}
\Hom(\Ker\Phi,\Im\Phi^!) = 0,
\qquad & \text{if $\Phi$ admits a right adjoint $\Phi^!$,}\\
\Hom(\Im\Phi^*,\Ker\Phi) = 0,
\qquad & \text{if $\Phi$ admits a left adjoint $\Phi^*$,}
\end{array}
$$
(evidently follows from the adjunction).

\begin{definition}
An exact functor $\Phi:\CB \to \CA$ is called {\sf right splitting}\/ if
$\Ker\Phi$ is a right admissible subcategory in $\CB$,
the restriction of $\Phi$ to $(\Ker\Phi)^\perp$ is fully faithful,
and $\Im\Phi$ is right admissible in $\CA$ (note that
$\Im\Phi = \Im(\Phi_{(\Ker\Phi)^\perp})$
is a triangulated subcategory of $\CA$).
An exact functor $\Phi:\CB \to \CA$ is called {\sf left splitting}\/ if
$\Ker\Phi$ is a left admissible subcategory in $\CB$,
the restriction of $\Phi$ to ${}^\perp(\Ker\Phi)$ is fully faithful,
and $\Im\Phi$ is left admissible in $\CA$.
\end{definition}

\begin{lemma}\label{adjspl}
A right {\rm(}resp.\ left{\rm)} splitting functor has
a right {\rm(}resp.\ left{\rm)} adjoint functor.
\end{lemma}
\begin{proof}
If $\Ker\Phi$ is right admissible then $(\Ker\Phi)^\perp$
is left admissible and we have a semiorthogonal decomposition
$\CB = \lan(\Ker\Phi)^\perp,\Ker\Phi\ran$ by lemmas~\ref{sos_sod}
and \ref{sod_adm}. Since $\Phi$ vanishes
on the second term and is fully faithful on the first term it follows that
$\Phi \cong j\circ\phi\circ i^*$, where $i:(\Ker\Phi)^\perp \to \CB$
and $j:\Im\Phi \to \CA$ are the inclusion functors, $i^*$ is a left
adjoint to $i$, and $\phi:(\Ker\Phi)^\perp \to \Im\Phi$ is an equivalence
of categories induced by $\Phi$. Therefore $\Phi^! := i\circ\phi^{-1}\circ j^!$
is right adjoint to $\Phi$ (functor $j^!$ right adjoint to $j$ exists
because $\Im\Phi$ is right admissible).
\end{proof}


\begin{theorem}\label{sf_th}
Let $\Phi:\CB \to \CA$ be an exact functor.
Then the following conditions are equivalent
$(1r) \LRA (2r) \LRA (3r) \LRA (4r)$ and $(1l) \LRA (2l) \LRA (3l) \LRA (4l)$,
where

\noindent$(1r)$ $\Phi$ is right splitting;

\noindent$(2r)$ $\Phi$ has a right adjoint functor $\Phi^!$
and the composition of the canonical morphism of functors
$\id_\CB \to \Phi^!\Phi$ with $\Phi$ gives an isomorphism
$\Phi \cong \Phi\Phi^!\Phi$;

\noindent$(3r)$ $\Phi$ has a right adjoint functor $\Phi^!$,
there are semiorthogonal decompositions
$$
\CB = \langle\Im\Phi^!,\Ker\Phi\rangle,
\qquad
\CA = \langle\Ker\Phi^!,\Im\Phi\rangle,
$$
and the functors $\Phi$ and $\Phi^!$ give quasiinverse equivalences
$\Im\Phi^! \cong \Im\Phi$;

\noindent$(4r)$ there exists a triangulated category $\CC$ and fully faithful
functors $\alpha:\CC \to \CA$, $\beta:\CC \to \CB$, such that $\alpha$ admits
a right adjoint, $\beta$ admits a left adjoint and $\Phi \cong \alpha\circ\beta^*$.

\noindent$(1l)$ $\Phi$ is left splitting;

\noindent$(2l)$ $\Phi$ has a left adjoint functor $\Phi^*$
and the composition of the canonical morphism of functors
$\Phi^*\Phi \to \id_\CB$ with $\Phi$ gives an isomorphism
$\Phi\Phi^*\Phi \cong \Phi$;

\noindent$(3l)$ $\Phi$ has a left adjoint functor $\Phi^*$,
there are semiorthogonal decompositions
$$
\CB = \langle\Ker\Phi,\Im\Phi^*\rangle,
\qquad
\CA = \langle\Im\Phi,\Ker\Phi^*\rangle,
$$
and the functors $\Phi$ and $\Phi^*$ give quasiinverse equivalences
$\Im\Phi^* \cong \Im\Phi$;

\noindent$(4l)$ there exists a triangulated category $\CC$ and fully faithful
functors $\alpha:\CC \to \CA$, $\beta:\CC \to \CB$, such that $\alpha$ admits
a left adjoint, $\beta$ admits a right adjoint and $\Phi \cong \alpha\circ\beta^!$.
%
%
%
%
%
%
%
%
%
%
\end{theorem}
\begin{proof}
$(1r) \RA (2r)$:
using the formula of lemma~\ref{adjspl} for $\Phi^!$
we deduce that $\Phi^!\Phi \cong i\phi^{-1}j^!j\phi i^* \cong i i^*$.
Composing with $\Phi$ we obtain
$\Phi\Phi^!\Phi \cong j\phi i^*i i^* \cong j \phi i^* \cong \Phi$.

$(2r) \RA (3r)$:
for any $B\in \CB$ let $K_B$ be the object defined from
the triangle
\begin{equation}\label{kb}
K_B \to B \to \Phi^!\Phi B.
\end{equation}
Applying the functor $\Phi$ to this triangle and using the assumption
we deduce that $\Phi(K_B)=0$, i.e.\ $K_B\in\Ker\Phi$. Thus any object $B$
can be included as the second vertex in a triangle with first vertex in
$\Ker\Phi$ and the third vertex in $\Im\Phi^!$. Since these categories are
semiorthogonal, we obtain the desired semiorthogonal decomposition for $\CB$.
Moreover, it follows from $(2r)$ that for $A\in\Im\Phi$ we have $A \cong \Phi\Phi^! A$,
hence we have an isomorphism of functors $\id \cong \Phi\Phi^!$ on $\Im\Phi$.
On the other hand, if $B\in\Im\Phi^!$ then $K_B = 0$ since $K_B$ is
the component of $B$ in $\Ker\Phi$ with respect to the semiorthogonal
decomposition $\CB = \lan\Im\Phi^!,\Ker\Phi\ran$. Therefore,
$\id \cong \Phi^!\Phi$ on $\Im\Phi^!$. Thus $\Phi$ and $\Phi^!$ are
quasiinverse equivalences between $\Im\Phi$ and $\Im\Phi^!$. Finally,
we note that for any $B\in\Im\Phi^!$, $A\in \CA$  we have
$$
\Hom_\CA(\Phi B,A) \cong \Hom_\CB(B,\Phi^!A) \cong \Hom_\CA(\Phi B,\Phi\Phi^! A)
$$
since $\Phi$ is fully faithful on $\Im\Phi^!$, hence $\Phi\Phi^!:\CA \to \Im\Phi$
is a right adjoint to the inclusion functor $\Im\Phi \to \CA$, hence $\Im\Phi$ is
right admissible, we have $\CA = \lan(\Im\Phi)^\perp,\Im\Phi\ran$ and it remains
to note that $(\Im\Phi)^\perp = \Ker\Phi^!$.

$(3r) \RA (4r)$: take $\CC = \Im\Phi$ with $\alpha$ being the inclusion functor
$\Im\Phi \to \CA$ and $\beta$ being the composition of the equivalence
$\Im\Phi \cong \Im\Phi^!$ and of the inclusion functor $\Im\Phi^! \to \CB$.
Then $\alpha$ admits a right adjoint because $\Im\Phi$ is right admissible in $\CA$
and $\beta$ admits a left adjoint because $\Im\Phi^!$ is left admissible in $\CB$
and we evidently have $\Phi\cong \alpha\circ\beta^*$.

$(4r) \RA (1r)$: $\Im\Phi = \alpha(\CC)$ is right admissible because $\alpha$
admits a right adjoint functor; on the other hand
$\Ker\Phi = \Ker(\beta^!) = {}^\perp\beta(\CC)$
is right admissible as the left orthogonal to $\beta(\CC)$ which is
left admissible because $\beta$ admits a left adjoint functor. Finally,
$\Phi = \alpha\circ\beta^*$ restricted to $(\Ker\Phi)^\perp = \beta(\CC)$ is
isomorphic to the composition of an equivalence $\beta(\CC) \cong \CC$ and
of a fully faithful functor $\alpha:\CC \to \CA$, hence fully faithful.

The equivalences  $(1l) \LRA (2l) \LRA (3l) \LRA (4l)$ are proved by similar arguments.
\end{proof}

\begin{corollary}
If $\Phi$ is a right {\rm(}resp.\ left{\rm)} splitting functor and $\Psi$
is its right {\rm(}resp.\ left{\rm)} adjoint then $\Psi$ is a left
{\rm(}resp.\ right{\rm)} splitting functor.
\end{corollary}
\begin{proof}
Compare $(3r)$ and $(3l)$ for $\Phi$ and $\Psi$.
\end{proof}

\begin{lemma}
If either $\CA$ or $\CB$ is saturated and $\Phi:\CB\to\CA$ is
a right {\rm(}resp.\ left{\rm)} splitting functor then
$\Phi$ is also a left {\rm(}resp.\ right{\rm)} splitting.
\end{lemma}
\begin{proof}
Assume that $\CB$ is saturated and $\Phi$ is right admissible.
Then $\Ker\Phi$ and $(\Ker\Phi)^\perp$ are saturated by lemma~\ref{adm_sat}.
Moreover, $\Im\Phi \cong (\Ker\Phi)^\perp$, hence $\Im\Phi$ is also saturated.
Hence by lemma~\ref{sat_adm} both $\Ker\Phi$ and $\Im\Phi$ are left admissible.
%
%
Moreover,
it is easy to see that the restriction of $\Phi$ to ${}^\perp(\Ker\Phi)$
is isomorphic to the composition of the restriction of $\Phi$ to
$(\Ker\Phi)^\perp$ with the mutation functor $L_{\Ker\Phi}$
$$
\xymatrix{
(\Ker\Phi)^\perp \ar[dr]_\Phi && {}^\perp(\Ker\Phi) \ar[dl]^\Phi \ar[ll]_{L_{\Ker\Phi}} \\
& \Im\Phi
}
$$
But the upper arrow $L_{\Ker\Phi}$ is fully faithful on ${}^\perp(\Ker\Phi)$ by lemma~\ref{rl_mut},
hence $\Phi$ is fully faithful on ${}^\perp(\Ker\Phi)$.
\end{proof}

We will also need an analog of the faithful base change theorem
for splitting functors.

\begin{proposition}\label{fbc_spl}
In the notations of proposition~$\ref{fbc_sod}$
if $\phi:T\to S$ is a faithful base change for a pair $(X,Y)$
over a smooth base scheme~$S$, $X$ and $Y$ are projective over $S$ and smooth,
and $K \in \D^b(X\times_S Y)$ is a kernel such that
the functor $\Phi_K:\D^b(X) \to \D^b(Y)$ is splitting then
the functor $\Phi_{K_T}:\D^b(X_T) \to \D^b(Y_T)$ is also splitting.
\end{proposition}
\begin{proof}
Analogous to the proof of Proposition~2.42 of \cite{K} using
criterion~(2r) or (2l) to check that the functors are splitting.
\end{proof}

The class of splitting functors is a good generalization of the class
of fully faithful functors having an adjoint. Recall that it was proved
by Orlov in \cite{O3} that any fully faithful functor having an adjoint
between derived categories of smooth projective varieties is isomorphic
to a kernel functor. It would be nice to prove the same result
for splitting functors.

\begin{conjecture}
A splitting functor between bounded derived categories of coherent sheaves
on smooth projective varieties is isomorphic to a kernel functor.
\end{conjecture}


\section{Lefschetz decompositions}

Assume that $X$ is an algebraic variety with a line bundle $\CO_X(1)$ on $X$.

\begin{definition}
A {\sf Lefschetz decomposition} of the derived category $\D^b(X)$ is
a semiorthogonal decomposition of $\D^b(X)$ of the form
\begin{equation}\label{dbx}
\D^b(X) = \langle \CA_0,\CA_1(1),\dots,\CA_{\ix-1}(\ix-1) \rangle,
\qquad
0 \subset \CA_{\ix-1} \subset \CA_{\ix-2} \subset \dots \subset \CA_1 \subset \CA_0 \subset \D^b(X),
\end{equation}
where $0 \subset \CA_{\ix-1} \subset \CA_{\ix-2} \subset \dots \subset \CA_1 \subset \CA_0 \subset \D^b(X)$
is a chain of admissible subcategories of $\D^b(X)$. Lefschetz decomposition
is called {\sf rectangular}\/ if $\CA_{\ix-1} = \CA_{\ix-2} = \dots = \CA_1 = \CA_0$.
\end{definition}

Let $\fa_k$ denote the right orthogonal to $\CA_{k+1}$ in $\CA_k$.
The categories $\fa_0,\fa_1,\dots,\fa_{\ix-1}$ will be called {\sf primitive}
categories of the Lefschetz decomposition~(\ref{dbx}).
By definition we have the following semiorthogonal decompositions:
\begin{equation}\label{sodak}
\CA_k = \langle\fa_k,\fa_{k+1},\dots,\fa_{\ix-1}\rangle.
\end{equation}
If the Lefschetz decomposition is rectangular then we have
$\fa_0 = \fa_1 = \dots = \fa_{\ix-2} = 0$ and $\fa_{\ix-1} = \CA_{\ix-1}$.

Assume that the bounded derived category of coherent
sheaves on $X$, $\D^b(X)$ admits a Lefschetz decomposition~(\ref{dbx})
with respect to $\CO_X(1)$. If $X$ is smooth and projective then
its derived category $\D^b(X)$ is saturated and admits a Serre functor.
Therefore for every $0\le k\le \ix-1$ the category $\CA_k$ is saturated
and has a Serre functor too. Moreover, for every $0\le k\le \ix-1$
the primitive category $\fa_k$ is also saturated and has a Serre functor.

Let $\alpha_k:\CA_k(k) \to \D^b(X)$ denote the embedding functor and let
$\alpha_k^*,\alpha^!:\D^b(X) \to \CA_k(k)$ be the left and the right
adjoint functors. Let $\FS_X$ denote a Serre functor of $\D^b(X)$,
$\FS_X(F) \cong F\otimes\omega_X[\dim X]$, and let $
\FS_0$ denote a Serre functor of $\CA_0$.

Consider the restriction of the functor $\alpha_0^*:\D^b(X) \to \CA_0$
to the subcategory $\CA_k(k+1) \subset \D^b(X)$. It follows from~(\ref{dbx})
that $\alpha_0^*(\CA_{k+1}(k+1))=0$, hence it factors through the quotient
$(\CA_k/\CA_{k+1})(k+1)$.

\begin{lemma}\label{a0sfak}
The functor $\alpha_0^*:(\CA_k/\CA_{k+1})(k+1) \to \CA_0$ is fully faithful.
\end{lemma}
\begin{proof}
It is clear that $\fa_k(k+1)$ is the right orthogonal to $\CA_{k+1}(k+1)$
in $\CA_k(k+1)$, hence we have to check that $\alpha_0^*$ is fully faithful
on $\fa_k(k+1)$. For this we note that
\begin{equation}\label{fakro}
\fa_k(k+1) \subset \langle\CA_{k+1}(k+1),\dots,\CA_{\ix-1}(\ix-1)\rangle^\perp =
\langle \CA_0,\CA_1(1),\dots,\CA_k(k)\rangle,
\end{equation}
since for $l>k+1$ we have
$\Hom(\CA_l(l),\fa_k(k+1)) = \Hom(\CA_l(l-1),\fa_k(k))$
and $\CA_l(l-1)\subset\CA_{l-1}(l-1)$, $\fa_k(k) \subset\CA_k(k)$, while
$\Hom(\CA_{k+1}(k+1),\fa_k(k+1)) = \Hom(\CA_{k+1},\fa_k) = 0$ by definition
of $\fa_k$. On the other hand, we have
\begin{equation}\label{faklo}
\fa_k(k+1) \subset {}^\perp\langle \CA_1(1),\dots,\CA_k(k)\rangle,
\end{equation}
since for $1\le l\le k$ we have
$\Hom(\fa_k(k+1),\CA_l(l)) = \Hom(\fa_k(k),\CA_l(l-1))$
and $\CA_l(l-1)\subset\CA_{l-1}(l-1)$, $\fa_k(k) \subset\CA_k(k)$.
It follows from (\ref{fakro}) that the functor $\alpha_0^*$
restricted to $\fa_k(k+1)$ is just the left mutation
of $\fa_k(k+1)$ through $\langle \CA_1(1),\dots,\CA_k(k)\rangle$.
But the left mutation through an admissible subcategory induces
an equivalence of its left orthogonal to its right orthogonal
by lemma~\ref{rl_mut}, and $\fa_k(k+1)$ lies in the left orthogonal
to $\langle \CA_1(1),\dots,\CA_k(k)\rangle$ by (\ref{faklo}).
\end{proof}

\begin{lemma}\label{soda01}
We have the following semiorthogonal decomposition of $\CA_0$
$$
\langle \alpha_0^*(\fa_{0}(1)),\alpha_0^*(\fa_{1}(2)),\dots,
\alpha_0^*(\fa_{\ix-1}(\ix))\rangle.
$$
\end{lemma}
\begin{proof}
For any $F\in\CA_0$, $F'\in\fa_l$ we have
$$
\Hom(\alpha_0^*(F'(l+1)),F) =
\Hom(F'(l+1),F) =
\Hom(F,\FS_X (F'(l+1)))^* =
\Hom(F,\alpha_0^!\FS_X (F'(l+1)))^*,
$$
therefore $(\alpha_0^*(\fa_l(l+1)))^\perp = {}^\perp(\alpha_0^!\FS_X(\fa_l(l+1)))$.
Thus for the semiorthogonality we should check that for any $k < l$ we have
$\Hom(\alpha_0^*(\fa_k(k+1)),\alpha_0^!\FS_X(\fa_l(l+1))) = 0$. For this we note
that the inclusion~(\ref{faklo}) (with $k$ replaced by $l$) implies that
$\fa_l(l+1) \subset \langle\CA_{l+1}(l+1),\dots,\CA_{\ix-1}(\ix-1),\FS_X^{-1}\CA_0\rangle$
by lemma~\ref{ser_loro}, hence
$$
\FS_X(\fa_l(l+1)) \subset \langle\FS_X(\CA_{l+1}(l+1)),\dots,\FS_X(\CA_{\ix-1}(\ix-1)),\CA_0\rangle.
$$
Comparing this with the inclusion~(\ref{fakro}) for $\fa_k(k+1)$
and taking into account that by lemma~\ref{ser_loro} we have
a semiorthogonal decomposition
$\D^b(X) = \langle\FS_X(\CA_{l+1}(l+1)),\dots,\FS_X(\CA_{\ix-1}(\ix-1)),\CA_0,\CA_1(1),\dots,\CA_l(l)\rangle$,
we deduce that
$\Hom(\alpha_0^*(\fa_k(k+1)),\alpha_0^!\FS_X(\fa_l(l+1))) =
\Hom(\fa_k(k+1),\FS_X(\fa_l(l+1)))$ which by the Serre duality
is dual to $\Hom(\fa_l(l+1),\fa_k(k+1)) = \Hom(\fa_l(l),\fa_k(k))$
which is zero since $\fa_l(l) \subset \CA_l(l)$ and $\fa_k(k) \subset \CA_k(k)$ .

Now assume that $F$ lies in the right orthogonal to the collection
$\langle \alpha_0^*(\fa_{0}(1)),\alpha_0^*(\fa_{1}(2)),\dots,
\alpha_0^*(\fa_{\ix-1}(\ix))\rangle$ in $\CA_0$. By adjunction
$\alpha_0(F)$ is in the right orthogonal to
$\langle \fa_{0}(1),\fa_{1}(2),\dots,\fa_{\ix-1}(\ix)\rangle$
in $\D^b(X)$. But
$\alpha_0(F)\in\CA_0 = \langle\CA_1(1),\dots,\CA_{\ix-1}(\ix-1)\rangle^\perp$,
therefore
$$
\alpha_0(F) \in \langle \fa_0(1),\CA_1(1),\fa_1(2),\CA_2(2),\dots,
\fa_{\ix-2}(\ix-1),\CA_{\ix-1}(\ix-1),\fa_{\ix-1}(\ix)\rangle^\perp.
$$
It remains to note that by definition of subcategories $\fa_0$, \dots, $\fa_{\ix-1}$
we have
$\langle\fa_0(1),\CA_1(1)\rangle = \CA_0(1)$,
$\langle \fa_1(2),\CA_2(2)\rangle = \CA_1(2)$, \dots,
$\langle \fa_{\ix-2}(\ix-1),\CA_{\ix-1}(\ix-1)\rangle = \CA_{\ix-2}(\ix-1)$, and
$\fa_{\ix-1}(\ix) = \CA_{\ix-1}(\ix)$,
so we see that
$\alpha_0(F) \in \langle\CA_0(1),\CA_1(2),\dots,\CA_{\ix-1}(\ix)\rangle^\perp$
which means that $F=0$ since $\langle\CA_0(1),\CA_1(2),\dots,\CA_{\ix-1}(\ix)\rangle$
is evidently a semiorthogonal decomposition of $\D^b(X)$.
\end{proof}

\begin{lemma}\label{use}
We have
$\alpha_0^*\left(\lan\CA_0(1),\dots,\CA_{r-1}(r)\ran\right) \in
\lan \alpha_0^*(\fa_{0}(1)),\dots,\alpha_0^*(\fa_{r-1}(r))\ran$.
\end{lemma}
\begin{proof}
We have $\CA_k(k+1) = \lan\fa_k(k+1),\CA_{k+1}(k+1)\ran$
and $\alpha_0^*(\CA_{k+1}(k+1)) = 0$ for any $0 \le k \le r-1$.
\end{proof}

\begin{lemma}\label{a0rgen}
Triangulated subcategory of $\D^b(X)$ generated by $\CA_0,\CA_0(1),\dots,\CA_0(r-1)$
coincides with $\lan\CA_0,\CA_1(1),\dots,\CA_{r-1}(r-1)\ran$.
\end{lemma}
\begin{proof}
It is clear that the latter category lies in the former. On the other hand,
it is clear that
$\CA_0,\CA_0(1),\dots,\CA_0(r-1) \subset
\lan\CA_r(r),\dots,\CA_{\ix-1}(\ix-1)\ran^\perp =
\lan\CA_0,\CA_1(1),\dots,\CA_{r-1}(r-1)\ran$.
\end{proof}

\section{Universal hyperplane section}\label{s5}

Assume that $X$ is a smooth projective variety
with an effective line bundle $\CO_X(1)$ on $X$ and
assume that we are given a Lefschetz decomposition~(\ref{dbx}) of its derived category.
Let $V^* \subset \Gamma(X,\CO_X(1))$ be a vector space of global sections.
Put $N = \dim V$. We assume that
\begin{equation}\label{ngi}
N > \ix.
\end{equation}

Consider the product $X\times\PP(V^*)$. Let $\CA_k(k)\boxtimes\D^b(\PP(V^*))$
denote the triangulated subcategory of $\D^b(X\times\PP(V^*))$ generated
by objects $F\boxtimes G$ with $F\in\CA_k(k) \subset \D^b(X)$ and
$G\in\D^b(\PP(V^*))$. Note that this category coincides with the category
$\lan \CA_k(k)\boxtimes\CO_{\PP(V^*)},\CA_k(k)\boxtimes\CO_{\PP(V^*)}(1),\dots,
\CA_k(k)\boxtimes\CO_{\PP(V^*)}(N-1)\ran$. Every term in the latter
decomposition is equivalent to $\CA_k$, hence saturated, hence admissible, therefore
$\CA_k(k)\boxtimes\D^b(\PP(V^*))$ is also admissible and saturated.
Moreover, it is clear that we have the following semiorthogonal decomposition
$$
\D^b(X\times\PP(V^*)) = \lan \CA_0\boxtimes\D^b(\PP(V^*)),
\CA_1(1)\boxtimes\D^b(\PP(V^*)),\dots,\CA_{\ix-1}(\ix-1)\boxtimes\D^b(\PP(V^*))\ran.
$$
Note also that by the K\"uneth formula we have
$$
\RHom_{X\times\PP(V^*)}(F_1\boxtimes G_1,F_2\boxtimes G_2) \cong
\RHom_X(F_1,F_2)\otimes\RHom_{\PP(V^*)}(G_1,G_2)
$$
for all $F_1,F_2\in\D^b(X)$, $G_1,G_2\in\D^b(\PP(V^*))$.

Consider the universal hyperplane section of $X$,
that is the zero locus $\CX_1 \subset X\times\PP(V^*)$
of the canonical section of a line bundle $\CO_X(1)\boxtimes\CO_{\PP(V^*)}(1)$.
Let $\pi:\CX_1\to X$ and $f:\CX_1 \to \PP(V^*)$ denote the projections,
and let $i:\CX_1 \to X \times \PP(V^*)$ denote the embedding.
Note that $\CX_1\subset X\times \PP(V^*)$ is a divisor of bidegree $(1,1)$
and we have the following resolution of its structure sheaf
\begin{equation}\label{ox1}
0 \to \CO_X(-1)\boxtimes\CO_{\PP(V^*)}(-1) \to
\CO_{X\times\PP(V^*)} \to i_*\CO_{\CX_1} \to 0.
\end{equation}

The following lemma is useful for calculations of $\Hom$'s between
decomposable objects in $\D^b(\CX_1)$.

\begin{lemma}\label{dechoms}
For any $F_1,F_2 \in \D^b(X)$, $G_1,G_2 \in \D^b(\PP(V^*))$ we have
an exact triangle
\begin{multline*}
\RHom_X(F_1,F_2(-1)) \otimes \RHom_{\PP(V^*)}(G_1,G_2(-1)) \to
\\ \to
\RHom_X(F_1,F_2) \otimes \RHom_{\PP(V^*)}(G_1,G_2) \to
\\ \to
\RHom_{\CX_1}(i^*(F_1\boxtimes G_1),i^*(F_2\boxtimes G_2)).
\end{multline*}
\end{lemma}
\begin{proof}
Tensoring resolution~(\ref{ox1}) by
$(F_1^*\otimes F_2)\boxtimes(G_1^*\otimes G_2)$
and applying $\RG$ we obtain the following exact triangle
\begin{multline*}
\RG(X\times\PP(V^*),(F_1^*\otimes F_2(-1))\boxtimes(G_1^*\otimes G_2(-1))) \to
\\ \to
\RG(X\times\PP(V^*),(F_1^*\otimes F_2)\boxtimes(G_1^*\otimes G_2)) \to
\\ \to
\RG(\CX_1,i^*((F_1^*\otimes F_2)\boxtimes(G_1^*\otimes G_2))).
\end{multline*}
Rewriting $\RG$ in terms of $\RHom$'s and applying K\"uneth formula
we obtain the desired triangle.
\end{proof}

\begin{corollary}\label{pisff}
The functor $\pi^*:\D^b(X) \to \D^b(\CX_1)$ is fully faithful.
Moreover, for any $F_1,F_2\in\D^b(X)$ and $1\le k\le N-2$ we have
$\RHom_{\CX_1}(\pi^*F_1,\pi^*F_2\otimes f^*\CO_{\PP(V^*)}(-k)) = 0$.
\end{corollary}
\begin{proof}
Take $G_1 = \CO_{\PP(V^*)}$, $G_2 = \CO_{\PP(V^*)}(-k)$.
Then we have isomorphisms
$i^*(F_1\boxtimes G_1) \cong \pi^*F_1$,
and
$i^*(F_2\boxtimes G_2) \cong \pi^*F_2\otimes f^*\CO_{\PP(V^*)}(-k)$.
Since
$\RHom_{\PP(V^*)}(\CO_{\PP(V^*)},\CO_{\PP(V^*)}(-k)) = 0$ for $1\le k \le N-1$
and
$\RHom_{\PP(V^*)}(\CO_{\PP(V^*)},\CO_{\PP(V^*)}) = \kk$,
the first term in the triangle of the lemma vanishes for $0\le k\le N-2$ and
the second term equals $\RHom_X(F_1,F_2)$ for $k=0$ and vanishes
for $1\le k\le N-1$ whereof we obtain the claim.
\end{proof}

\begin{lemma}\label{sodx1}
For any $1\le k\le \ix-1$ the functor
$\CA_k(k)\boxtimes\D^b(\PP(V^*)) \subset \D^b(X\times\PP(V^*)) \exto{i^*} \D^b(\CX_1)$
is fully faithful,
and the collection
$$
(\CA_1(1)\boxtimes\D^b(\PP(V^*)),\dots,\CA_{\ix-1}(\ix-1)\boxtimes\D^b(\PP(V^*)))
\subset \D^b(\CX_1)
$$
is semiorthogonal.
\end{lemma}
\begin{proof}
Let $1\le k\le l\le \ix-1$, take $F_1\in\CA_l(l)$, $F_2\in\CA_k(k)$,
$G_1,G_2\in\D^b(\PP(V^*))$ and consider the triangle of lemma~\ref{dechoms}.
Its first term vanishes since $F_1\in\CA_l(l)$ and
$F_2(-1)\in\CA_k(k-1) \subset \CA_{k-1}(k-1)$.
Therefore, in the case $k=l$ we see that the functor
$i^*:\CA_k(k)\boxtimes\D^b(\PP(V^*)) \to \D^b(\CX_1)$ is fully faithful.
On the other hand, for $1\le k < l\le \ix-1$ the second term vanishes as well,
since $F_1\in\CA_l(l)$ and $F_2\in\CA_k(k)$. Therefore the above collection
is semiorthogonal.

It remains to check that categories  $\CA_k(k)\boxtimes\D^b(\PP(V^*))$ are
admissible in $\D^b(\CX_1)$. For this we note that they are saturated,
hence admissible in $\D^b(\CX_1)$.
\end{proof}

Let $\CC$ denote the right orthogonal to the subcategory
$\langle \CA_1(1)\boxtimes\D^b(\PP(V^*)),\dots,
\CA_{\ix-1}(\ix-1)\boxtimes\D^b(\PP(V^*))\rangle$
in $\D^b(\CX_1)$,
\begin{equation}\label{defcc}
\CC = \langle \CA_1(1)\boxtimes\D^b(\PP(V^*)),\dots,
\CA_{\ix-1}(\ix-1)\boxtimes\D^b(\PP(V^*))\rangle^\perp_{\D^b(X)}.
\end{equation}
Let $\gamma:\CC \to \D^b(\CX_1)$ denote the inclusion functor.
Since the subcategories $\CA_1(1)\boxtimes\D^b(\PP(V^*))$, \dots,
$\CA_{\ix-1}(\ix-1)\boxtimes\D^b(\PP(V^*))$ are admissible it follows that
$\CC$ is left admissible, hence the functor $\gamma$ has a left adjoint
functor $\gamma^*:\D^b(\CX_1) \to \CC$.

Note that the subcategory
$\langle \CA_1(1)\boxtimes\D^b(\PP(V^*)),\dots,
\CA_{\ix-1}(\ix-1)\boxtimes\D^b(\PP(V^*))\rangle \subset \D^b(\CX_1)$
is $\PP(V^*)$-linear.
In particular, the functor $F\mapsto F\otimes f^*\CO_{\PP(V^*)}(1)$,
$\D^b(\CX_1) \to \D^b(\CX_1)$ restricts to an endofunctor of $\CC$
which we denote simply by $F \mapsto F(1)$.


%

Consider the composition of functors $\pi_*\circ\gamma:\CC \to \D^b(X)$.

\begin{lemma}\label{impig}
The image of the functor $\pi_*\circ\gamma$ is contained
in the strictly full subcategory $\CA_0\subset\D^b(X)$.
\end{lemma}
\begin{proof}
If $F\in\CA_k(k)$, $1\le k\le \ix-1$, and $F'\in\CC$ then we have
$\Hom(F,\pi_*(\gamma(F'))) = \Hom(\pi^*F,\gamma(F')) = 0$
since $\pi^*F\in\CA_k(k)\boxtimes\D^b(\PP(V^*))$.
Thus $\pi_*(\gamma(F'))$ is contained in the right orthogonal
to the subcategory $\langle\CA_1(1),\dots,\CA_{\ix-1}(\ix-1)\rangle$,
which by~(\ref{dbx}) coincides with $\CA_0$.
\end{proof}

Consider the functor $\gamma^*\circ\pi^*:\D^b(X) \to \CC$ which is left adjoint
to $\pi_*\circ\gamma$.
In proposition~\ref{gps_ff} below we will show that the restriction
of this functor to the subcategory $\CA_0 \subset \D^b(X)$ is fully faithful.
We start with two lemmas.

For any object $F\in\D^b(X)$ consider the decomposition of $\pi^*F\in\D^b(\CX_1)$
\begin{equation}\label{def_f1}
F_\PCC \to \pi^*F \to F_\CC
\end{equation}
with $F_\CC \in \CC$, $F_\PCC \in \PCC = \langle \CA_1(1)\boxtimes\D^b(\PP(V^*)),\dots,
\CA_{\ix-1}(\ix-1)\boxtimes\D^b(\PP(V^*))\rangle$. Then it is clear that
\begin{equation}\label{def_fcc}
F_\CC = \gamma\gamma^*\pi^*F.
\end{equation}

\begin{lemma}\label{hompis}
If $\RHom(F,\pi_*(F'_\PCC(0,k))) = 0$ then we have a canonical isomorphism
$$
\Hom_{\D^b(\CX_1)}(\pi^*F,\pi^*F'(0,k)) \cong
\Hom_\CC(\gamma^*\pi^*F,(\gamma^*\pi^*F')(k)).
$$
\end{lemma}
\begin{proof}
Applying the functor $\RHom(\pi^*F,-)$ to the exact triangle
$F'_\PCC(0,k) \to \pi^*F'(0,k) \to F'_\CC(0,k)$
and taking into account the isomorphism
$\RHom(\pi^*F,F'_\PCC(0,k)) \cong \RHom(F,\pi_*(F'_\PCC(0,k))) = 0$
we deduce
$\Hom(\pi^*F,\pi^*F'(0,k)) \cong \Hom(\pi^*F,F'_\CC(0,k))$.
It remains to note that
$$
\Hom_{\D^b(\CX_1)}(\pi^*F,F'_\CC(0,k)) \cong
\Hom_{\D^b(\CX_1)}(\pi^*F,\gamma\gamma^*\pi^*F'(0,k)) \cong
\Hom_\CC(\gamma^*\pi^*F,\gamma^*\pi^*F'(k)).
$$
\end{proof}

Recall
a semiorthogonal decomposition
$\CA_0 = \langle \alpha_0^*(\fa_{0}(1)),\dots,\alpha_0^*(\fa_{\ix-1}(\ix))\rangle$
constructed in lemma~\ref{soda01}.

\begin{lemma}\label{f1}
Let $F\in \langle\alpha_0^*(\fa_0(1)),\dots,\alpha_0^*(\fa_k(k+1))\rangle
\subset \CA_0\subset\D^b(X)$. Then
$$
F_\PCC \in
\left\langle\begin{array}{r}
\CA_1(1)\boxtimes\CO_{\PP(V^*)}(-k),
\CA_1(1)\boxtimes\CO_{\PP(V^*)}(1-k),\dots,
\CA_1(1)\boxtimes\CO_{\PP(V^*)}(-1)\\
\CA_2(2)\boxtimes\CO_{\PP(V^*)}(1-k),\dots,
\CA_2(2)\boxtimes\CO_{\PP(V^*)}(-1)\\
\vdots\qquad\qquad\qquad\\
\CA_{k}(k)\boxtimes\CO_{\PP(V^*)}(-1)
\end{array}\right\rangle.
$$
\end{lemma}
\begin{proof}
Consider a decomposition of $\pi^*F$ with respect to the semiorthogonal
decomposition
$\D^b(\CX_1) = \lan \CC,\CA_1(1)\boxtimes\D^b(\PP(V^*)),\dots,
\CA_{\ix-1}(\ix-1)\boxtimes\D^b(\PP(V^*))\ran$.
Then $F_\PCC$ is glued from its components in the subcategories
$\CA_1(1)\boxtimes\D^b(\PP(V^*))$, \dots, $\CA_{\ix-1}(\ix-1)\boxtimes\D^b(\PP(V^*))$.
First of all let us compute the component of $\pi^*F$ in
$\CA_{\ix-1}(\ix-1)\boxtimes\D^b(\PP(V^*))$. It is given by applying to $\pi^*F$
the right adjoint functor to the inclusion functor
$\CA_{\ix-1}(\ix-1)\boxtimes\D^b(\PP(V^*)) \to \D^b(\CX_1)$.
To compute this we take $F_1\in\CA_{\ix-1}(\ix-1)$, $G_1\in\D^b(\PP(V^*))$,
$F_2 = F$, $G_2 = \CO_{\PP(V^*)}$ and consider the triangle
of lemma~\ref{dechoms}. The second term of this triangle vanishes because
$F_1\in\CA_{\ix-1}(\ix-1)$ and $F_2\in\CA_0$
and the first term vanishes because
$$
\RHom_X(F_1,F_2(-1)) = \RHom_X(F_1(1),F_2) = \RHom_{\CA_0}(\alpha_0^*(F_1(1)),F_2)
$$
and since $F_1\in\CA_{\ix-1}(\ix-1) = \fa_{\ix-1}(\ix-1)$ we have
$\alpha_0^*(F_1(1)) \in \alpha_0^*(\fa_{\ix-1}(\ix))$ which is orthogonal
to the category $\langle\alpha_0^*(\fa_0(1)),\dots,\alpha_0^*(\fa_k(k+1))\rangle$
by lemma~\ref{soda01}. Therefore the component of the object $\pi^*F$ in the category
$\CA_{\ix-1}(\ix-1)\boxtimes\D^b(\PP(V^*))$ is zero.
Similar arguments show that the components of $\pi^*F$ in
$\CA_{\ix-2}(\ix-2)\boxtimes\D^b(\PP(V^*))$, \dots,
$\CA_{k+1}(k+1)\boxtimes\D^b(\PP(V^*))$ are also zero.

Now let us compute the component of $\pi^*F$
in $\CA_{k}(k)\boxtimes\D^b(\PP(V^*))$.
It is given by applying to $\pi^*F$
the right adjoint functor to the inclusion functor
$\CA_{k}(k)\boxtimes\D^b(\PP(V^*)) \to \D^b(\CX_1)$.
To compute this we take $F_1\in\CA_{k}(k)$, $G_1\in\D^b(\PP(V^*))$,
$F_2 = F$, $G_2 = \CO_{\PP(V^*)}$ and consider the triangle
of lemma~\ref{dechoms}.
Note again that the second term vanishes because $F_1\in\CA_{k}(k)$
and $F_2\in\CA_0$, hence
$$
\RHom_{\CX_1}(F_1\boxtimes G,\pi^*F) \cong
\RHom_X(F_1,F(-1)) \otimes \RHom_{\PP(V^*)}(G,\CO_{\PP(V^*)}(-1))[1],
$$
and since the embedding of $\CA_{k}(k)\boxtimes\D^b(\PP(V^*))$ into
$\D^b(\CX_1)$ is fully faithful by lemma~\ref{sodx1}, we conclude that
the component of $\pi^*F$ in $\CA_{k}(k)\boxtimes\D^b(\PP(V^*))$
equals $\alpha_k\alpha_{k}^!(F(-1))\boxtimes\CO_{\PP(V^*)}(-1)[1]$.

Now let us compute the component of $\pi^*F$
in $\CA_{k-1}(k-1)\boxtimes\D^b(\PP(V^*))$.
It is given by applying to the fiber of the morphism
$\alpha_k\alpha_{k}^!(F(-1))\boxtimes\CO_{\PP(V^*)}(-1)[1] \to \pi^*F$
the right adjoint functor to the inclusion functor
$\CA_{k-1}(k-1)\boxtimes\D^b(\PP(V^*)) \to \D^b(\CX_1)$.
To compute this we take $F_1\in\CA_{k-1}(k-1)$, $G_1\in\D^b(\PP(V^*))$,
and either $F_2 = \alpha_k\alpha_{k}^!(F(-1))$, $G_2 = \CO_{\PP(V^*)}(-1)[1]$,
or $F_2 = \pi^*F$, $G_2 = 0$.
Repeating the above arguments we see that the projection of the target equals
$\alpha_{k-1}\alpha_{k-1}^!(F(-1))\boxtimes\CO_{\PP(V^*)}(-1)[1]$,
and the projection of source is the cone of the morphism
$\alpha_k\alpha_{k}^!(F(-1))(-1)\boxtimes\CO_{\PP(V^*)}(-2)[1] \to
\alpha_{k-1}\alpha_{k-1}^!\alpha_k\alpha_{k}^!(F(-1))\boxtimes\CO_{\PP(V^*)}(-1)[1]$.
Both projections are contained in $\CA_{k-1}(k-1)\boxtimes
\langle \CO_{\PP(V^*)}(-2),\CO_{\PP(V^*)}(-1)\rangle$,
hence the same is true for the corresponding component of $F_\PCC$.

Proceeding in the same manner we deduce the rest of the lemma.
\end{proof}

Now we can prove

\begin{proposition}\label{gps_ff}
The restriction of the functor $\gamma^*\circ\pi^*:\D^b(X) \to \CC$
to the subcategory $\CA_0\subset\D^b(X)$ is fully faithful.
\end{proposition}
\begin{proof}
Take $F,F'\in\CA_0$. Then
$F'_\PCC \in \pi^*(\D^b(X))\otimes\lan\CO_{\PP(V^*)}(1-\ix),\dots,\CO_{\PP(V^*)}(-1)\ran$
by lemma~\ref{f1}.
On the other hand, by corollary~\ref{pisff} the pushforward functor
$\pi_*:\D^b(\CX_1) \to \D^b(X)$ takes the subcategory
$\pi^*(\D^b(X))\otimes\lan\CO_{\PP(V^*)}(2-N),\dots,\CO_{\PP(V^*)}(-1)\ran$
to zero. Since $N > \ix$ by assumption~(\ref{ngi}), we conclude
that $\pi_*F'_\PCC = 0$.
Therefore, by lemma~\ref{hompis} we have
$\Hom(\gamma^*\pi^*F,\gamma^*\pi^*F') \cong \Hom(\pi^*F,\pi^*F')$
which is isomorphic to $\Hom(F,F')$ by corollary~\ref{pisff}.
\end{proof}

The following corollary is not needed below, however we put it here
as an illustration.

\begin{corollary}
The functors $\pi_*\gamma:\CC \to \D^b(X)$ and $\gamma^*\pi^*:\D^b(X) \to \CC$
are splitting.
\end{corollary}
\begin{proof}
The first functor is right adjoint to the second one, hence it suffices to check
that $\gamma^*\pi^*$ is splitting. Since $\Im(\pi_*\gamma) \subset \CA_0$
by lemma~\ref{impig} and $\gamma^*\pi^*$ is fully faithful on $\CA_0$
we deduce that
$\Ker(\gamma^*\pi^*) = {}^\perp\CA_0 = \lan\CA_1(1),\dots,\CA_{\ix-1}(\ix-1)\ran$
hence is admissible.
Moreover, $\Im(\gamma^*\pi^*)$ is equivalent to $\CA_0$ and $\CA_0$ is saturated,
hence $\Im(\gamma^*\pi^*)$ is admissible.
\end{proof}

Let
\begin{equation}\label{defj}
\jx = N - 1 - \max\{ k\ |\ \CA_k = \CA_0\}
\end{equation}
(note that $\jx > 0$ by (\ref{ngi})).
Then $\fa_k = 0$ for $k < N-1-\jx$.
Consider the subcategories
\begin{equation}\label{defcb}
\CB_k = \gamma^*(\pi^*(
\langle\alpha_0^*(\fa_0(1)),\dots,\alpha_0^*(\fa_{N-k-2}(N-k-1))\rangle
)) \subset \CB_0 = \gamma^*(\pi^*(\CA_0)),
\end{equation}
where we put $\fa_l = 0$ for $l\ge \ix$ for convenience.
Note that
\begin{equation}\label{bk0}
\CB_k = 0
\qquad
\text{for $k\ge \jx$}.
\end{equation}
If the initial Lefschetz decomposition is rectangular then $\jx = N - \ix$
and $\CB_{\jx-1} = \CB_{\jx-2} = \dots = \CB_1 = \CB_0$.

\begin{lemma}
The chain of subcategories
$$
0 = \CB_{\jx-1} \subset \CB_{\jx-2} \subset \dots \subset \CB_1 \subset \CB_0 \subset \CC
$$
is a chain of admissible subcategories in $\CC$.
\end{lemma}
\begin{proof}
The category $\lan\alpha_0^*(\fa_0(1)),\dots,\alpha_0^*(\fa_{N-k-2}(N-k-1))\ran$
is generated by a semiorthogonal collection (see lemma~\ref{soda01})
of admissible subcategories of $\CA_0$, hence admissible in $\CA_0$,
hence saturated. Therefore its image under fully faithful functor
$\gamma^*\pi^*:\CA_0 \to \CC$ is admissible.
\end{proof}

Almost the same arguments show the following

\begin{proposition}\label{soc_cb}
The collection
$\langle \CB_{\jx-1}(1-\jx),\CB_{\jx-2}(2-\jx),\dots,\CB_1(-1),\CB_0\rangle$
is semiorthogonal in $\CC$.
\end{proposition}
\begin{proof}
Take $F\in\CA_0$,
$F' \in \langle\alpha_0^*(\fa_0(1)),\dots,\alpha_0^*(\fa_{N-k-2}(N-k-1))\rangle$.
Then by lemma~\ref{f1} we have an inclusion
$F'_\PCC \in \pi^*(\D^b(X))\boxtimes\langle\CO(0,k+2-N),\dots,\CO(0,-1)\rangle$,
hence $\pi_*(F'_\PCC(0,-k)) = 0$ by corollary~\ref{pisff}.
Therefore, by lemma~\ref{hompis} we have an isomorphism
$\Hom(\gamma^*\pi^*F,\gamma^*\pi^*F'(0,-k)) = \Hom(\pi^*F,\pi^*F'(0,-k))$
which equals zero for $1\le k\le \jx-1$ by corollary~\ref{pisff}
(note that $\jx-1 \le N-2$ by~(\ref{defj})).
\end{proof}

In the following section we prove that the semiorthogonal collection
of the proposition generates $\CC$. For the proof we use some additional
assumptions, though the fact must be true without them. It would be interesting
to find a direct proof.

We conclude the section with a couple of lemmas that will be useful later.

\begin{lemma}\label{imkerpig}
We have
$\Im\pi_*\gamma = \CA_0$, $\Ker\gamma^*\pi^* = \lan\CA_1(1),\dots,\CA_{\ix-1}(\ix-1)\ran$,
and similarly
$\Im\gamma^*\pi^* = \CB_0$, $\Ker\pi_*\gamma \supset \lan\CB_{\jx-1}(1-\jx),\dots,\CB_1(-1)\ran$.
\end{lemma}
\begin{proof}
Since $\Im\pi_*\gamma \subset \CA_0$ by lemma~\ref{impig} and
the functor $\gamma^*\pi^*$ (which is left adjoint to $\pi_*\gamma$)
is fully faithful on $\CA_0$ by lemma~\ref{gps_ff},
we deduce that $\Im \pi_*\gamma = \CA_0$,
$\Ker\gamma^*\pi^* = {}^\perp\CA_0 = \lan\CA_1(1),\dots,\CA_{\ix-1}(\ix-1)\ran$,
and $\Im\gamma^*\pi^* = \gamma^*\pi^*(\CA_0) = \CB_0$.
Therefore $\Ker\pi_*\gamma = \CB_0^\perp \supset \lan\CB_{\jx-1}(1-\jx),\dots,\CB_1(-1)\ran$
by proposition~\ref{soc_cb}.
\end{proof}

\begin{lemma}\label{use1}
We have $\gamma^*\pi^*(\lan\CA_{0}(1),\dots,\CA_{r-2}(r-1)\ran^\perp) \subset \CB_{N-r}$.
\end{lemma}
\begin{proof}
Since $\lan\CA_1(1),\dots,\CA_{\ix-1}(\ix-1)\ran = \Ker(\gamma^*\pi^*)$
we have $\gamma^*\pi^*(F) = \gamma^*\pi^*\alpha_0^*(F)$ for any $F\in\D^b(X)$.
On the other hand, by lemma~\ref{use} we have
$\alpha_0^*\left(\lan\CA_{0}(1),\dots,\CA_{r-2}(r-1)\ran^\perp\right) \subset
\lan\alpha_0^*(\fa_0(1)),\dots,\alpha_0^*(\fa_{r-2}(r-1))\ran$
and by definition of $\CB_{N-r}$ we have
$\gamma^*\pi^*(\lan\alpha_0^*(\fa_0(1)),\dots,\alpha_0^*(\fa_{r-2}(r-1))\ran) = \CB_{N-r}$.
%
\end{proof}

\begin{remark}
Consider any smooth base scheme $S$ (not necessarily compact)
and assume that $p:X \to S$ is an algebraic $S$-variety,
projective over $S$ with an $S$-linear Lefschetz decomposition.
Then all results of this section can be proved by essentially
the same arguments. We only should replace $\RHom_X$ by $p_*\RCHom_X$,
and the Serre functor of $X$ by the relative Serre functor of $X$
over $S$ (see \cite{K}).
\end{remark}

%
%

\section{Homological projective duality}

Recall the assumptions of the previous section:
we have a smooth projective variety $X$ with an effective line bundle $\CO_X(1)$,
a Lefschetz decomposition~(\ref{dbx}) of its derived category,
a vector space of global sections $V^* \subset \Gamma(X,\CO_X(1))$
such that (\ref{ngi}) holds, i.e. $N = \dim V > \ix$
($\ix$ is the number of terms in the Lefschetz decomposition).
Assume also that the space $V^*$ generates $\CO_X(1)$, so that
we have a regular morphism $f:X \to \PP(V)$.

Recall that we denoted by $\CX_1 \subset X\times\PP(V^*)$ the universal hyperplane section
of $X$ and by $\CC$ the right orthogonal to the subcategory
$\langle \CA_1(1)\boxtimes\D^b(\PP(V^*)),\dots,
\CA_{\ix-1}(\ix-1)\boxtimes\D^b(\PP(V^*))\rangle$
in $\D^b(\CX_1)$,
$$
\CC = \langle \CA_1(1)\boxtimes\D^b(\PP(V^*)),\dots,
\CA_{\ix-1}(\ix-1)\boxtimes\D^b(\PP(V^*))\rangle^\perp_{\D^b(\CX_1)}.
$$
Note that the category $\CC$ is a module category over the tensor category
$\D^b(\PP(V^*))$: if $F\in\CC$ and $G\in\D^b(\PP(V^*))$ then
$F\otimes f^*G \in \CC$.

Now we need some additional assumptions. Assume that $\CC$ is {\em geometrical},
meaning that $\CC$ is equivalent to the derived category of coherent sheaves
on some algebraic variety $Y$.
Let $\Phi:\D^b(Y) \to \D^b(\CX_1)$ denote the composition of the equivalence
$\D^b(Y) \to \CC$ with the inclusion functor $\gamma:\CC \to \D^b(\CX_1)$.
Further assume that the module structure on $\CC$ is {\em geometrical},
meaning that there is an algebraic morphism $g:Y \to \PP(V^*)$, such that
there is an isomorphism of bifunctors
$$
\Phi(F\otimes g^*G) \cong \Phi(F)\otimes f^* G,
\qquad F\in\D^b(Y),\ G\in\D^b(\PP(V^*)).
$$
In other words, the functor $\Phi$ is assumed to be $\PP(V^*)$-linear.
Note also that the functor $\Phi:\D^b(Y) \to \D^b(\CX_1)$ is fully faithful,
hence by Orlov's Theorem \cite{O3} it can be represented by a kernel
on $Y\times\CX_1$. Moreover, it is easy to see that $\PP(V^*)$-linearity
of the functor $\Phi$ implies that the kernel of $\Phi$ is supported
set-theoretically on the fiber product $Y\times_{\PP(V^*)}\CX_1$.
Actually, it is natural to conjecture that the kernel is supported
even {\em scheme-theoretically} on the fiber product (this must be
a relative version of the Orlov's Theorem). However, we don't address
this question here, taking this as an additional assumption.
Finally, note that
$$
Y\times_{\PP(V^*)}\CX_1 = Q(X,Y)
:= (X\times Y) \times_{\scriptscriptstyle\PP(V)\times\PP(V^*)}Q,
$$
where
$Q = \{(v,H)\in\PP(V)\times\PP(V^*)\ |\ v\in H\}$
is the incidence quadric.

\begin{definition}
An algebraic variety $Y$ with a projective morphism $g:Y\to\PP(V^*)$
is called {\sf Homologically Projectively Dual} to $f:X\to \PP(V)$
with respect to a Lefschetz decomposition~$(\ref{dbx})$, if
there exists an object $\CE\in\D^b(Q(X,Y))$ such that the functor
$\Phi = \Phi_\CE:\D^b(Y) \to \D^b(\CX_1)$ is fully faithful
and gives the following semiorthogonal decomposition
\begin{equation}\label{dbx1}
\D^b(\CX_1) = \langle \Phi(\D^b(Y)),\CA_1(1)\boxtimes\D^b(\PP(V^*)),\dots,
\CA_{\ix-1}(\ix-1)\boxtimes\D^b(\PP(V^*))\rangle.
\end{equation}
\end{definition}

In the next section we will reveal relation of the Homological Projective
Duality to the classical projective duality, and now we will state and
prove the main theorem about Homologically Projectively Dual varieties.

For every linear subspace $L \subset V^*$ we consider
the corresponding linear sections of $X$ and $Y$:
$$
X_L = X\times_{\PP(V)}\PP(L^\perp),
\qquad
Y_L = Y\times_{\PP(V^*)}\PP(L),
$$
where $L^\perp \subset V$ is the orthogonal subspace to $L\subset V^*$.

\begin{definition}
A subspace $L\subset V^*$ is called {\sf admissible}, if
\begin{enumerate}
\renewcommand{\theenumi}{\alph{enumi}}
\renewcommand{\labelenumi}{(\theenumi)}
\item $\dim X_L = \dim X - \dim L$, and
\item $\dim Y_L = \dim Y + \dim L - N$
\end{enumerate}
\end{definition}

The main result of this paper is the following

\begin{theorem}\label{themain}
If\/ $Y$ is Homologically Projectively Dual to $X$ then

\noindent$(i)$
$Y$ is smooth and $\D^b(Y)$ admits a dual Lefschetz decomposition
$$
\D^b(Y) = \lan \CB_{\jx-1}(1-\jx),\dots,\CB_{1}(-1),\CB_0\ran,\qquad
0 \subset \CB_{\jx-1} \subset \dots \subset \CB_1 \subset \CB_0 \subset \D^b(Y)
$$
with the same set of primitive subcategories:
$\CB_k = \langle\fa_0,\dots,\fa_{N-k-2}\rangle$;

\noindent$(ii)$
for any admissible linear subspace $L\subset V^*$, $\dim L = r$,
there exist a triangulated category $\CC_L$ and semiorthogonal decompositions
$$
\begin{array}{lll}
\D^b(X_L) &=& \langle \CC_L,\CA_{r}(1),\dots,\CA_{\ix-1}(\ix-r)\rangle\\
\D^b(Y_L) &=& \langle \CB_{\jx-1}(N-r-\jx),\dots,\CB_{N-r}(-1),\CC_L\rangle.
\end{array}
$$
\end{theorem}

This theorem can be illustrated by the following picture:

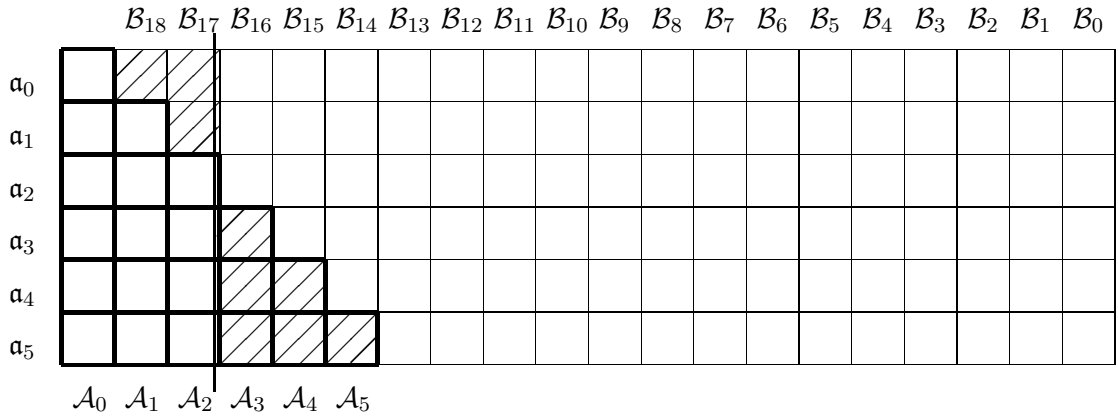
\begin{figure}[h]
\caption{Illustration of theorem~\ref{themain}: $N = 20$, $\ix = 6$, $\jx = 19$, $r = 3$}
\unitlength=7mm
\begin{picture}(25,7.5)
\linethickness{0.01mm}
\multiput(3,0)(0,1){7}{\line(1,0){20}}
\multiput(3,0)(1,0){21}{\line(0,1){6}}
\linethickness{0.5mm}
\put(3,0){\line(1,0){6}}
\put(3,1){\line(1,0){6}}
\put(3,2){\line(1,0){5}}
\put(3,3){\line(1,0){4}}
\put(3,4){\line(1,0){3}}
\put(3,5){\line(1,0){2}}
\put(3,6){\line(1,0){1}}
\put(3,0){\line(0,1){6}}
\put(4,0){\line(0,1){6}}
\put(5,0){\line(0,1){5}}
\put(6,0){\line(0,1){4}}
\put(7,0){\line(0,1){3}}
\put(8,0){\line(0,1){2}}
\put(9,0){\line(0,1){1}}
\linethickness{0.3mm}
\put(5.9,-0.5){\line(0,1){6.8}}
\put(2,5.2){$\fa_0$}
\put(2,4.2){$\fa_1$}
\put(2,3.2){$\fa_2$}
\put(2,2.2){$\fa_3$}
\put(2,1.2){$\fa_4$}
\put(2,0.2){$\fa_5$}
\put(3.2,-0.8){$\CA_0$}
\put(4.2,-0.8){$\CA_1$}
\put(5.2,-0.8){$\CA_2$}
\put(6.2,-0.8){$\CA_3$}
\put(7.2,-0.8){$\CA_4$}
\put(8.2,-0.8){$\CA_5$}
\put(22.2,6.4){$\CB_0$}
\put(21.2,6.4){$\CB_1$}
\put(20.2,6.4){$\CB_2$}
\put(19.2,6.4){$\CB_3$}
\put(18.2,6.4){$\CB_4$}
\put(17.2,6.4){$\CB_5$}
\put(16.2,6.4){$\CB_6$}
\put(15.2,6.4){$\CB_7$}
\put(14.2,6.4){$\CB_8$}
\put(13.2,6.4){$\CB_9$}
\put(12.2,6.4){$\CB_{10}$}
\put(11.2,6.4){$\CB_{11}$}
\put(10.2,6.4){$\CB_{12}$}
\put(9.2,6.4){$\CB_{13}$}
\put(8.2,6.4){$\CB_{14}$}
\put(7.2,6.4){$\CB_{15}$}
\put(6.2,6.4){$\CB_{16}$}
\put(5.2,6.4){$\CB_{17}$}
\put(4.2,6.4){$\CB_{18}$}
%
%
\put(4.0,5.45){\line(1,1){0.55}}
\put(4.0,5.0){\line(1,1){1.0}}
\put(4.5,5.0){\line(1,1){1.0}}
\put(5.0,5.0){\line(1,1){1.0}}
\put(5.0,4.5){\line(1,1){1.0}}
\put(5.0,4.0){\line(1,1){1.0}}
\put(5.45,4.0){\line(1,1){0.55}}
\put(5.95,2.45){\line(1,1){0.55}}
\put(6.0,2.0){\line(1,1){1.0}}
\put(6.0,1.5){\line(1,1){1.0}}
\put(6.0,1.0){\line(1,1){1.0}}
\put(6.0,0.5){\line(1,1){1.5}}
\put(6.0,0.0){\line(1,1){2.0}}
\put(6.5,0.0){\line(1,1){1.5}}
\put(7.0,0.0){\line(1,1){1.0}}
\put(7.5,0.0){\line(1,1){1.0}}
\put(8.0,0.0){\line(1,1){1.0}}
\put(8.45,0.0){\line(1,1){0.55}}

\end{picture}
\end{figure}
\bigskip
The bold Young diagram on the left represents the Lefschetz decomposition
of $\D^b(X)$ and the rest of the rectangle represents the dual Lefschetz
decomposition of $\D^b(Y)$. We divide the picture by a vertical line
after the column with number $r=3$. The hatched area on the right
of this line represents the part of the Lefschetz decomposition
of $\D^b(X)$ which is present in $\D^b(X_L)$, and the hatched area
on the left of this line represents the part of the Lefschetz
decomposition of $\D^b(Y)$ which is present in $\D^b(Y_L)$.

In fact the Lefschetz decomposition of $\D^b(Y)$ was constructed
in the previous section. Indeed, by definition of Homological Projective
Duality the category $\D^b(Y)$ is equivalent to $\CC$ and we have constructed
a Lefschetz collection $\CB_{\jx-1}(1-\jx),\dots,\CB_1(-1),\CB_0 \subset \CC$
in proposition~\ref{soc_cb}. The additional claim of the first part
of the theorem is that this collection generates $\D^b(Y) \cong \CC$.

The claim of the second part of the theorem can be reformulated as follows.
Derived categories of orthogonal admissible linear sections of Homologically
Projectively Dual varieties admit semiorthogonal decompositions, one part
of which comes from the Lefschetz decompositions of the ambient varieties,
and the additional parts are equivalent. This behavior, analogous
to the Lefschetz theory for cohomology of linear sections, was
a motivation for our terminology.


%

\subsection{Universal families of linear sections}

To prove the main theorem it is convenient to consider
the {\sf universal families}\/ of linear sections of $X$ and $Y$.
All $r$-dimensional subspaces $L \subset V^*$ are parameterized
by the Grassmannian $\Gr(r,V^*)$ which we denote for short by $\BP_r$.
Let $\CL_r \subset V^*\otimes\CO_{\BP_r}$ be the tautological rank $r$
subbundle on the Grassmannian $\BP_r$ and let
$\CL_r^\perp:=(V^*\otimes\CO_{\BP_r}/\CL_r)^* \subset V\otimes\CO_{\BP_r}$
be the orthogonal subbundle.
Then the universal families of linear sections of $X$ and $Y$ are
$$
\arraycolsep=2pt
\begin{array}{lll}
\CX_r & =
(X\times\BP_r)\times_{\PP(V)\times\BP_r}\PP_{\BP_r}(\CL_r^\perp) &
\subset X\times\BP_r, \\
\CY_r & =
(Y\times\BP_r)\times_{\PP(V^*)\times\BP_r}\PP_{\BP_r}(\CL_r) &
\subset Y\times\BP_r,
\end{array}
$$
It is clear that both $\CX_r$ and $\CY_r$ are fibred over $\BP_r = \Gr(r,V^*)$
with fibers $X_L$ and $Y_L$ over a point corresponding to a subspace $L \subset V^*$.

Consider the fiber product $\CX_r\times_{\BP_r}\CY_r$ and the projection
$\pi_r:\CX_r\times_{\BP_r}\CY_r \to X\times Y$. Since for any vector subspace
$L\subset V^*$ the product $\PP(L^\perp)\times\PP(L)$ is contained
in the incidence quadric $Q \subset \PP(V)\times\PP(V^*)$ it follows
that $\pi_r$ factors via a map
$\zeta_r:\CX_r\times_{\BP_r}\CY_r \to Q(X,Y) \subset X\times Y$.

Consider the object $\CE_r = \zeta_r^*\CE \in \D^b(\CX_r\times_{\BP_r}\CY_r)$
as a kernel on $\CX_r\times\CY_r$. It gives the following kernel functors
$\Phi_r = \Phi_{\CE_r}:\D^b(\CY_r) \to \D^b(\CX_r)$ and
$\Phi_r^! = \Phi_{\CE_r}^!:\D^b(\CX_r) \to \D^b(\CY_r)$.
We will check below that the functors $\Phi_r$ are splitting
for all $r$ and that there exist the following semiorthogonal
decompositions
$$
\arraycolsep = 2pt
\begin{array}{lll}
\D^b(\CX_r) &=& \langle \CC_r,
\CA_r(1)\boxtimes\D^b(\BP_r),\dots,\CA_{\ix-1}(\ix-r)\boxtimes\D^b(\BP_r)\rangle\\
\D^b(\CY_r) &=& \langle \CB_{\jx-1}(N-r-\jx)\boxtimes\D^b(\BP_r),\dots,
\CB_{N-r}(-1)\boxtimes\D^b(\BP_r),\CC_r\rangle.
\end{array}
$$
where $\CC_r = \Im\Phi_r$. After that we deduce from this
the main theorem~\ref{themain} using the faithful base change
theorem~\ref{phitsod}.

For the proof of the above decompositions we use induction in $r$.
Note that for $r=1$ we have $\BP_1 = \PP(V^*)$, $\CY_1 = Y$,
$\CX_1$ is the universal hyperplane section of $X$ and $\CE_1 = \CE$,
hence the base of induction is given by the definition of
Homological Projective Duality ($\CC_1 = \CC = \D^b(Y)$).

To compare the universal families $\CX_{r-1}$, $\CY_{r-1}$ and
$\CX_r$, $\CY_r$ we take
for a base scheme
$$
\BS_r = \Fl(r-1,r;V^*) \subset \Gr(r-1,V^*)\times\Gr(r,V^*) = \BP_{r-1}\times\BP_r,
$$
the partial flag variety. The scheme $\BS_r$ parameterizes
flags $L_{r-1}\subset L_r\subset V^*$ such that $\dim L_{r-1} = r-1$,
and $\dim L_r = r$.
Let $\phi:\BS_r\to\BP_{r-1}$ and
$\psi:\BS_r\to\BP_r$ denote the natural projections.
Let $\TCL_{r-1} = \phi^*\CL_{r-1}$, $\TCL_{r} = \psi^*\CL_{r}$,
$\TCL^\perp_{r-1} = \phi^*\CL_{r-1}^\perp$,
$\TCL^\perp_{r} = \psi^*\CL^\perp_{r}$.
Then we have the universal flags of subbundles
\begin{equation}\label{tclrss}
\TCL_{r-1} \subset \TCL_r \subset V^*\otimes\CO_{\BS_r}, \qquad
\TCL_r^\perp \subset \TCL_{r-1}^\perp \subset V\otimes\CO_{\BS_r},
\end{equation}
Denote
$$
\arraycolsep=2pt
\begin{array}{llllll}
\TCX_{r-1} & =
\CX_{r-1}\times_{\BP_{r-1}}\BS_r  &
\subset X\times\BS_r, \qquad\qquad &
\TCX_r & =
\CX_r\times_{\BP_r}\BS_r &
\subset X\times\BS_r, \smallskip\\
\TCY_{r-1} & =
\CY_{r-1}\times_{\BP_{r-1}}\BS_r  &
\subset Y\times\BS_r, &
\TCY_r & =
\CY_r\times_{\BP_r}\BS_r  &
\subset Y\times\BS_r.
\end{array}
$$
Note that
$$
\arraycolsep=2pt
\begin{array}{llll}
\TCX_{r-1} &
= (X\times\BS_r) \times_{\PP(V)\times\BS_r}\PP_{\BS_r}(\TCL_{r-1}^\perp),
\qquad\qquad &
\TCX_r &
= (X\times\BS_r) \times_{\PP(V)\times\BS_r}\PP_{\BS_r}(\TCL_{r}^\perp),
\smallskip\\
\TCY_{r-1} &
= (Y\times\BS_r) \times_{\PP(V^*)\times\BS_r}\PP_{\BS_r}(\TCL_{r-1}),
&
\TCY_r &
= (Y\times\BS_r) \times_{\PP(V^*)\times\BS_r}\PP_{\BS_r}(\TCL_{r}).
\end{array}
$$
Therefore the embeddings~(\ref{tclrss}) induce embeddings
$\xi:\TCX_r \to \TCX_{r-1}$ and $\eta:\TCY_{r-1} \to \TCY_r$.
Consider the following commutative diagrams
(the squares marked with $\ecart$ are exact cartesian by lemma~\ref{ec},
because the maps $\phi$ and $\psi$ are flat)
\begin{equation}\label{yayb}
\vcenter{\xymatrix{
\BS_r \ar[d]_\psi \ar@{}[dr]|{\ecart} &
\TCY_r \ar[l]_{g_r} \ar[d]_\psi &
\TCY_{r-1} \ar[d]^\phi \ar[l]_\eta \ar[dl]_{\hat{\psi}} \ar[r]^{g_{r-1}} &
\BS_r \ar[d]_\phi \ar@{}[dl]|{\ecart} \\
\BP_r &
\CY_r \ar[l]_{g_r} &
\CY_{r-1} \ar[r]^{g_{r-1}}  &
\BP_{r-1}
}}
\qquad\qquad
\vcenter{\xymatrix{
\BS_r \ar[d]_\psi \ar@{}[dr]|{\ecart} &
\TCX_r \ar[l]_{f_r} \ar[d]_\psi \ar[r]^\xi \ar[dr]^{\hat{\phi}} &
\TCX_{r-1} \ar[d]^\phi \ar[r]^{f_{r-1}} &
\BS_r \ar[d]_\phi \ar@{}[dl]|{\ecart}  \\
\BP_r &
\CX_r \ar[l]_{f_r} &
\CX_{r-1} \ar[r]^{f_{r-1}}  &
\BP_{r-1}
}}
\end{equation}
where $f_{r-1}$, $f_r$, $g_{r-1}$ and $g_r$ are the natural projections
and $\hat{\phi} = \phi\circ\xi$, $\hat{\psi} = \psi\circ\eta$.

Let $\TE_{r-1}\in\D(\TCX_{r-1}\times_{\BS_r}\TCY_{r-1})$ and
$\TE_r \in \D(\TCX_r\times_{\BS_r}\TCY_r)$ denote the pullbacks
of the objects $\CE_{r-1}$ and $\CE_r$ via the projections
$\TCX_{r-1}\times_{\BS_r}\TCY_{r-1} \to \CX_{r-1}\times_{\BP_{r-1}}\CY_{r-1}$,
$\TCX_r\times_{\BS_r}\TCY_r \to \CX_r\times_{\BP_r}\CY_r$.
Then we have the corresponding kernel functors
$\TPhi_{r-1}$, $\TPhi_r$ e.t.c between the derived categories
of $\TCX_{r-1}$, $\TCY_{r-1}$, $\TCX_r$ and $\TCY_r$.

The induction step is based on relation of the functors
$\Phi_{r-1}$, $\Phi_r$, $\TPhi_{r-1}$ and $\TPhi_r$
to the base change functors $\psi^*$, $\psi_*$, $\phi^*$, $\phi_*$
and to the functors of the pushforward and pullback
via $\xi$ and $\eta$. The relation to $\psi^*$, $\psi_*$, $\phi^*$ and $\phi_*$
is given by lemma~\ref{phit}. The relation to $\xi$ and $\eta$
in a sense is the key point of the proof.
We prove that $\xi^*\TPhi_{r-1} \cong \TPhi_r\eta_*$
and that the ``difference'' between $\xi_*\TPhi_r$ and $\TPhi_{r-1}\eta^!$
is given by a very simple functor.

Other results in this section (e.g.\ the above semiorthogonal
decompositions) are proved by similar arguments using (either
ascending or descending) induction in $r$.

The section is organized as follows. We start with some preparations
concluding with a description of the relation of the functors
$\TPhi_{r-1}$ and $\TPhi_r$ to the pushforward and pullback
via $\xi$ and $\eta$. Then we use induction in $r$ to prove
the semiorthogonal decompositions.

\subsection{Preparations}\label{indsetup}

\begin{lemma}\label{ysm}
$Y$ is smooth.
\end{lemma}
\begin{proof}
By definition of Homological Projective Duality $\D^b(Y)$ is
a full subcategory of $\D^b(\CX_1)$. On the other hand,
in the following lemma we prove that $\CX_1$ is smooth,
hence its derived category $\D^b(\CX_1)$ is $\Ext$-bounded
by lemma~\ref{snav_sm}. Therefore $\D^b(Y)$ is also $\Ext$-bounded,
hence $Y$ is smooth again by lemma~\ref{snav_sm}.
\end{proof}

Recall that $\BP_r = \Gr(r,V^*)$ is the Grassmannian
parameterizing linear sections of $X$ and $Y$, and $\CX_r$, $\CY_r$
are the universal families over $\BP_r$ of linear sections.

\begin{lemma}\label{xyksm}
The projections $\CX_r \to X$, $\CY_r \to Y$ and
$\zeta_r:\CX_r\times_{\BP_r}\CY_r \to Q(X,Y)$ are smooth.
In particular, $\CX_r$ and $\CY_r$ are smooth and
$$
\begin{array}{l}
\dim\CX_r = \dim X + \dim\BP_r - r,\qquad
\dim\CY_r = \dim Y + \dim\BP_r + r - N,\\
\dim\CX_r\times_{\BP_r}\CY_r = \dim X + \dim Y + \dim\BP_r - N.
\end{array}
$$
Moreover, the maps $f_r:\CX_r \to \BP_r$ and $g_r:\CY_r \to \BP_r$ are projective.
\end{lemma}
\begin{proof}
Note that we have the following isomorphisms
$$
\CX_r = \Gr_X(r,\CV_X),\qquad
\CY_r = \Gr_Y(r-1,\CV_Y),\qquad
\CX_r\times_{\BP_r}\CY_r = \Gr_{Q(X,Y)}(r-1,\CV_Q)
$$
with the relative Grassmannians, where the bundles $\CV_X$, $\CV_Y$
are defined from exact sequences
$$
0 \to \CV_X \to V^*\otimes\CO_X \to \CO_X(1) \to 0,\qquad
0 \to \CO_Y(-1) \to V^*\otimes\CO_Y \to \CV_Y \to 0,
$$
and $\CV_Q$ is the middle cohomology bundle of the complex
$$
\CO_{Q(X,Y)}(0,-1) \to V^*\otimes\CO_{Q(X,Y)} \to \CO_{Q(X,Y)}(1,0).
$$
From this and lemma~\ref{ysm} we easily deduce the smoothness and
compute the dimensions. It is also clear that the fibers
of the projections $\CX_r \to \BP_r$, and $\CY_r \to \BP_r$ are
linear sections of $X$ and $Y$ corresponding to subspaces $L \in \BP_r$,
so they are projective.
\end{proof}

\begin{lemma}\label{zl1}
In parts $(i)$ and $(ii)$ below $k$ stands either for $r$ or for $r-1$.

\noindent
$(i)$ $\TCX_k$ is the zero locus of a section
of vector bundle $\CO_X(1)\boxtimes\TCL_k^*$ on $X\times\BS_r$;

\noindent
$(ii)$
$\TCY_k$ is the zero locus of a section
of vector bundle $\CO_Y(1)\boxtimes\TCL_k^{\perp*}$ on $Y\times\BS_r$;

\noindent
$(iii)$
$\TCX_r$ is the zero locus of a section of line bundle
$\CO_X(1)\otimes(\TCL_r/\TCL_{r-1})^*$ on $\TCX_{r-1}$;

\noindent
$(iv)$
$\TCY_{r-1}$ is the zero locus of a section of line bundle
$\CO_Y(1)\otimes(\TCL_{r-1}^\perp/\TCL_r^\perp)^* \cong
\CO_Y(1)\otimes(\TCL_r/\TCL_{r-1})$ on $\TCY_r$.\\
All these sections are regular.
\end{lemma}
\begin{proof}
The parts $(i)$ and $(ii)$ evidently follow from the definition
of $\TCX_k \subset X\times\BS_r$. The parts $(iii)$ and $(iv)$
follow from the exact sequences
$$
0 \to (\TCL_r/\TCL_{r-1})^* \to \TCL_r^* \to \TCL_{r-1}^* \to 0,
\qquad
0 \to (\TCL_{r-1}^\perp/\TCL_r^\perp)^* \to
(\TCL_{r-1}^\perp)^* \to (\TCL_r^\perp)^* \to 0.
$$
Finally, it follows from lemma~\ref{xyksm} that
$\dim\TCX_{r-1} = \dim X + \dim\BS_r - (r-1)$,
$\dim\TCX_r = \dim X + \dim\BS_r - r$,
since the base changes $\BS_r \to \BP_{r-1}$ and $\BS_r \to \BP_r$
are flat. Therefore the sections in the parts~$(i)$ and $(iii)$ are regular.
The sections in the parts~$(ii)$ and $(iv)$ are regular by similar reasons.
\end{proof}

\begin{lemma}\label{soc_xyr}
For any $r\le k\le \ix-1$ the functors
$\CA_k\boxtimes\D^b(\BP_r) \subset \D^b(X\times\BP_r) \exto{i^*} \D^b(\CX_r)$
are fully faithful, and the collection
$(\CA_r(1)\boxtimes\D^b(\BP_r),\dots,\CA_{\ix-1}(\ix-r)\boxtimes\D^b(\BP_r)) \subset \D^b(\CX_r)$
is semiorthogonal.

Similarly, for any $N-r\le k\le \jx-1$ the functors
$\CB_k\boxtimes\D^b(\BP_r) \subset \D^b(Y\times\BP_r) \exto{i^*} \D^b(\CY_r)$
are fully faithful, and the collection
$(\CB_{\jx-1}(N-r-\jx)\boxtimes\D^b(\BP_r),\dots,\CB_{N-r}(-1)\boxtimes\D^b(\BP_r)) \subset \D^b(\CY_r)$
is semiorthogonal.
\end{lemma}
\begin{proof}
Analogous to the proof of lemma~\ref{sodx1}
using the Koszul resolutions of $\CO_{\CX_r}$ on $X\times\BP_r$
and of $\CO_{\CY_r}$ on $Y\times\BP_r$.
\end{proof}

Now we describe the maps $\psi$ and $\phi$.

\begin{lemma}\label{im1}
The maps $\phi$, $\psi$, $\hat{\phi}$ and $\hat{\psi}$ are projectivizations
of vector bundles. Explicitly,
$\phi$ is the projectivization of $V^*\otimes\CO/\CL_r$,
$\psi$ is the projectivization of $\CL_r^*$,
$\hat\phi$ is the projectivization of $\CV_X/\CL_{r-1}$, and
$\hat\psi$ is the projectivization of $(\CL_r/\CO_{\PP(V^*)}(-1))^*$,
where the embedding $\CO_{\PP(V^*)}(-1) \to \CL_r$ is induced
by the projection $\CY_r \to \PP_{\BP_r}(\CL_r) \to \PP(V^*)$.
\end{lemma}
\begin{proof}
By definition of $\BS_r$ the fiber of $\phi$ is the set of all
lines in $V/L_{r-1}$ and the fiber of $\psi$ is the set of all
hyperplanes in $L_r$. Similarly, the fiber of $\hat\phi$ is the set
of all lines in $V/L_{r-1}$ contained in $\CV_x$, where $\CV_x$ is the fiber at $x$
of the vector bundle $\CV_X$ on $X$ defined in the proof of lemma~\ref{xyksm},
and the fiber of $\hat\psi$ is the set of all hyperplanes in $L_r$ passing
through a point $y\in\PP(L_r)$.
\end{proof}

Applying results of~\cite{O2} we deduce the following.

\begin{corollary}\label{cim1}
The functors $\phi^*$, $\psi^*$, ${\hat\phi}^*$ and ${\hat\psi}^*$
are fully faithful and we have
$$
\phi_*\phi^* \cong \id,\qquad
\psi_*\psi^* \cong \id,\qquad
{\hat\phi}_*{\hat\phi}^* \cong \id,\qquad
{\hat\psi}_*{\hat\psi}^* \cong \id.
$$
\end{corollary}

Recall that we have defined the objects $\CE_r$ on
$\CX_r\times_{\BP_r}\CY_r$ as the pullbacks of $\CE\in\D^b(Q(X,Y))$
via the map $\CX_r\times_{\BP_r}\CY_r \to Q(X,Y)$,
and the objects $\TE_{r-1}$ and $\TE_r$ as the pullbacks
of $\CE_{r-1}$ and $\CE_r$ via the maps
$\phi:\TCX_{r-1}\times_{\BS_r}\TCY_{r-1} \to
\CX_{r-1}\times_{\BP_{r-1}}\CY_{r-1}$ and
$\psi:\TCX_r\times_{\BS_r}\TCY_r \to \CX_r\times_{\BP_r}\CY_r$.
The functors $\Phi_{r-1}$, $\Phi_r$, $\Phi_{r-1}^!$, $\Phi_r^!$,
$\TPhi_{r-1}$, $\TPhi_r$, $\TPhi_{r-1}^!$, and $\TPhi_r^!$ are
kernel functors of the first and second type corresponding
to the kernels $\CE_{r-1}$, $\CE_r$, $\TE_{r-1}$ and $\TE_r$ respectively.

\begin{lemma}\label{b_adj}
The functors
$\Phi_{r-1}$, $\Phi_r$, $\TPhi_{r-1}$ and $\TPhi_r$
have right adjoint functors
$\Phi_{r-1}^!$, $\Phi_r^!$, $\TPhi_{r-1}^!$, and $\TPhi_r^!$,
and left adjoint functors
$\Phi_{r-1}^*$, $\Phi_r^*$, $\TPhi_{r-1}^*$ and $\TPhi_r^*$.
All these functors take bounded derived categories to
bounded derived categories.
%
%
%
\end{lemma}
\begin{proof}
First of all note that $\CE_r = \zeta_r^*\CE$ is bounded because
$\zeta_r:\CX_r\times\BP_r\CY_r \to Q(X,Y)$ is smooth by lemma~\ref{xyksm}.
Since $\CX_r$ and $\CY_r$ are smooth by lemma~\ref{xyksm} it follows that
the pushforward of $\CE_r$ to $\CX_r\times\CY_r$ is a perfect complex.
In particular it has finite $\Tor$ and $\Ext$-amplitude over $\CX_r$ and $\CY_r$.
On the other hand, the projections of $\CX_r\times_{\BP_r}\CY_r$ to the factors
are projective because the projections of $\CX_r$ and $\CY_r$ to $\BP_r$
are projective by lemma~\ref{xyksm}. Therefore by lemma~\ref{phi_bounded}
the functor $\Phi_r^!$ is right adjoint to $\Phi_r$ and by lemma~\ref{ladj}
there exists a left adjoint functor $\Phi_r^*$ to $\Phi_r$.
Moreover, all these functors take bounded derived categories
to bounded derived categories.
The same arguments prove the rest of the lemma.
%
%
\end{proof}

\begin{proposition}\label{phiPhi}
We have functorial isomorphisms
$$
\phi^*\Phi_{r-1} \cong \TPhi_{r-1}\phi^*,\quad
\psi^*\Phi_{r} \cong \TPhi_{r}\psi^*,\qquad
\phi_*\TPhi_{r-1} \cong \Phi_{r-1}\phi_*,\quad
\psi_*\TPhi_{r} \cong \Phi_{r}\psi_*.
$$
\end{proposition}
\begin{proof}
Note that the base changes $\phi$ and $\psi$ are smooth, hence they are
faithful by lemma~\ref{fisf} and we conclude by lemma~\ref{phit}.
\end{proof}

Now we go to the relation of the functors $\TPhi_{r-1}$ and $\TPhi_r$
to the pushforward and pullback via $\xi$ and $\eta$.

Consider the following diagram
$$
\xymatrix{
\TCX_r\times_{\BS_r}\TCY_{r-1} \ar[r]^{\eta} \ar[d]_{\xi} &
\TCX_r\times_{\BS_r}\TCY_r \ar[d]_{\xi} \\
\TCX_{r-1}\times_{\BS_r}\TCY_{r-1} \ar[r]^{\eta} &
\TCX_{r-1}\times_{\BS_r}\TCY_r \ar[r]^-{\tpi} &
X\times Y
}
$$
where $\tpi$ is the composition
$\TCX_{r-1}\times_{\BS_r}\TCY_r \subset
\TCX_{r-1}\times\TCY_r \subset
(X\times\BS_r)\times(Y\times\BS_r) \to
X\times Y$.

\begin{lemma}
The maps $\xi$ and $\eta$ in the above diagram are divisorial embeddings and
\begin{equation}\label{xietashriek}
\xi^!(-) = \xi^*(-)\otimes\CO(1,0)\otimes(\TCL_r/\TCL_{r-1})^*[-1],
\quad
\eta^!(-) = \eta^*(-)\otimes\CO(0,1)\otimes(\TCL_r/\TCL_{r-1})[-1].
\end{equation}
Moreover, we have the following scheme-theoretical equalities
\begin{equation}\label{capcup}
\begin{array}{l}
(\TCX_{r-1}\times_{\BS_r}\TCY_{r-1}) \cap (\TCX_r\times_{\BS_r}\TCY_r) =
\TCX_r\times_{\BS_r}\TCY_{r-1},\\
(\TCX_{r-1}\times_{\BS_r}\TCY_{r-1}) \cup (\TCX_r\times_{\BS_r}\TCY_r) =
\tpi^{-1}(Q(X,Y)).
\end{array}
\end{equation}
and the following square is exact carthesian
\begin{equation}\label{xaybs}
\vcenter{\xymatrix{
(\TCX_{r-1}\times_{\BS_r}\TCY_{r-1}) \cup (\TCX_r\times_{\BS_r}\TCY_r)
\ar[rr]^-i \ar[d]^{\tzeta} &&
\TCX_{r-1}\times_{\BS_r}\TCY_r \ar[d]^{\tpi} \\
Q(X,Y) \ar[rr]^i &&
X\times Y \ar@{}[ull]|{\ecart}
}}
\end{equation}
\end{lemma}
\begin{proof}
Consider the projections of
$\TCX_r\times_{\BS_r}\TCY_{r-1}$,
$\TCX_{r-1}\times_{\BS_r}\TCY_{r-1}$,
$\TCX_r\times_{\BS_r}\TCY_r$ and
$\TCX_{r-1}\times_{\BS_r}\TCY_r$
to $X\times Y$. It is easy to check that
their fibers over a point $(x,y) \in (X,Y)$ are subsets
of the flag variety $\Fl(r-1,r;V^*) = \BS_r$ consisting of all flags
$L_{r-1} \subset L_r$ satisfying the following incidence conditions
$$
\arraycolsep=2pt
\begin{array}{ccccc}
y & \subset & L_{r-1} \\
&& \cap \\
&& L_r & \subset & \CV_x
\end{array}
\qquad,\qquad
\begin{array}{ccccc}
y & \subset & L_{r-1} & \subset & \CV_x \\
&& \cap \\
&& L_r
\end{array}
\qquad,\qquad
\begin{array}{ccccc}
&& L_{r-1} \\
&& \cap \\
y & \subset & L_r & \subset & \CV_x
\end{array}
\qquad\text{and}\qquad
\begin{array}{ccccc}
&& L_{r-1} & \subset & \CV_x \\
&& \cap \\
y & \subset & L_r
\end{array}
$$
respectively. In particular, the first three fibers are empty
if $(x,y)\not\in Q(X,Y)$. On the other hand, over $Q(X,Y)$
the first three fibers are irreducible, have dimension
$(r-1)(N-r)-1$, $(r-1)(N-r)$ and $(r-1)(N-r)$ respectively,
and the first of them is the intersection of the other two.
On the contrary, the fourth fiber is irreducible and $(r-1)(N-r)$-dimensional
if $(x,y)\not\in Q(X,Y)$ and for $(x,y)\in Q(X,Y)$ it coincides with the union
of the second and the third fibers (if $y\subset \CV_x$ and
$y\not\subset L_{r-1}$ then $L_r = \langle y, L_{r-1}\rangle \subset \CV_x$).
It follows that images of $\xi$ and $\eta$ have pure codimension~$1$.
Since they are also zero loci of line bundles by lemma~\ref{zl1}~$(iii)$
and $(iv)$, we conclude that $\xi$ and $\eta$ are divisorial embeddings.

The above arguments also prove the first equality of~(\ref{capcup})
on the scheme-theoretical level and the second equality on the set-theoretical
level. Taking into account that the LHS of the second equality is the zero
locus of the line bundle $\CO_X(1)\otimes\CO_Y(1)$ by definition of $Q(X,Y)$,
and that the RHS of the equality is the zero locus of the line bundle
$(\CO_X(1)\otimes(\TCL_r/\TCL_{r-1})^*)\otimes
(\CO_Y(1)\otimes(\TCL_r/\TCL_{r-1}))$
by lemma~\ref{zl1}~$(iii)$ and $(iv)$, and noting that these bundles
are isomorphic, we deduce that the second equality is also true
on the scheme-theoretical level.
Finally, we note that the square~(\ref{xaybs}) is exact cartesian
by lemma~\ref{ec}~$(iii)$ since $X\times Y$ and $Q(X,Y)$ are Cohen-Macaulay.
\end{proof}

Consider the pullback $\HE = \tzeta^*\CE$ of $\CE$ from $Q(X,Y)$ to
$(\TCX_{r-1}\times_{\BS_r}\TCY_{r-1}) \cup (\TCX_r\times_{\BS_r}\TCY_r)$.
Let us also denote $\CO(k,l) := \CO_X(k)\otimes\CO_Y(l)$ for brevity.
The following lemma gives a relation of $\TE_{r-1}$ and $\TE_r$.

\begin{lemma}
We have the following exact triangle on
$\TCX_{r-1}\times_{\BS_r}\TCY_r$:
\begin{equation}\label{ieab2}
\xi_*\TE_r(0,-1)\otimes(\TCL_r/\TCL_{r-1})^* \to
i_*\HE \to
\eta_*\TE_{r-1}.
\end{equation}
Moreover, we have an isomorphism on $\TCX_r\times_{\BS_r}\TCY_{r-1}$:
\begin{equation}\label{eaeb}
\eta^*\TE_r \cong \xi^*\TE_{r-1}.
\end{equation}
and an isomorphism on $\TCX_{r-1}\times_{\BS_r}\TCY_r$
\begin{equation}\label{ieab}
i_*\HE  \cong  \tpi^*i_*\CE.
\end{equation}
\end{lemma}
\begin{proof}
Since the square~(\ref{xaybs}) is exact cartesian we have
$i_*\HE = i_*\tzeta^*\CE = \tpi^*i_*\CE$ which gives us~(\ref{ieab}).
Triangle (\ref{ieab2}) can be obtained by tensoring the resolution
(the twist of the left term is determined by lemma~\ref{zl1}~$(iv)$)
$$
0 \to \CO_{\TCX_r\times_{\BS_r}\TCY_r}(0,-1)\otimes(\TCL_r/\TCL_{r-1})^*
\to
\CO_{(\TCX_{r-1}\times_{\BS_r}\TCY_{r-1}) \cup (\TCX_r\times_{\BS_r}\TCY_r)}
\to
\CO_{\TCX_{r-1}\times_{\BS_r}\TCY_{r-1}} \to 0
$$
with $\HE$ and applying $i_*$, since the pullback of $\HE$ to
$\TCX_{r-1}\times_{\BS_r}\TCY_{r-1}$ and $\TCX_r\times_{\BS_r}\TCY_r$
coincides with $\TE_{r-1}$ and $\TE_r$ respectively.
Finally, (\ref{eaeb}) is evident, because both sides are isomorphic
to the pullback of $\HE$.
\end{proof}

\begin{corollary}\label{phixieta1}
We have the following exact triangles of functors between
$\D^b(\TCX_{r-1})$ and $\D^b(\TCY_r)$:
\begin{equation}\label{ft1}
\TPhi_{r-1}\eta^! \to \xi_*\TPhi_r \to
\Phi_{i_*\HE(0,1)\otimes(\TCL_r/\TCL_{r-1})},
\end{equation}
\begin{equation}\label{ft3}
\Phi^*_{i_*\HE(0,1)\otimes(\TCL_r/\TCL_{r-1})} \to
\TPhi_r^*\xi^* \to \eta_*\TPhi_{r-1}^*
\end{equation}
and the following canonical isomorphism of functors
from $\D^b(\TCY_{r-1})$ to $\D^b(\TCX_r)$:
\begin{equation}\label{fi1}
\xi^*\TPhi_{r-1} \cong \TPhi_r\eta_*.
\end{equation}
\end{corollary}
\begin{proof}
Twisting triangle (\ref{ieab2}) by $\CO(0,1)\otimes(\TCL_r/\TCL_{r-1})$,
considering its terms as kernels, and taking into account the second formula
of~(\ref{xietashriek}) we obtain triangles~(\ref{ft1}) and (\ref{ft3})
by lemma~\ref{etf}.
Finally, isomorphism of kernels~(\ref{eaeb}) gives an isomorphism
of functors~(\ref{fi1}).
%
%
\end{proof}

\begin{lemma}\label{zl2}
$(i)$ The map $\CX_r\times_{\BP_r}\CY_r \to \CY_r\times_Y Q(X,Y)$ induced
by the projection $\CX_r\times_{\BP_r}\CY_r \to \CY_r$ and by the map
$\zeta_r:\CX_r\times_{\BP_r}\CY_r \to Q(X,Y)$ is a closed embedding
and its image is a zero locus of a regular section of the vector bundle
$(\CL_r/\CO_Y(-1))^*\otimes\CO_X(1)$ on $\CY_r\times_Y Q(X,Y)$.

$(ii)$ The map $\CX_r\times_{\BP_r}\CY_r \to \CX_r\times_X Q(X,Y)$ induced
by the projection $\CX_r\times_{\BP_r}\CY_r \to \CX_r$ and by the map
$\zeta_r:\CX_r\times_{\BP_r}\CY_r \to Q(X,Y)$ is a closed embedding
and its image is a zero locus of a regular section of the vector bundle
$(\CL_r^\perp/\CO_X(-1))^*\otimes\CO_Y(1)$ on $\CX_r\times_X Q(X,Y)$.
\end{lemma}
\begin{proof}
Recall the notation of the proof of lemma~\ref{xyksm}. It is clear that we have
$$
\begin{array}{c}
\CX_r\times_{\BP_r}\CY_r =
\Gr_{Q(X,Y)}(r-1,\CV_Q),
\smallskip\\
\CY_r\times_Y Q(X,Y) =
\Gr_{Q(X,Y)}(r-1,\CV_Y),
\quad\text{and}\quad
\CX_r\times_X Q(X,Y) =
\Gr_{Q(X,Y)}(r,\CV_X).
\end{array}
$$
Moreover, $\CV_Q = \Ker(\CV_Y \to \CO_X(1)) = \Ker(\CV_X \to \CO_Y(1))$.
These considerations make the claims of the lemma evident.
\end{proof}

\begin{lemma}\label{imkerphir}
We have

\noindent$(i)$
$\Im \Phi_r \subset
[\langle\CA_{r}(1),\dots,\CA_{\ix-1}(\ix-r)\rangle\boxtimes\D^b(\BP_r)]^\perp
\subset \D^b(\CX_{r})$;

\noindent$(ii)$
$\langle\CB_{\jx-1}(N-r-\jx),\dots,\CB_{N-r}(-1)\rangle\boxtimes\D^b(\BP_r)
\subset \Ker \Phi_r \subset \D^b(\CY_r)$.
\end{lemma}
\begin{proof}
$(i)$ Consider the following commutative diagram
$$
\xymatrix{
\CY_r \ar@{=}[d] && \CX_r\times_{\BP_r}\CY_r \ar[ll]_q \ar[rr]^p \ar[d]^j && \CX_r \ar[dd]^\pi \\
\CY_r \ar[d]_\pi \ar@{}[rrd]|\ecart && \CY_r\times_Y Q(X,Y) \ar[d]^\pi \ar[ll]_q \\
Y && Q(X,Y) \ar[ll]_q \ar[rr]^p && X
}
$$
(the square marked with $\ecart$ is exact cartesian by lemma~\ref{ec} since the map
$\pi:\CY_r \to Y$ is smooth by lemma~\ref{xyksm}).
It is clear that the functor $\pi_*\circ\Phi_r:\D^b(\CY_r) \to \D^b(X)$
is a kernel functor with kernel $j_*\CE_r$ on $\CY_r\times_Y Q(X,Y)$.
Further, by definition $\CE_r = (\pi\circ j)^*\CE \cong j^*\pi^*\CE$,
hence
$j_*\CE_r \cong j_*j^*\pi^*\CE \cong
\pi^*\CE\otimes j_*\CO_{\CX_r\times_{\BP_r}\CY_r}$.
On the other hand, by lemma~\ref{zl2} we have a Koszul resolution
\begin{multline*}
0 \to
\Lambda^{r-1}(\CL_r/\CO_Y(-1))\otimes\CO_X(1-r) \to \dots \to
\Lambda^2(\CL_r/\CO_Y(-1))\otimes\CO_X(-2) \to
\\ \to
(\CL_r/\CO_Y(-1))\otimes\CO_X(-1) \to
\CO_{\CY_r\times_Y Q(X,Y)} \to
j_*\CO_{\CX_r\times_{\BP_r}\CY_r} \to 0.
\end{multline*}
Considering
$\pi^*\CE\otimes\Lambda^t(\CL_r/\CO_Y(-1))\otimes\CO_X(-t) \cong
\pi^*(\CE\otimes\CO_X(-t))\otimes q^*\Lambda^t(\CL_r/\CO_Y(-1))$
as a kernel on $\CY_r\times_Y Q(X,Y)$ we note that the corresponding
kernel functor $\Psi_t:\D^b(\CY_r) \to \D^b(X)$ takes form
\begin{multline*}
G \mapsto
p_*\pi_*(\pi^*(\CE\otimes\CO_X(-t))\otimes q^*\Lambda^t(\CL_r/\CO_Y(-1))\otimes q^* G) \cong
\\ \cong
p_*(\CE\otimes\CO_X(-t))\otimes \pi_*q^*(\Lambda^t(\CL_r/\CO_Y(-1))\otimes G) \cong
p_*(\CE\otimes q^*\pi_*(\Lambda^t(\CL_r/\CO_Y(-1))\otimes G)\otimes\CO_X(-t) =
\\ \cong
\Phi_\CE(\pi_*(\Lambda^t(\CL_r/\CO_Y(-1))\otimes G))\otimes\CO_X(-t).
\end{multline*}
(in the second isomorphism we used exactness of the square marked with $\ecart$ symbol).
Note that $\Phi_\CE(\pi_*(\Lambda^t(\CL_r/\CO_Y(-1))\otimes G)) \in \CA_0$
for any $G\in\D^b(\CY_r)$ by lemma~\ref{impig},
hence the image of the functor $\Psi_t$ is contained in the subcategory $\CA_0(-t) \subset \D^b(X)$.
It follows that $\Im (\pi_*\Phi_r) \subset \lan\CA_0(1-r),\dots,\CA_0(-1),\CA_0\ran$.
But the latter subcategory of $\D^b(X)$ coincides with
$\lan \CA_0(1-r),\dots,\CA_{r-1}\ran = \lan\CA_{r}(1),\dots,\CA_{\ix-1}(\ix-r)\ran^\perp$ by lemma~\ref{a0rgen}.
Therefore we have $\Im(\pi_*\Phi_r) \subset \lan\CA_{r}(1),\dots,\CA_{\ix-1}(\ix-r)\ran^\perp$.
Since the functor $\Phi_r$ is $\BP_r$-linear it follows that
$\Im(\Phi_r) \subset [\lan\CA_{r}(1),\dots,\CA_{\ix-1}(\ix-r)\ran\boxtimes\D^b(\BP_r)]^\perp$.

$(ii)$
Similarly, consider the following commutative diagram
$$
\xymatrix{
\CY_r \ar[dd]_\pi && \CX_r\times_{\BP_r}\CY_r \ar[ll]_q \ar[rr]^p \ar[d]^j && \CX_r \ar@{=}[d]  \\
&& \CX_r\times_X Q(X,Y) \ar[d]^\pi \ar[rr]^p && \CX_r  \ar[d]^\pi \ar@{}[lld]|\ecart \\
Y && Q(X,Y) \ar[ll]_q \ar[rr]^p && X
}
$$
(the square marked with $\ecart$ is exact cartesian by lemma~\ref{ec} since the map
$\pi:\CX_r \to X$ is smooth by lemma~\ref{xyksm}).
and functor $\Phi_r\circ\pi^*:\D^b(Y) \to \D^b(\CX_r)$ which is
a kernel functor with kernel
$\pi^*\CE\otimes j_*\CO_{\CX_r\times_{\BP_r}\CY_r}$
on $\CX_r\times_X Q(X,Y)$.
On the other hand, by lemma~\ref{zl2} we have a Kozsul resolution
\begin{multline*}
0 \to
\Lambda^{N-r-1}(\CL_r^\perp/\CO_X(-1))\otimes\CO_Y(1+r-N) \to \dots \to
\Lambda^2(\CL_r^\perp/\CO_X(-1))\otimes\CO_Y(-2) \to
\\ \to
(\CL_r^\perp/\CO_X(-1))\otimes\CO_Y(-1) \to
\CO_{\CX_r\times_X Q(X,Y)} \to
j_*\CO_{\CX_r\times_{\BP_r}\CY_r} \to 0.
\end{multline*}
The kernel functors $\Psi_t:\D^b(Y) \to \D^b(X)$ given by the terms of this resolution
tensored by $\pi^*\CE$ take form
$$
\Psi_t:G \mapsto
\pi^*\Phi_\CE(G(-t)) \otimes \Lambda^t(\CL^\perp_r/\CO_X(-1)),
\qquad
t=0,\dots,N-r-1.
$$
If $G \in \lan\CB_{\jx-1}(N-r-\jx),\dots,\CB_{N-r}(-1)\ran$ then
for any $t=0,1,\dots,N-r-1$ we have
$$
G(-t) \in \lan\CB_{\jx-1}(1-\jx),\dots,\CB_{1}(-1)\ran
$$
But $\lan\CB_{\jx-1}(1-\jx),\dots,\CB_{1}(-1)\ran \subset \Ker\Phi_\CE$
by lemma~\ref{imkerpig}, hence
for $G \in \lan\CB_{\jx-1}(N-r-\jx),\dots,\CB_{N-r}(-1)\ran$ and $0\le t \le N-r-1$
we have $\Psi_t(G) = 0$.
Therefore
$\lan\CB_{\jx-1}(N-r-\jx),\dots,\CB_{N-r}(-1)\ran \subset \Ker(\Phi_r\circ\pi^*)$,
and the claim follows since $\Phi_r$ is $\BP_r$-linear.
\end{proof}

\begin{lemma}\label{phiisemain}
We have
$\Phi^*_{i_*\HE(0,1)}([\lan\CA_{r-1}(1),\dots,\CA_{\ix-1}(\ix-r+1)\ran\boxtimes\D^b(\BS_r)]^\perp)
\subset
\CB_{N-r}(-1)\boxtimes\D^b(\BS_r)$.
\end{lemma}
\begin{proof}
First of all we note that the claim is equivalent to
$$
\Phi^*_{i_*\HE}([\lan\CA_{r-1}(1),\dots,\CA_{\ix-1}(\ix-r+1)\ran\boxtimes\D^b(\BS_r)]^\perp)
\subset
\CB_{N-r}\boxtimes\D^b(\BS_r)
$$
which we will prove.
Consider the following commutative diagram
$$
\xymatrix{
\TCY_r \ar@{}[drr]|\ecart \ar[d]_j && \TCX_{r-1}\times_{\BS_r}\TCY_r \ar[ll]_q \ar[rr]^p \ar[d]^j && \TCX_{r-1} \ar@{=}[d]  \\
Y\times\BS_r \ar[d]_\pi && \TCX_{r-1}\times Y \ar[d]^\pi \ar[ll]_q \ar[rr]^p && \TCX_{r-1}  \ar[d]^\pi \\
Y && X\times Y \ar[ll]_q \ar[rr]^p && X \ar@{}[ull]|\ecart
}
$$
(the squares marked with $\ecart$ are exact cartesian:
the first by lemma~\ref{ec}~$(iii)$ and the second by lemma~\ref{ec}~$(i)$).
Note also that $\pi\circ j = \tpi$, hence
$i_*\HE \cong j^*\pi^*i_*\CE$ by~(\ref{ieab}). It follows that
$$
\Phi_{i_*\HE}(G) \cong
p_*(i_*\HE \otimes q^*G) \cong
p_*j_*(j^*\pi^*i_*\CE \otimes q^*G) \cong
p_*(\pi^*i_*\CE \otimes j_*q^*G) \cong
p_*(\pi^*i_*\CE \otimes q^*j_*G) =
\Phi_{\pi^*i_*\CE}(j_*G),
$$
that is $\Phi_{i_*\HE} \cong \Phi_{\pi^*i_*\CE}\circ j_*$,
whereof by adjunction we have
$\Phi^*_{i_*\HE} \cong j^*\circ\Phi^*_{\pi^*i_*\CE}$.
Therefore it suffices to check that
$$
\Phi^*_{\pi^*i_*\CE}([\lan\CA_{r-1}(1),\dots,\CA_{\ix-1}(\ix-r+1)\ran\boxtimes\D^b(\BS_r)]^\perp) \subset
\CB_{N-r}\boxtimes\D^b(\BS_r)
\qquad\text{in $\D^b(Y\times\BS_r)$}.
$$
The LHS is evidently $\BS_r$-linear, hence this is equivalent to
$$
\pi_*\Phi^*_{\pi^*i_*\CE}([\lan\CA_{r-1}(1),\dots,\CA_{\ix-1}(\ix-r+1)\ran\boxtimes\D^b(\BS_r)]^\perp) \subset \CB_{N-r}
\qquad\text{in $\D^b(Y)$}.
$$
But
$$
\Phi_{\pi^*i_*\CE}(\pi^*G) =
p_*(\pi^*i_*\CE\otimes q^*\pi^*G) \cong
p_*(\pi^*i_*\CE\otimes \pi^*q^*G) \cong
p_*\pi^*(i_*\CE\otimes q^*G) \cong
\pi^*p_*(i_*\CE\otimes q^*G) =
\pi^*\Phi_{i_*\CE}(G),
$$
that is $\Phi_{\pi^*i_*\CE}\circ \pi^* \cong \pi^*\circ\Phi_{i_*\CE}$
whereof by adjunction we deduce
$\pi_\#\circ\Phi^*_{\pi^*i_*\CE} \cong \Phi^*_{i_*\CE}\circ\pi_\#$,
where $\pi_\#$ is the left adjoint functor to $\pi^*$.
Now note that for the projection $\pi:Y\times\BS_r \to Y$ we have
$\pi_\#(G) = \pi_*(G\otimes\omega_{\BS_r}[\dim\BS_r])$.
On the other hand, for the projection $\pi:\TCX_{r-1} \to X$ factors as
$\TCX_{r-1} \to X\times\BS_r \to X$, hence by lemma~\ref{zl1}~$(i)$ we have
$$
\pi_\#(F) = \pi_*(F\otimes\CO_X(r-1)\otimes\omega_{\BS_r}\otimes\det\TCL_{r-1}^*[\dim\BS_r-r+1]).
$$
Since line bundles $\omega_{\BS_r}$ and $\det\TCL_{r-1}^*$ are pullbacks from $\BS_r$,
it remains to check that
$$
\Phi^*_{i_*\CE}(\lan\CA_{r-1}(r),\dots,\CA_{\ix-1}(\ix)\ran^\perp) \subset \CB_{N-r}
\qquad\text{in $\D^b(Y)$}.
$$
But it is clear that
$\lan\CA_{r-1}(r),\dots,\CA_{\ix-1}(\ix)\ran^\perp = \lan\CA_0(1),\dots,\CA_{r-2}(r-1)\ran$
and in the notation of section~\ref{s5} we have $\Phi^*_{i_*\CE} = \gamma^*\pi^*$.
Now we see that the desired inclusion is proved in lemma~\ref{use1}.
\end{proof}


\begin{lemma}\label{phphph}
For all $r$ we have
$$
\arraycolsep=2pt
\begin{array}{lll}
\hfill      \Phi^*_{i_*\HE(0,1)\otimes(\TCL_r/\TCL_{r-1})}\circ\TPhi_{r-1} & = 0, &
\text{if $r\le N - \jx$}\\
\TPhi_r\circ\Phi^*_{i_*\HE(0,1)\otimes(\TCL_r/\TCL_{r-1})} & = 0, &
\text{if $r\ge \ix+1$}\\
\TPhi_r\circ\Phi^*_{i_*\HE(0,1)\otimes(\TCL_r/\TCL_{r-1})}\circ\TPhi_{r-1} &  = 0, \qquad\qquad &
\text{in other cases.}
\end{array}
$$
\end{lemma}
\begin{proof}
Note that by lemma~\ref{imkerphir}~$(i)$ we have
$\Im \TPhi_{r-1} \subset [\lan\CA_{r-1}(1),\dots,\CA_{\ix-1}(\ix-r+1)\ran\boxtimes\D^b(\BS_r)]^\perp$ and
if $r\ge \ix+1$ then the RHS equals $\D^b(\TCX_{r-1})$.
On the other hand, by lemma~\ref{imkerphir}~$(ii)$ we have
$\CB_{N-r}(-1)\otimes\D^b(\BS_r) \subset \Ker\TPhi_r$ and
if $r\le N - \jx$ then $N-r \ge \jx$ hence $\CB_{N-r} = 0$ by~(\ref{bk0}).
Thus it suffices to check that for all $r$ we have
$$
\Phi^*_{i_*\HE(0,1)\otimes(\TCL_r/\TCL_{r-1})}([\lan\CA_{r-1}(1),\dots,\CA_{\ix-1}(\ix-r+1)\ran\boxtimes\D^b(\BS_r)]^\perp)
\subset \CB_{N-r}(-1)\boxtimes\D^b(\BS_r).
$$
This easily follows from lemma~\ref{phiisemain}
(the functor $\Phi^*_{i_*\HE(0,1)\otimes(\TCL_r/\TCL_{r-1})}$
differs from $\Phi^*_{i_*\HE(0,1)}$ by the $(\TCL_r/\TCL_{r-1})$-twist,
and the line bundle $(\TCL_r/\TCL_{r-1})$ is a pullback from $\BS_r$).
\end{proof}

\subsection{The induction arguments}

\begin{proposition}\label{phir_spl}
The functors $\Phi_r$ are left splitting functors for all $r$.
\end{proposition}
\begin{proof}
We use induction in $r$. The functor $\Phi_1$ is fully faithful and
its image is left admissible by definition of Homological Projective Duality
(we have $\CY_1 = Y$ and $\Phi_1 = \Phi_\CE:\D^b(Y) \to \CC \subset \D^b(\CX_1)$),
hence $\Phi_1$ is left splitting.
Assume that $\Phi_{r-1}$ is left splitting. Then the functor $\TPhi_{r-1}$
is left splitting by proposition~\ref{fbc_spl}.
Now consider the functor $\TPhi_r\TPhi_r^*\TPhi_r\eta_*$.
Composing an isomorphism of functors~(\ref{fi1}) with $\TPhi_r\TPhi_r^*$
we obtain an isomorphism
$$
\TPhi_r\TPhi_r^*\TPhi_r\eta_* \cong \TPhi_r\TPhi_r^*\xi^*\TPhi_{r-1}.
$$
Composing exact triangle of functors~(\ref{ft3}) with $\TPhi_r$ (on the left)
and $\TPhi_{r-1}$ (on the right) we obtain an exact triangle of functors
$$
\TPhi_r\Phi^*_{i_*\HE(0,1)\otimes(\TCL_r/\TCL_{r-1})}\TPhi_{r-1} \to
\TPhi_r\TPhi_r^*\xi^*\TPhi_{r-1} \to
\TPhi_r\eta_*\TPhi_{r-1}^*\TPhi_{r-1}.
$$
Composing~(\ref{fi1}) with $\TPhi_{r-1}^*\TPhi_{r-1}$
we obtain an isomorphism
$$
\TPhi_r\eta_*\TPhi_{r-1}^*\TPhi_{r-1} \cong
\xi^*\TPhi_{r-1}\TPhi_{r-1}^*\TPhi_{r-1}.
$$
Using the induction assumption, criterion~\ref{sf_th}~$(2l)$ and
isomorphism (\ref{fi1}) we deduce
$$
\xi^*\TPhi_{r-1}\TPhi_{r-1}^*\TPhi_{r-1} \cong
\xi^*\TPhi_{r-1} \cong
\TPhi_r\eta_*.
$$
On the other hand,
$\TPhi_r\Phi^*_{i_*\HE(0,1)\otimes(\TCL_r/\TCL_{r-1})}\TPhi_{r-1} = 0$
by lemma~\ref{phphph}. Summarizing, we deduce that
$$
\TPhi_r\TPhi_r^*\TPhi_r\eta_* \cong \TPhi_r\eta_*.
$$
Finally since $\psi_*\eta_*\hat\psi^* = \hat\psi_*\hat\psi^* = \id$
by corollary~\ref{cim1}, we have
$$
\Phi_r\Phi_r^*\Phi_r \cong
\Phi_r\Phi_r^*\Phi_r\psi_*\eta_*\hat\psi^* \cong
\psi_*\TPhi_r\TPhi_r^*\TPhi_r\eta_*\hat\psi^* \cong
\psi_*\TPhi_r\eta_*\hat\psi^* \cong
\Phi_r\psi_*\eta_*\hat\psi^* \cong
\Phi_r
$$
by proposition~\ref{phiPhi}. Therefore $\Phi_r$ is left splitting
by theorem~\ref{sf_th}.
\end{proof}

\begin{corollary}
We have the following semiorthogonal collections
in $\D^b(\CX_r)$ and $\D^b(\CY_r)$:
$$
\begin{array}{lll}
\langle \Im\Phi_r,\CA_r(1)\boxtimes\D^b(\BP_r),\dots,\CA_{\ix-1}(\ix-r)\boxtimes\D^b(\BP_r)\rangle &\subset& \D^b(\CX_r)\\
\langle \CB_{\jx-1}(N-r-\jx)\boxtimes\D^b(\BP_r),\dots,\CB_{N-r}(-1)\boxtimes\D^b(\BP_r),\Im\Phi_r^*\rangle &\subset& \D^b(\CY_r).
\end{array}
$$
Moreover, the functors $\Phi_r$ and $\Phi_r^*$ induce an equivalence
$\Im\Phi_r \cong \Im\Phi_r^*$.
\end{corollary}
\begin{proof}
Combine theorem~\ref{sf_th} with proposition~\ref{imkerphir} and lemma~\ref{soc_xyr}.
\end{proof}

It remains to check that these semiorthogonal collections are full,
i.e.\ that they generate the derived categories of $\CX_r$ and $\CY_r$.
We begin with the case of $\CY_r$.

\begin{proposition}\label{comply}
For all $r$ we have
$$
\D^b(\CY_r) = \langle \CB_{\jx-1}(N-r-\jx)\boxtimes\D^b(\BP_r),\dots,\CB_{N-r}(-1)\boxtimes\D^b(\BP_r),\Im\Phi_r^*\rangle.
$$
\end{proposition}
\begin{proof}
Note that for $s\ge \jx$ we have $\CB_s = 0$ by (\ref{bk0}).
Hence for $r\le N - \jx$ the RHS coincides with $\Im\Phi_r^*$, and since
we already know that $\Phi_r$ is splitting, for $r\le N - \jx$ it suffices
to check that $\Ker\Phi_r = 0$, that is that $\Phi_r$ is fully faithful.
The arguments are the same as in the proof of proposition~\ref{phir_spl},
the only difference is that we use the first equality of lemma~\ref{phphph}
instead of the third.

We use induction in $r$. In the case $r=1$ we have $\CY_1 = Y$ and $\Phi_1$
is fully faithful by definition of Homological Projective Duality.
Now let $1 < r \le N - \jx$ and assume that $\Phi_{r-1}$ is fully faithful.
Then the functor $\TPhi_{r-1}$ is fully faithful by proposition~\ref{fbc_sod}.
Now consider the functor $\TPhi^*_r\TPhi_r\eta_*$.
Composing an isomorphism of functors~(\ref{fi1}) with $\TPhi_r^*$
we obtain an isomorphism
$$
\TPhi^*_r\TPhi_r\eta_* \cong
\TPhi^*_r\xi^*\TPhi_{r-1}.
$$
Composing exact triangle of functors~(\ref{ft3}) with $\TPhi_{r-1}$
we obtain an exact triangle of functors
$$
\Phi^*_{i_*\HE(0,1)\otimes(\TCL_r/\TCL_{r-1})}\TPhi_{r-1} \to
\TPhi_r^*\xi^*\TPhi_{r-1} \to
\eta_*\TPhi_{r-1}^*\TPhi_{r-1}.
$$
Since the functor $\TPhi_{r-1}$ is fully faithful we have
$$
\eta_*\TPhi_{r-1}^*\TPhi_{r-1} \cong \eta_*.
$$
On the other hand,
$\Phi^*_{i_*\HE(0,1)\otimes(\TCL_r/\TCL_{r-1})}\TPhi_{r-1} = 0$
by lemma~\ref{phphph} since $r\le N - \jx$. Summarizing, we deduce that
$$
\TPhi^*_r\TPhi_r\eta_* \cong \eta_*.
$$
Finally since $\psi_*\eta_*\hat\psi^* = \hat\psi_*\hat\psi^* = \id$
by corollary~\ref{cim1}, we have
$$
\Phi^*_r\Phi_r \cong
\Phi^*_r\Phi_r\psi_*\eta_*\hat{\psi}^* \cong
\psi_*\TPhi^*_r\TPhi_r\eta_*\hat{\psi}^* \cong
\psi_*\eta_*\hat{\psi}^* \cong
\id_{\CY_r}
$$
by proposition~\ref{phiPhi}. Therefore $\Phi_r$ is fully faithful.

For $r\ge N-\jx$ we also use induction in $r$.
However the arguments are slightly different in this case.
Assume that $r\ge N-\jx+1$ and the claim for $r-1$ is true.
Assume that $G$ is in the right orthogonal to the category
$\langle \CB_{\jx-1}(N-r-\jx)\boxtimes\D^b(\BP_r),\dots,\CB_{N-r}(-1)\boxtimes\D^b(\BP_r),\Im\Phi_r^*\rangle$.
By proposition~\ref{phiPhi} we have $\TPhi_r(\psi^*G) = \psi^*\Phi_r(G)$.
But $\Phi_r(G) = 0$ since $G \in (\Im\Phi^*_r)^\perp$,
hence $\TPhi_r(\psi^*G) = 0$, thus $\psi^*G \in (\Im\TPhi^*_r)^\perp$.
On the other hand, for any $F \in \CB_{N-r+t}(-1-t)$, $H \in \D^b(\BS_r)$ we have
\begin{multline*}
\psi_*\RCHom(\psi^*F\otimes H,\psi^*G) \cong
\psi_*(\RCHom(\psi^*F,\psi^*G)\otimes H^*) \cong
\psi_*(\psi^*\RCHom(F,G)\otimes H^*) \cong
\\ \cong
\RCHom(F,G)\otimes \psi_*(H^*) \cong
\RCHom(F\otimes (\psi_*(H^*))^*,G)
\end{multline*}
whereof we deduce that $\psi^*G \in (\CB_{N-r+t}(-1-t) \otimes \D^b(\BS_r))^\perp$.
Combining these two inclusions we see that
$\psi^* G \in
\langle \CB_{\jx-1}(N-r-\jx)\boxtimes\D^b(\BS_r),\dots,\CB_{N-r}(-1)\boxtimes\D^b(\BS_r),
\Im\TPhi_r^*\rangle^\perp$.
Recalling lemma~\ref{phiisemain} we note that
$\psi^*G$ is in the right orthogonal to the image of the first two terms
of the triangle~(\ref{ft3}) applied to the subcategory
$[\lan\CA_{r-1}(1),\dots,\CA_{\ix-1}(\ix-r+1)\ran\boxtimes\D^b(\BS_r)]^\perp$.
Therefore $\psi^*G$ is in the right orthogonal to the image
of $\eta_*\TPhi_{r-1}^*$ applied to the subcategory
$[\lan\CA_{r-1}(1),\dots,\CA_{\ix-1}(\ix-r+1)\ran\boxtimes\D^b(\BS_r)]^\perp$.
But from lemma~\ref{imkerphir} it follows by adjunction that
$[\lan\CA_{r-1}(1),\dots,\CA_{\ix-1}(\ix-r+1)\ran\boxtimes\D^b(\BS_r)] \subset \Ker\TPhi_{r-1}^*$.
Therefore $\psi^*G$ is in the right orthogonal to the image of the functor
$\eta_*\TPhi_{r-1}^*$. By adjunction we deduce
$\eta^!\psi^*G \in (\Im\TPhi^*_{r-1})^\perp$.


On the other hand, by lemma~\ref{zl1}~$(iv)$ we have the following resolution
$$
0 \to \CO_{\TCY_r}(-1)\otimes(\TCL_r/\TCL_{r-1})^* \to \CO_{\TCY_r} \to \eta_*\CO_{\TCY_{r-1}} \to 0
$$
which implies that
$\eta_*(\CB_{N-r+t}(-1-t)\boxtimes\D^b(\BS_r)) \subset
\langle\CB_{N-r+t}(-2-t),\CB_{N-r+t}(-1-t)\rangle\boxtimes\D^b(\BS_r)$ for any $t$,
and since we have $\CB_{N-r+t} \subset \CB_{N-r+t-1}$ we conclude that
$$
\eta_*\langle\CB_{\jx-1}(N-r-\jx+1),\dots,\CB_{N-r+1}(-1)\rangle\boxtimes\D^b(\BS_r) \subset
\langle\CB_{\jx-1}(N-r-\jx),\dots,\CB_{N-r}(-1)\rangle\boxtimes\D^b(\BS_r).
$$
Therefore
$\eta^!\psi^*G \in [\langle\CB_{\jx-1}(N-r-\jx+1),\dots,\CB_{N-r+1}(-1)\rangle\boxtimes\D^b(\BS_r)]^\perp$.
Summarizing we deduce that $\eta^!\psi^*G = 0$ by the induction hypothesis.
But $\eta^!\psi^*G = \eta^*\psi^*G \otimes \CO_Y(1)\otimes(\TCL_r/\TCL_{r-1})[-1]$
by lemma~\ref{zl1}~$(iv)$,
and $\eta^*\psi^* = {\hat\psi}^*$ is a fully faithful functor by corollary~\ref{cim1},
hence $G = 0$.
\end{proof}

\begin{corollary}\label{ecy}
We have $\D^b(Y) = \lan\CB_{\jx-1}(1-\jx),\dots,\CB_1(-1),\CB_0\ran$.
\end{corollary}
\begin{proof}
Note that $\BP_N$ is a point, $\CY_N = Y$, $\CX_N = \emptyset$,
and apply proposition~\ref{comply}
for $r=N$.
\end{proof}

Fullness of the semiorthogonal decomposition of $\D^b(\CX_r)$ will be proved
by a decreasing induction. The base of induction is given by the following

\begin{lemma}\label{imphin1}
We have $\Im\Phi_{N-1} = \D^b(\CX_{N-1})$.
\end{lemma}
\begin{proof}
Note that
$\CX_{N-1} = X$ and the projection $\pi:\CX_{N-1} \to X$ is the identity.
Therefore, the functor $\pi_*\circ\Phi_{N-1}$ considered in the proof
of lemma~\ref{imkerphir}~$(i)$ coincides with $\Phi_{N-1}$.
Further, note that $\CY_{N-1} \cong \PP_Y(\Ker(V\otimes\CO_Y \to \CO_Y(1)))$,
and, moreover, the projection to $Y$ coincides with $\pi$ and
the Grothendieck ample line bundle coincides with $\CL_{N-1}^{\perp*}$.
It easily follows that for $0\le t,s\le N-2$ we have
$$
\pi_*(\Lambda^t(\CL_{N-1}/\CO_Y(-1))\otimes S^s(\CL_{N-1}^\perp)) \cong
\begin{cases}
\CO_Y[-t], & \text{if $t=s$}\\
0, & \text{otherwise}
\end{cases}
$$
Hence in the notation of the proof of lemma~\ref{imkerphir}~$(i)$ we have
$$
\Psi_t(\pi^*G\otimes S^s(\CL_{N-1}^\perp)) \cong
\begin{cases}
\Phi_\CE(G)\otimes\CO_X(-t)[-t], & \text{if $t=s$}\\
0, & \text{otherwise}
\end{cases}
$$
It follows that
$$
\Phi_{N-1}(\pi^*G\otimes S^s(\CL_{N-1}^\perp)) \cong
\Phi_\CE(G)\otimes\CO_X(-s).
$$
But $\Im\Phi_\CE = \CA_0$ by lemma~\ref{imkerpig},
hence $\Im\Phi_{N-1}$ contains $\CA_0(-s)$ for $0\le s \le N-2$.
It remains to note that $\Im\Phi_{N-1}$ is triangulated subcategory of $\D^b(X)$
since $\Phi_{N-1}$ is a splitting functor by proposition~\ref{phir_spl},
and on the other hand by lemma~\ref{a0rgen} we have
$\lan\CA_0(2-N),\dots,\CA_0\ran = \lan\CA_0(2-N),\dots,\CA_{\ix-1}(\ix+1-N)\ran$
(note that $N > \ix$ by assumption~(\ref{ngi})), and the latter category
evidently coincides with $\D^b(X) = \D^b(\CX_{N-1})$.
\end{proof}

\begin{proposition}\label{complx}
For all $r$ we have
$$
\D^b(\CX_r) = \lan \Im\Phi_r,\CA_r(1)\boxtimes\D^b(\BP_r),\dots,\CA_{\ix-1}(\ix-r)\boxtimes\D^b(\BP_r)\ran.
$$
\end{proposition}
\begin{proof}
The arguments are analogous to those used in the proof of proposition~\ref{comply}.

Note that for $s\ge \ix$ we have $\CA_s = 0$.
Hence for $r\ge \ix$ the RHS coincides with $\Im\Phi_r$, and since
we already know that $\Phi_r$ is splitting, for $r\ge \ix$ it suffices
to check that $\Ker\Phi^*_r = 0$, that is that $\Phi^*_r$ is fully faithful.
For this we use descending induction in $r$.

In the case $r=N-1$ we know that $\Im\Phi_{N-1} = \D^b(\CX_{N-1})$
by lemma~\ref{imphin1}. Now let $\ix < r \le N - 1$ and assume that
$\Phi_{r}$ is fully faithful. Then the functor $\TPhi_{r}$ is fully
faithful by proposition~\ref{fbc_sod}.
Now consider the functor $\xi^*\TPhi_{r-1}\TPhi^*_{r-1}$.
Composing exact triangle of functors~(\ref{ft3}) with $\TPhi_{r}$
we obtain an exact triangle of functors
$$
\TPhi_{r}\Phi^*_{i_*\HE(0,1)\otimes(\TCL_r/\TCL_{r-1})} \to
\TPhi_{r}\TPhi_r^*\xi^* \to
\TPhi_{r}\eta_*\TPhi_{r-1}^*.
$$
Composing isomorphism of functors~(\ref{fi1}) with $\TPhi_{r-1}^*$
we obtain an isomorphism
$$
\TPhi_{r}\eta_*\TPhi_{r-1}^* \cong
\xi^*\TPhi_{r-1}\TPhi_{r-1}^*.
$$
Since the functor $\TPhi_{r}$ is fully faithful we have
$$
\TPhi_{r}\TPhi^*_{r}\xi^* \cong \xi^*.
$$
On the other hand,
$\TPhi_{r}\Phi^*_{i_*\HE(0,1)\otimes(\TCL_r/\TCL_{r-1})} = 0$
by lemma~\ref{phphph} since $r\ge \ix+1$. Summarizing, we deduce that
$$
\xi^*\TPhi_{r-1}\TPhi_{r-1}^* \cong \xi^*.
$$
Finally, since $\hat{\phi}_*\xi^*\phi^* = \hat{\phi}_*\hat{\phi}^* = \id$
by corollary~\ref{cim1}, we have
$$
\Phi_{r-1}\Phi_{r-1}^* \cong
\hat{\phi}_*\xi^*\phi^*\Phi_{r-1}\Phi_{r-1}^* \cong
\hat{\phi}_*\xi^*\TPhi_{r-1}\TPhi_{r-1}^*\phi^* \cong
\hat{\phi}_*\xi^*\phi^* \cong
\id_{\CX_{r-1}}
$$
by proposition~\ref{phiPhi}. Therefore $\Phi^*_{r-1}$ is fully faithful.

For $r < \ix$ we also use induction in $r$.
However the arguments are slightly different in this case.
Assume that $r\le \ix$ and the claim for $r$ is true.
Assume that $F$ is in the left orthogonal to $\Im\Phi_{r-1}$
and in the right orthogonal to the category
$\lan \CA_{r-1}(1)\boxtimes\D^b(\BP_{r-1}),\dots,\CA_{\ix-1}(\ix-r+1)\boxtimes\D^b(\BP_{r-1})\ran$.
Then by the same arguments as in the proof of proposition~\ref{comply} we check that
$\phi^* F$ is in the left orthogonal to $\Im\TPhi_{r-1}$ and
in the right orthogonal to the category
$\lan \CA_{r-1}(1)\boxtimes\D^b(\BS_r),\dots,\CA_{\ix-1}(\ix-r+1)\boxtimes\D^b(\BS_r)\ran$.
Note that by adjunction it follows from lemma~\ref{phiisemain} that
$$
\Phi_{i_*\HE(0,1)\otimes(\TCL_r/\TCL_{r-1})}([\CB_{N-r}(-1)\boxtimes\D^b(\BS_r)]^\perp) \subset
\lan\CA_{r-1}(1)\boxtimes\D^b(\BS_{r-1}),\dots,\CA_{\ix-1}(\ix-r+1)\boxtimes\D^b(\BS_{r-1})\ran^{\perp\perp}.
$$
We deduce that $\phi^*F$ is in the left orthogonal to the image of the first
and the third terms of the triangle~(\ref{ft1}) applied to the subcategory
$[\CB_{N-r}(-1)\boxtimes\D^b(\BS_r)]^\perp$.
Therefore $\phi^*F$ is in the left orthogonal to the image
of $\xi_*\TPhi_r$ applied to the subcategory
$[\CB_{N-r}(-1)\boxtimes\D^b(\BS_r)]^\perp$.
But $[\CB_{N-r}(-1)\boxtimes\D^b(\BS_r)] \subset \Ker\TPhi_r$
by lemma~\ref{imkerphir}, hence $\phi^*F$ is in the left orthogonal
to the image of the functor $\xi_*\TPhi_r$.
and by adjunction we deduce that
$\xi^*\phi^*F \in {}^\perp(\Im\TPhi_{r})$.

On the other hand, by lemma~\ref{zl1}~$(iii)$ we have the following resolution
$$
0 \to \CO_{\TCX_{r-1}}(-1)\otimes(\TCL_r/\TCL_{r-1}) \to \CO_{\TCX_{r-1}} \to \xi_*\CO_{\TCX_r} \to 0
$$
which implies that
$\xi_*(\CA_{r+t}(1+t)\boxtimes\D^b(\BS_r)) \subset
\langle\CA_{r+t}(t),\CA_{r+t}(1+t)\rangle\boxtimes\D^b(\BS_r)$ for any $t$,
and since we have $\CA_{r+t} \subset \CB_{r+t-1}$ we conclude that
$$
\xi_*\lan\CA_r(1)\boxtimes\D^b(\BS_r),\dots,\CA_{\ix-1}(\ix-r)\boxtimes\D^b(\BS_r)\ran \subset
\lan\CA_{r-1}\boxtimes\D^b(\BS_{r-1}),\dots,\CA_{\ix-1}(\ix-r)\boxtimes\D^b(\BS_{r-1})\ran.
$$
Therefore
$\xi^!\phi^*F(-1) \in \lan\CA_r(1)\boxtimes\D^b(\BS_r),\dots,\CA_{\ix-1}(\ix-r)\boxtimes\D^b(\BS_r)\ran^\perp$.
Finally, we have $\xi^!\phi^*F(-1) \cong \xi^*\phi^*G\otimes(\TCL_r/\TCL_{r-1})^*[-1]$
by lemma~\ref{zl1}~$(iii)$.
Summarizing we deduce that $\xi^*\phi^*F = 0$ by the induction hypothesis.
But $\xi^*\phi^* = \hat{\phi}^*$ is fully faithful by corollary~\ref{cim1},
hence $F = 0$.
\end{proof}

\subsection{Proof of the main theorem}

In this subsection we prove theorem~\ref{themain} and describe
some of its generalizations.

First of all, the first claim of the theorem is proved in lemma~\ref{ysm}
and corollary~\ref{ecy}. For the second claim, let $L\subset V^*$ be
an admissible subspace, $\dim L = r$.
Then the map $\lambda:\Spec\kk \to \BP_r$ induced by $L$ is
a faithful base change for the pair $(\CX_r,\CY_r)$.
by lemmas~\ref{ec}~$(iii)$, \ref{zl1}, \ref{xyksm} and
the definition of admissible subspace.
Therefore we can apply the faithful base change theorem~\ref{phitsod}.
Then theorem~\ref{themain} follows from proposition~\ref{complx}
and proposition~\ref{comply}.
\qed

\begin{remark}\label{nonadm}
An interesting question is what can we say about derived categories
of $X_L$ and $Y_L$ when $L$ is not admissible. In fact one can show
that replacing the {\em naive}\/ linear sections $X_L$ and $Y_L$
by the {\em derived}\/ linear sections (i.e.\ considering
the {\em derived}\/ fiber products
$X{\mathop{\times}\limits^{\mathbb{L}}}_{\PP(V)}\PP(L^\perp)$ and
$Y{\mathop{\times}\limits^{\mathbb{L}}}_{\PP(V^*)}\PP(L)$)
one can manage to get the same semiorthogonal decompositions
for their derived categories.
\end{remark}

The same argument with $\Spec\kk$ replaced by any base scheme $T$
proves a relative version of theorem~\ref{themain}:

\begin{theorem}\label{refthemain}
Assume that $T$ is a smooth algebraic variety and
$\CL \subset V^*\otimes\CO_T$ is a vector subbundle,
$\rank \CL = r$, such that the corresponding families
of linear sections of $X$ and $Y$
$$
\begin{array}{lll}
\CX_\CL & =
(X\times T)\times_{\PP(V)\times T}\PP_{T}(\CL^\perp) &
\subset X\times T, \\
\CY_\CL & =
(Y\times T)\times_{\PP(V^*)\times T}\PP_{T}(\CL) &
\subset Y\times T,
\end{array}
$$
and their fiber product $\CX_\CL\times_T \CY_\CL$ have expected dimension
$$
\begin{array}{l}
\dim\CX_\CL = \dim X + \dim T - r,\qquad
\dim\CY_\CL = \dim Y + \dim T + r - N,\\
\dim\CX_\CL\times_{T}\CY_\CL = \dim X + \dim Y + \dim T - N.
\end{array}
$$
Then there exist a triangulated category $\CC_\CL$ and semiorthogonal decompositions
$$
\arraycolsep = 2pt
\begin{array}{lll}
\D^b(\CX_\CL) &=& \langle \CC_\CL,
\CA_r(1)\boxtimes\D^b(T),\dots,\CA_{\ix-1}(\ix-r)\boxtimes\D^b(T)\rangle\\
\D^b(\CY_\CL) &=& \langle \CB_{\jx-1}(N-r-\jx)\boxtimes\D^b(T),\dots,
\CB_{N-r}(-1)\boxtimes\D^b(T),\CC_\CL\rangle.
\end{array}
$$
\end{theorem}

\begin{remark}\label{relhpd}
Another relative version of theorem~\ref{themain} can be obtained
as follows. Consider a base scheme~$S$ (not necessarily compact),
assume that $X$ and $Y$ are algebraic varieties over $S$,
replace the assumptions of projectivity of the maps
$f:X \to \PP(V)$ and $g:Y\to\PP(V^*)$ in the definition
of Homological Projective Duality by projectivity
of the maps $X \to S\times\PP(V)$ and $Y \to S\times\PP(V^*)$,
and assume that we are given a Lefschetz decomposition of $\D^b(X)$
which is $S$-linear.
We will say that $Y$ is
{\sf Homologically Projectively Dual to $X$ relatively over $S$}\/
if there exists an object $\CE \in \D^b(Q(X,Y)\times_{S\times S}S)$
(the fiber product is taken with respect to the canonical map
$Q(X,Y) \subset X\times Y \to S\times S$ and
the diagonal embedding $S \stackrel\Delta\to S\times S$)
such that the functor $\Phi_\CE$ is fully faithful and gives
semiorthogonal decomposition~(\ref{dbx1}).
One can prove by the same arguments that theorem~\ref{themain}
and theorem~\ref{refthemain} are true in this case as well.
\end{remark}

\section{Properties of homological projective duality}

We believe that phenomenon of Homological Projective Duality
deserves to be thoroughly investigated. In this section we will
discuss some basic properties of Homological Projective Duality.

The first natural question is when a Homologically Projectively Dual
variety for a given algebraic variety $X$ exists. From the definition
of the Homological Projective Duality it follows that it always exists
{\em on a categorical level}, we always know the derived category
of the Homologically Projectively Dual variety. On the other hand,
the question of existence of a Homologically Projectively Dual variety
{\em on a geometrical level}\/ seems to be of a philosophical nature.
Indeed, in some sense every sufficiently good triangulated category
can be considered as the derived category of coherent sheaves
on a {\em noncommutative algebraic variety}. In fact, this is one
of the ways to understand what a noncommutative algebraic variety is.
From this point of view the question of existence of a Homologically
Projectively Dual variety as a usual commutative variety seems to be
not very natural and it is difficult to expect a nice answer
(especially if we remember that the notion of Homological Projective
Duality depends on a choice of a line bundle and a Lefschetz decomposition).

The next question is whether a Homologically Projectively Dual variety is unique.
Certainly this is true if it is a Fano variety by the Reconstruction Theorem
of Bondal and Orlov~\cite{BO4}. However, in general it doesn't need to be Fano,
so there are examples of several different Homologically Projectively Dual varieties.

Another important question is how one should {\em construct} a Homological Projectively
Dual variety for a given variety $X$. A natural approach is to consider a moduli space
of objects in $\D^b(X)$ with a given class in $K_0(X)$ and {\em supported on hyperplane
sections}\/ of $X$. However, there are two problems on this way. The first one is
of technical nature --- we don't have a good theory of moduli spaces of objects
in triangulated categories yet (moduli spaces in a ``good theory'' should depend
on a choice of stability conditions and should be ``noncommutative'' in general).
The second problem is more complicated --- how to choose a correct class in $K_0(X)$.
There is a trivial restriction on this class --- it should be orthogonal to
subcategories $\CA_1(1)$, \dots $\CA_{\ix-1}(\ix-1)$ of the Lefschetz decomposition of $X$.
Sometimes, these restrictions determine unique class in $K_0(X)$ up to a multiplicity.
However, examples considered in~\cite{K} show that the choice of correct multiplicity
turns out to be quite mysterious.

Now we turn to more specific questions.

\subsection{Disjoint unions and products}

If algebraic variety $X$ is a disjoint union, $X = X' \sqcup X''$
then its derived category is a completely orthogonal direct sum,
$\D^b(X) = \D^b(X')\oplus\D^b(X'')$. If we are given
Lef\-schetz decompositions
$\D^b(X') = \lan \CA'_0,\CA'_1(1),\dots,\CA'_{\ix'-1}(\ix'-1)\ran$,
$\D^b(X'') = \lan \CA''_0,\CA''_1(1),\dots,\CA''_{\ix''-1}(\ix''-1)\ran$,
then we have a Lefschetz decomposition
$$
\D^b(X) = \lan \CA'_0\oplus\CA''_0,(\CA'_1\oplus\CA''_1)(1),\dots,(\CA'_{\ix-1}\oplus\CA''_{\ix-1})(\ix-1)\ran,
\qquad \ix = \max\{\ix',\ix''\}.
$$
It is natural question, if we know a Homologically Projectively Dual varieties
to $X'$ and $X''$, what will be Homologically Projectively Dual to $X'\sqcup X''$?
The answer is quite simple

\begin{proposition}
If $Y'$ and $Y''$ are Homologically Projectively Dual to $X'$ and $X''$ respectively
then $Y = Y'\sqcup Y''$ is Homologically Projectively Dual to $X = X' \sqcup X''$.
\end{proposition}
\begin{proof}
Note that the universal hyperplane section of $X$ can be represented as
$\CX_1 = \CX'_1 \sqcup \CX''_1$. Combining semiorthogonal decompositions
$\D^b(\CX'_1) = \lan\D^b(Y'),\CA'_1(1)\boxtimes\D^b(\PP(V^*)),\dots,\CA'_{\ix'-1}(\ix'-1)\boxtimes\D^b(\PP(V^*))\ran$
and
$\D^b(\CX''_1) = \lan\D^b(Y''),\CA''_1(1)\boxtimes\D^b(\PP(V^*)),\dots,\CA''_{\ix''-1}(\ix''-1)\boxtimes\D^b(\PP(V^*))\ran$
we get
$$
\D^b(\CX_1) = \lan\D^b(Y')\oplus\D^b(Y''),
(\CA'_1\oplus\CA''_1)(1)\boxtimes\D^b(\PP(V^*)),\dots,
(\CA'_{\ix-1}\oplus\CA''_{\ix-1})(\ix-1)\boxtimes\D^b(\PP(V^*))\ran
$$
which shows that $Y = Y'\sqcup Y''$ is Homologically Projectively Dual to $X = X' \sqcup X''$.
\end{proof}

In other word, Homological Projective duality {\em commutes with disjoint unions}.

Now assume that $X = X'\times F$ and take $\CO_X(1) = p^*\CO_{X'}(1)$, where
$p:X \to X'$ is the projection along~$F$. Then we have a Lefschetz decomposition
$$
\D^b(X) = \lan \CA'_0\boxtimes\D^b(F),(\CA'_1\boxtimes\D^b(F))(1),\dots,
(\CA'_{\ix'-1}\boxtimes\D^b(F))(\ix'-1)\ran.
$$

\begin{proposition}
If $Y'$ is Homologically Projectively Dual to $X'$ then
$Y = Y'\times F$ is Homologically Projectively Dual to $X = X' \times F$.
\end{proposition}
\begin{proof}
Note that the universal hyperplane section of $X$ can be represented as
$\CX_1 = \CX'_1 \times F$. Tensoring semiorthogonal decomposition
$\D^b(\CX'_1) = \lan\D^b(Y'),\CA'_1(1)\boxtimes\D^b(\PP(V^*)),\dots,\CA'_{\ix'-1}(\ix'-1)\boxtimes\D^b(\PP(V^*))\ran$
with $\D^b(F)$
we get
$$
\D^b(\CX_1) = \lan\D^b(Y')\boxtimes\D^b(F),
(\CA'_1\boxtimes\D^b(F))(1)\boxtimes\D^b(\PP(V^*)),\dots,
(\CA'_{\ix'-1}\boxtimes\D^b(F))(\ix'-1)\boxtimes\D^b(\PP(V^*))\ran
$$
which shows that $Y = Y'\times F$ is Homologically Projectively Dual to $X = X' \times F$.
\end{proof}

\subsection{Duality}

In this subsection we are going to check that relation of Homological Projective Duality
is a duality indeed.

\begin{theorem}
If $g:Y \to \PP(V^*)$ is Homologically Projectively Dual to $f:X \to \PP(V)$ then
$f:X \to \PP(V)$ is Homologically Projectively Dual to $g:Y \to \PP(V^*)$.
\end{theorem}
\begin{proof}
Indeed, let $\D^b(X) = \lan\CA_0,\CA_1(1),\dots,\CA_{\ix-1}(\ix-1)\ran$
be the Lefschetz decomposition of $\D^b(X)$ and
let $\D^b(Y) = \lan\CB_{\jx-1}(1-\jx),\dots,\CB_1(-1),\CB_0\ran$
be the dual Lefschetz decomposition of $\D^b(Y)$
given by theorem~\ref{themain}. Dualizing, we obtain a Lefschetz decomposition
$\D^b(Y) = \lan \CB_0^*,\CB_1^*(1),\dots,\CB_{\jx-1}^*(\jx-1)\ran$.
Let us show that $X$ is Homologically Projectively Dual to $Y$
with respect to this Lefschetz decomposition.

Indeed, consider $\CX_{N-1}$ and $\CY_{N-1}$.
Note that $\BP_{N-1} = \Gr(N-1,V^*) \cong \PP(V)$,
$\CX_{N-1} \cong X$ (its embedding into $X\times\Gr(N-1,V^*) = X\times\PP(V)$
is given by the graph of $f$) and $\CY_{N-1} \subset Y\times\PP(V)$ is
the universal hyperplane section of $Y$.
Dualizing the decomposition of proposition~\ref{comply} with $r = N - 1$
we obtain a semiorthogonal decomposition
$$
\D^b(\CY_{N-1}) =
\lan \D^b(\CX_{N-1})^*,
\CB_1^*(1)\boxtimes\D^b(\PP(V)),\dots,
\CB_{\jx-1}^*(\jx-1)\boxtimes\D^b(\PP(V))\ran.
$$
Moreover, the embedding functor $\D^b(\CX_{N-1}) \to \D^b(\CY_{N-1})$
is obtained by conjugation with the duality functor
of the functor $\Phi_{N-1}^*$.
Note that $\CY_{N-1}\times_{\PP(V)}\CX_{N-1} = Q(X,Y)$ and it is easy
to check that the embedding functor $\D^b(\CX_{N-1}) \to \D^b(\CY_{N-1})$
is a kernel functor with kernel scheme-theoretically supported on $Q(X,Y)$.
Thus $X$ is Homologically Projectively Dual to $Y$.
\end{proof}

\subsection{Dimension of the dual variety}

It is natural to ask, what can we say about the dimension of the Homologically
Projectively Dual variety to a given variety $X$. This question can be answered
precisely in the special case of {\em rectangular}\/ Lefschetz decomposition
of $X$ (i.e.\ when $\CA_0 = \CA_1 = \dots = \CA_{\ix-1}$).

\begin{proposition}
If $g:Y \to \PP(V^*)$ is Homologically Projectively Dual to $f:X \to \PP(V)$
with respect to a rectangular Lefschetz decomposition with $\ix$ terms,
then the number of terms $\jx$ in the dual Lefschetz decomposition of $\D^b(Y)$
and the dimension of $Y$ equal
$$
\jx = N - \ix,
\qquad
\dim Y = \dim X + N - 2\ix,
$$
where $N = \dim V$.
\end{proposition}
\begin{proof}
The formula for $\jx$ follows immediately from~\ref{defj}. To get the formula
for the dimension we note that in the case of rectangular Lefschetz decompositions
for any $\ix$-dimensional admissible subspace $L \subset V^*$ we have
by theorem~\ref{themain} an equivalence of categories $\D^b(X_L) = \D^b(Y_L)$.
It follows that $\dim X_L = \dim Y_L$. But $\dim X_L = \dim X - \ix$
and $\dim Y_L = \dim Y - (N - \ix)$.
\end{proof}

In general however it seems that it is impossible to give an explicit
formula for $\jx$ and $\dim Y$. However, we can prove an inequality

\begin{lemma}
If $g:Y \to \PP(V^*)$ is Homologically Projectively Dual to $f:X \to \PP(V)$
with respect to a Lefschetz decomposition with $\ix$ terms,
then the number of terms $\jx$ in the dual Lefschetz decomposition of $\D^b(Y)$
and the dimension of $Y$ satisfy
$$
\jx \ge  N - \ix,
\qquad
\dim Y \ge \dim X + N - 2\ix,
$$
where $N = \dim V$.
\end{lemma}
\begin{proof}
The formula for $\jx$ also follows from~\ref{defj}. To get the formula
for the dimension we note that for any $\ix$-dimensional admissible subspace $L \subset V^*$
we have by theorem~\ref{themain} a fully faithful functor $\D^b(X_L) \to \D^b(Y_L)$.
It follows that $\dim X_L \le \dim Y_L$. But $\dim X_L = \dim X - \ix$
and $\dim Y_L = \dim Y - (N - \ix)$.
\end{proof}



Another question is the relation of the number of terms of
a Lefschetz decomposition to the dimension of $X$.

\begin{proposition}\label{inddim}
Assume that $X$ is connected, $\D^b(X) = \lan\CA_0,\CA_1(1),\dots,\CA_{\ix-1}(\ix-1)\ran$
is a Lefschetz decomposition and $\CA_{\ix-1} \ne 0$.
Then $\dim X \ge \ix-1$ and equality is possible only
if $f:X \to \PP(V)$ is birational onto
a linear subspace $\PP^{\ix-1} \subset \PP(V)$ and $\CA_{\ix-2} = \CA_{\ix-1}$.
If moreover $\CA_0 = \CA_{\ix-1}$ then $X = \PP^{\ix-1}$.
\end{proposition}
\begin{proof}
Consider generic subspace $L \subset V^*$ of dimension $\ix-1$.
Then $X_L$ is a complete intersection of $\ix-1$ hyperplanes in $X$.
By theorem~\ref{themain} the restriction functor $\CA_{\ix-1} \subset \D^b(X) \to \D^b(X_L)$
is fully faithful. Therefore $\D^b(X_L) \ne 0$, so $X_L$ is not empty.
Since this is true for generic $L$ we conclude that $\dim X \ge \ix-1$.

Assume that $\dim X = \ix-1$ and consider generic subspace $L \subset V^*$
of dimension $\ix-2$. Then $X_L$ is a complete intersection of $\ix-2$
hyperplanes in $X$, so $X_L$ is a smooth connected curve for generic~$L$.
By theorem~\ref{themain}
the restriction functors $\CA_{\ix-2} \subset \D^b(X) \to \D^b(X_L)$
and $\CA_{\ix-1} \subset \D^b(X) \to \D^b(X_L)$ are fully faithful
and $\lan\CA_{\ix-2},\CA_{\ix-1}(1)\ran$ is a semiorthogonal collection
in $\D^b(X_L)$. But by lemma~\ref{ldcurve} the only smooth connected curve
admitting a nontrivial Lefschetz decomposition is $\PP^1$ and the decomposition
necessarily takes form $\D^b(\PP^1) = \lan \CO_{\PP^1}(k),\CO_{\PP^1}(k+1)\ran$.
Therefore $f$ restricted to $X_L$ is an isomorphism onto
a line $\PP^1 \subset \PP(V)$. When $L$ varies these lines
span an open subset of a linear subspace $\PP^{\ix-1}\subset\PP(V)$,
and $f$ is an isomorphism on the preimage of this open subset.
Hence $f$ is birational onto $\PP^{\ix-1}$. Moreover,
taking any indecomposable object $E \in \CA_{\ix-1}$ we deduce
that its restriction to $D^b(X_L)$ is isomorphic (up to a shift)
to $\CO_{\PP^1}(k)$ (the only indecomposable object
in $\lan\CO_{\PP^1}(k)\ran$). Hence $E$ is exceptional
and $(E,E(1),\dots,E(\ix-1))$ is an exceptional collection on $X$.


Finally, if $\CA_0 = \CA_{\ix-1}$ then
$\D^b(X) = \lan E,E(1),\dots,E(\ix-1)\ran$,
hence the Grothendieck group $K_0(X)$ is a free abelian group of rank $\ix$.
On the other hand, it is easy to see that if the birational map $X \to \PP^{\ix-1}$
is not trivial then the rank of $K_0(X)$ is strictly greater then $\ix$.
Hence $X = \PP^{\ix-1}$.
%
\end{proof}

\begin{remark}
If $X = \PP(W) \subset \PP(V)$, a linear subspace, $\dim W = \ix$,
considered with the Lefschetz decomposition
$\D^b(X) = \lan f^*\CO_{\PP(V)}(k),\dots,f^*\CO_{\PP(V)}(k+\ix-1)\ran$,
then it is proved in lemma~\ref{hpd_p1} below that $Y = \PP(W^\perp) \subset \PP(V^*)$
is Homologically Projectively Dual to $X$,
where $W^\perp \subset V^*$ is the orthogonal subspace.
\end{remark}

\begin{corollary}
If $g:Y \to \PP(V^*)$ is Homologically Projectively Dual to
$f:X \to \PP(V)$ then
either $\dim X + \dim Y = \dim V - 2$,
or $\dim X + \dim Y \ge \dim V$.
Moreover, the first case is possible only
when $X = \PP(W) \subset \PP(V)$ and $Y = \PP(W^\perp) \subset \PP(V^*)$
where $W \subset V$ is a vector subspace.
\end{corollary}
\begin{proof}
Note that $\ix + \jx \ge \dim V$, hence if $\dim X \ge \ix$ and $\dim Y \ge \jx$
then $\dim X + \dim Y \ge \dim V$. Assume that $\dim X = \ix - 1$.
If $\CA_0 \ne \CA_{\ix-1} = \CA_{\ix-2}$ then $\ix + \jx \ge \dim V + 2$,
hence $\dim X + \dim Y \ge (\ix - 1) + (\jx - 1) = \ix + \jx - 2 \ge \dim V$.
Finally, if $\CA_0 = \CA_{\ix-1}$ then $X = \PP(W)$ and $Y = \PP(W^\perp)$
by proposition~\ref{inddim} and lemma~\ref{hpd_p1}.
\end{proof}

\subsection{Homological Projective Duality and classical projective duality}

Given a projective morphism $f:X \to \PP(V)$ we denote by $X^\vee \subset \PP(V^*)$
the set of all points $H\in\PP(V^*)$ such that the corresponding hyperplane section
$X_H$ of $X$ is singular. It is clear that $X^\vee$ is a Zariski closed
subset in~$\PP(V^*)$. Note that if $f:X \to \PP(V)$ is an embedding
then $X^\vee$ is the classical projectively dual variety to~$X$.

The main result of this subsection is the following

\begin{theorem}\label{hpd_cpd}
Assume that $g:Y \to \PP(V^*)$ is Homologically Projectively Dual
to $f:X \to \PP(V)$. Then the set
$\sing(g) := \{\text{critical values of $g$}\}$
coincides with $X^\vee$, the classical projectively dual variety of $X$.
\end{theorem}
\begin{proof}
Consider the universal hyperplane section $\CX_1$ of $X$ and the maps
$f_1:\CX_1 \to \PP(V^*)$ and $g:Y \to \PP(V^*)$. Note that by definition of
Homological Projective Duality we have a semiorthogonal decomposition~(\ref{dbx1}).
Note also that $X^\vee = \sing(f_1)$ is the set of critical values of the map $f_1$.
Thus we have to check that $\sing(f_1) = \sing(g)$.

First of all, assume that $\sing(g) \not\subset \sing(f_1)$.
Let $H\in\PP(V^*)$ be a point in $\sing(g)$ such that $H\not\in \sing(f_1)$.
Then it is clear that there exists a smooth hypersurface $D \subset \PP(V^*)$
such that $H\in D$, $Y_D := Y\times_{\PP(V^*)} D$ has a singularity over $H$,
and $\dim Y_D = \dim Y - 1$. Let $T = D \setminus \sing(f_1)$.
Then $H \in T$, $Y_T := Y\times_{\PP(V^*)} T$ has a singularity over $H$, and $\dim Y_T = \dim Y - 1$.
On the other hand, $f_1$ is smooth over $T$, hence $\CX_{1T} = \CX_1\times_{\PP(V^*)}T$ is smooth
and both $\CX_{1T}$ and $\CX_{1T}\times_T Y_T$ have expected dimension.
Therefore the base change $T \to \PP(V^*)$ is faithful for the pair $(\CX_1,Y)$
and we obtain by the faithful base change theorem~\ref{phitsod} a semiorthogonal
decomposition
\begin{equation}\label{dbx1t}
\D^b(\CX_{1T}) = \lan \D^b(Y_T),\CA_1(1)\boxtimes\D^b(T),\dots,\CA_{\ix-1}(\ix-1)\boxtimes\D^b(T)\ran
\end{equation}
But category $\D^b(\CX_{1T})$ is $\Ext$-bounded since $\CX_{1T}$ is smooth,
while category $\D^b(Y_T)$ is not $\Ext$-bounded since $Y_T$ is singular (see lemma~\ref{snav_sm}).
This is a contradiction, which shows that we must have an embedding $\sing(g) \subset \sing(f_1)$.

Similarly, assume that $\sing(f_1) \not\subset \sing(g)$.
Let $H\in\PP(V^*)$ be a point in $\sing(f_1)$ such that $H\not\in \sing(g)$.
Then it is clear that there exists a smooth hypersurface $D \subset \PP(V^*)$
such that $H\in D$, $\CX_{1D} := \CX_1\times_{\PP(V^*)} D$ has a singularity over $H$,
and $\dim \CX_{1D} = \dim \CX_1 - 1$. Let $T = D \setminus \sing(g)$.
Then $H \in T$, $\CX_{1T} := \CX_1\times_{\PP(V^*)} T$ has a singularity over $H$, and $\dim \CX_{1T} = \dim \CX_1 - 1$.
On the other hand, $g$ is smooth over $T$, hence $Y_T = Y\times_{\PP(V^*)}T$ is smooth
and both $Y_T$ and $\CX_{1T}\times_T Y_T$ have expected dimension.
Therefore the base change $T \to \PP(V^*)$ is faithful for the pair $(\CX_1,Y)$
and we again obtain by the faithful base change theorem~\ref{phitsod} a semiorthogonal
decomposition~(\ref{dbx1t}). Now we note that category $\D^b(Y_T)$ is $\Ext$-bounded
since $Y_T$ is smooth, and categories $\CA_1(1)\boxtimes\D^b(T)$, \dots,
$\CA_{\ix-1}(\ix-1)\boxtimes\D^b(T)$ are $\Ext$-bounded because $T$ is smooth.
Therefore category $\D^b(\CX_{1T})$ is $\Ext$-bounded by lemma~\ref{abt}.
But this is a contradiction with the fact that $\CX_{1T}$ is singular.
\end{proof}

\subsection{Homological projective duality and triangulated categories of singularities}

Recall that in \cite{O4} to every algebraic variety $X$ there was associated
a triangulated category $\D_\sg(X) := \D^b(X)/\D^\perf(X)$, the quotient
category of the bounded derived category of coherent sheaves
by the subcategory of perfect complexes, which was called
the {\sf triangulated category of singularities}\/ of $X$.
This definition easily generalizes to any triangulated category.

\begin{definition}[\cite{O5}]
Let $\D$ be a triangulated category. An object $F\in\D$ is
{\sf homologically finite}\/ if for any $G\in\D$ the set
$\{n\in\ZZ\ |\ \Hom(F,G[n])\ne 0 \}$ is finite.
\end{definition}

The full subcategory of $\D$ consisting of homologically finite objects
is denoted by $\D_\hf$. It is a triangulated subcategory.
The quotient category $\D_\sg := \D/\D_\hf$ is called
the {\sf triangulated category of singularities}\/ of $\D$.

\begin{lemma}[\cite{O5}]\label{sod_sg}
If $\D = \lan\D_1,\D_2,\dots,\D_m\ran$ is a semiorthogonal decomposition
of $\D$ then its triangulated category of singularities has the following
semiorthogonal decomposition
$$
\D_\sg = \lan\D_{1\sg},\D_{2\sg},\dots,\D_{m\sg}\ran.
$$
\end{lemma}

\begin{theorem}
If $g:Y\to\PP(V^*)$ is Homologically Projectively Dual to $f:X\to\PP(V)$ and
$L\subset V^*$ is an admissible subspace then $\D_\sg(X_L) \cong \D_\sg(Y_L)$.
\end{theorem}
\begin{proof}
Since $X$ is smooth by assumptions and $Y$ is smooth by theorem~\ref{themain}
it follows from lemma~\ref{sod_sg} that
$$
\CA_{0\sg} = \CA_{1\sg} = \dots = \CA_{(\ix-1)\sg} = 0
\qquad\text{and}\qquad
\CB_{0\sg} = \CB_{1\sg} = \dots = \CB_{(\jx-1)\sg} = 0,
$$
where $\D^b(X) = \lan\CA_0,\CA_1(1),\dots,\CA_{\ix-1}(\ix-1)\ran$ and
$\D^b(Y) = \lan\CB_{\jx-1}(1-\jx),\dots,\CB_1(-1),\CB_0\ran$ are
the Lefschetz decompositions of $X$ and $Y$ respectively.
Using again theorem~\ref{themain} and lemma~\ref{sod_sg} we deduce that
$\D_\sg(X_L) \cong (\CC_L)_\sg \cong \D_\sg(Y_L)$.
\end{proof}

\section{Projective bundles}

Let $S$ be a smooth (not necessarily compact) base scheme
with a vector bundle $E$ of rank $\ix$.
Let $X = \PP_S(E)$ be a projectivization of this vector bundle
with the projection $p:X \to S$, and let $\CO_X(1)$ be the Grothendieck
line bundle on $\PP_S(E)$ over $S$ (such that $p_*\CO_X(1) \cong E^*$).
Let $V^* \subset \Gamma(S,E^*) = \Gamma(X,\CO_X(1))$ be a space of
global sections generating $E^*$, and let
$f:X \to \PP(V)$ be the corresponding morphism.
Let $\CA_0 = p^*(\D^b(S)) \subset \D^b(X)$.
Then by the result of Orlov~\cite{O2} we have
a semiorthogonal decomposition
\begin{equation}\label{odec}
\D^b(X) = \langle\CA_0,\CA_1(1),\dots,\CA_{\ix-1}(\ix-1)\rangle
\qquad\text{with}\quad
\CA_0=\CA_1=\dots=\CA_{\ix-1} = p^*(\D^b(S)).
\end{equation}

Let $\CX_1 \subset X\times\PP(V^*)$ be the universal hyperplane section of $X$.
Since $X\times\PP(V^*)$ is the projectivization of the pullback of $E$
to $S\times\PP(V^*)$, it follows that the fiber of $\CX_1$ over the generic
point of $S\times\PP(V^*)$ is a hyperplane in the projectivization of the fiber
of $E$ over the corresponding point of $S$, and over a certain closed subset
of $S\times\PP(V^*)$ the fiber of $\CX_1$ coincides with the whole
projectivization of the fiber of $E$.
This closed subset $Y \subset S\times\PP(V^*)$ is the zero locus of the section
of the vector bundle $E^*\boxtimes\CO_{\PP(V^*)}(1)$, corresponding to
the identity in
$\Gamma(S\times\PP(V^*),E^*\boxtimes\CO_{\PP(V^*)}(1)) =
\Gamma(S,E^*)\otimes\Gamma(\PP(V^*),\CO_{\PP(V^*)}(1)) \cong
V^*\otimes V$.
Note that

\begin{lemma}
We have $Y \cong \PP_S(E^\perp)$, where $E^\perp = \Ker(V^*\otimes\CO_S \to E^*)$.
In particular, $Y$ is smooth and
\begin{equation}\label{yreg}
\codim_{S\times\PP(V^*)} Y = \ix.
\end{equation}
\end{lemma}
\begin{proof}
The fiber of $Y$ over a point $s \in S$ consists of all sections $H \in V^* = H^0(S,E^*)$,
that vanish at $s$, i.e.\ which are contained in the fiber of $E^\perp$ over $S$.
\end{proof}

Let $f:\CX_1 \to S\times\PP(V^*)$ be the canonical projection.
Let $g:Y \to S\times\PP(V^*)$ be the embedding.
Let $Z = Y\times_{(S\times\PP(V^*))}\CX_1$ be the fiber product,
and denote by $\phi:Z \to Y$ and $i:Z \to \CX_1$ the projections.
So we have the following cartesian square
$$
\xymatrix{
\llap{$\PP_Y(E) \cong \,$} Z \ar[rr]^i \ar[d]_\phi && \CX_1 \ar[d]^f \\ Y \ar[rr]^g && S\times\PP(V^*)
}
$$
Let $\CE = i_*\CO_Z$ and consider the corresponding kernel functor
$\Phi_\CE:\D^b(Y) \to \D^b(\CX_1)$. Note that
the functor $\Phi_\CE$ is $(S\times\PP(V^*))$-linear.

The main result of this section is the following

\begin{theorem}\label{dbx1pe}
In the above notation and assumptions we have a semiorthogonal decomposition
$$
\D^b(\CX_1) = \langle\Phi_\CE(\D^b(Y)),
\CA_1(1)\boxtimes\D^b(\PP(V^*)),\dots,\CA_{\ix-1}(\ix-1)\boxtimes\D^b(\PP(V^*))
\rangle.
$$
\end{theorem}

Comparing this theorem with the definition of relative
Homological Projective Duality (see remark~\ref{relhpd})
we obtain the following

\begin{corollary}\label{hpd_pe}
If $E$ is generated by global sections then $Y = \PP_S(E^\perp)$
is Homologically Projectively Dual to $X = \PP_S(E)$ relatively over $S$.
\end{corollary}

We start the proof with some preparations.

\begin{lemma}
The subscheme $Z \subset \CX_1$ is a zero locus of a section
of the vector bundle $\Omega_{X/S}(1)\otimes\CO_{\PP(V^*)}(1)$ on $\CX_1$,
where $\Omega_{X/S}$ is the sheaf of relative differentials.
\end{lemma}
\begin{proof}
Note that it follows from the definitions that
$\CX_1$ is a zero locus of a section of the line bundle
$\CO_X(1)\boxtimes\CO_{\PP(V^*)}(1)$ on $X\times\PP(V^*)$ and
$Z$ is a zero locus of a section of the vector bundle
$E^*\boxtimes\CO_{\PP(V^*)}(1)$ on $X\times\PP(V^*)$.
Moreover, it is clear that the canonical epimorphism
$E^*\boxtimes\CO_{\PP(V^*)}(1) \to \CO_X(1)\boxtimes\CO_{\PP(V^*)}(1)$
takes the latter section to the former. Therefore, the latter section
restricted to the zero locus of the former section is contained
in the kernel of the epimorphism which is isomorphic to
$\Omega_{X/S}(1)\boxtimes\CO_{\PP(V^*)}(1)$, and its zero locus
coincides with $Z$.
\end{proof}

\begin{corollary}\label{ioz}
The sheaf $\CE = i_*\CO_Z$ is quasiisomorphic
to the Koszul complex,
$$
i_*\CO_Z \cong
\Kosz(T_{X/S}(-1)\otimes\CO_{\PP(V^*)}(-1)) :=
\Lambda^\bullet(T_{X/S}(-1)\otimes\CO_{\PP(V^*)}(-1)).
$$
In particular, $i_*\CO_Z$ is a perfect complex on $\CX_1$
and $i:Z \to \CX_1$ has finite $\Tor$-dimension.
\end{corollary}
\begin{proof}
Since $\CX_1$ has pure codimension $1$ in $X\times\PP(V^*)$ it follows
from~(\ref{yreg}) that $Z$ has codimension $\ix-1$ in $\CX_1$,
therefore the corresponding section of the vector bundle
$\Omega_{X/S}(1)\otimes\CO_{\PP(V^*)}(1)$ is regular
and $i_*\CO_Z$ is quasiisomorphic to the Koszul complex.
In particular it is a perfect complex.
Finally, it follows from the projection formula that
$i_*i^*F \cong F\otimes i_*\CO_Z$ for any $F\in\D^b(\CX_1)$,
and since $i_*$ is exact and conservative it follows that $\Tor$-dimension
of $i$ is finite.
\end{proof}

\begin{proposition}\label{phiceisff}
The functor $\Phi_\CE$ is fully faithful.
\end{proposition}
\begin{proof}
Note that $\CE$ is perfect on $Z$, $Z$ is projective over $Y$ and over $\CX_1$,
and $\CE$ has finite $\Ext$-amplitude over $\CX_1$ by corollary~\ref{ioz}.
Therefore by lemma~\ref{ladj} the kernel functor
$\Phi_{\CE^\#}:\D^b(\CX_1)\to\D^b(Y)$ with the kernel
\begin{multline*}
\CE^\# =
\RCHom((\phi\times i)_*\CO_Z,\omega_{\CX_1\times Y/Y}[\dim\CX_1]) \cong
\\ \cong
(\phi\times i)_*\RCHom(\CO_Z,(\phi\times i)^!(\omega_{\CX_1\times Y/Y}[\dim\CX_1])) \cong
(\phi\times i)_*\omega_{Z/Y}[\dim Z/Y]
\end{multline*}
is left adjoint to $\Phi_\CE$. Further, the composition
$\Phi_{\CE^\#}\circ\Phi_\CE:\D^b(Y) \to D^b(Y)$ is given by
$$
F \mapsto \phi_*(i^*i_*\phi^*F \otimes \omega_{Z/Y}[\dim Z/Y]).
$$
Let $K\in\D^b(Z\times Z)$ denote the kernel of the functor
$i^*i_*:\D^b(Z) \to \D^b(Z)$. It is clear that
$(1\times i)_*K \in \D^b(Z\times\CX_1)$
is a kernel of the functor $i_*i^*i_*:\D^b(Z) \to \D^b(\CX_1)$.
But
$$
i_*i^*i_*F \cong i_*F\otimes i_*\CO_Z \cong i_*(F\otimes i^*i_*\CO_Z),
$$
therefore $(1\times i)_*K \cong \Gamma^i_*i^*i_*\CO_Z$,
where $\Gamma^i:Z \to Z\times\CX_1$ is the graph of $i$.
On the other hand, from corollary~\ref{ioz} it follows that
$i^*i_*\CO_Z$ is isomorphic to
$$
i^*\Kosz(T_{X/S}(-1)\otimes\CO_{\PP(V^*)}(-1)) \cong
\oplus\, i^*\Lambda^t T_{X/S}(-t) \otimes \CO_{\PP(V^*)}(-t)[t] \cong
\oplus\, \Lambda^t T_{Z/Y}(-t) \otimes \CO_{\PP(V^*)}(-t)[t].
$$
Therefore
$$
\CH^{-t}(K) \cong \Delta_*(\Lambda^t T_{Z/Y}(-t) \otimes \CO_{\PP(V^*)}(-t))
\qquad
\text{for all $t$.}
$$
where $\Delta:Z \to Z\times Z$ is the diagonal embedding.

Consider the functor given by the kernel $\CH^{-t}(K)$:
\begin{multline}\label{chtf}
\Phi_{\CH^{-t}(K)}:\D^b(Z) \to \D^b(Z),
\qquad
F \mapsto
\phi_*(\phi^*F \otimes \Lambda^t T_{Z/Y}(-t) \otimes \CO_{\PP(V^*)}(-t) \otimes \omega_{Z/Y}[\dim Z/Y]) \cong
\\ \cong
F \otimes \phi_*(\Lambda^t T_{Z/Y}(-t) \otimes \CO_{\PP(V^*)}(-t) \otimes \omega_{Z/Y}[\dim Z/Y]).
\end{multline}
Note that
\begin{multline*}
\phi_*(\Lambda^t T_{Z/Y}(-t) \otimes \CO_{\PP(V^*)}(-t) \otimes \omega_{Z/Y}[\dim Z/Y]) \cong
\\ \cong
\phi_*\RCHom(\Omega^t_{Z/Y}(t),\phi^!(\CO_{\PP(V^*)}(-t))) \cong
\RCHom(\phi_*(\Omega^t_{Z/Y}(t)),\CO_{\PP(V^*)}(-t))),
\end{multline*}
but for $t\ne 0$ we have $\phi_*(\Omega^t_{Z/Y}(t)) = 0$, while
for $t = 0$ we have $\phi_*(\Omega^t_{Z/Y}(t)) = \CO_Y$.
It follows that the functor~(\ref{chtf}) is zero for $t\ne 0$ and is identity for $t=0$.
Using a devissage argument we deduce that the functor
$$
F \mapsto
\phi_*(\phi^*F \otimes K \otimes \omega_{Z/Y}[\dim Z/Y]) =
\phi_*(i^*i_*\phi^*F \otimes \omega_{Z/Y}[\dim Z/Y])
$$
is the identity. But as we have seen above this is the composition
of $\Phi_\CE$ with its left adjoint. Therefore $\Phi_\CE$ is fully faithful.
\end{proof}

\begin{lemma}\label{fgip}
The functor $f^*g_*:\D^b(Y) \to \D^b(\CX_1)$ is a kernel functor,
$f^*g_*\cong\Phi_{K(f,g)}$, and its kernel fits into exact triangle
$K(f,g) \to \CO_Z \to \CO_Z\otimes(\alpha^*\CO_X(-1)\otimes\phi^*g^*\CO_{\PP(V^*)}(-1))[2]$.
In particular, for any $F\in\D^b(Y)$ we have an exact triangle in $\D^b(\CX_1)$
$$
f^*g_*F \to i_*\phi^*F \to i_*\phi^*(F\otimes g^*\CO_{\PP(V^*)}(-1))\otimes \CO_X(-1)[2]
$$
\end{lemma}
\begin{proof}
Consider the following commutative diagram
$$
\xymatrix{
&&& \CX_1 \ar[dl]^\alpha \ar[ddl]^f \\
Z \ar@{}[drr]|\ecart \ar[urrr]^i \ar[rr]^{\alpha\circ i} \ar[d]_\phi && X\times\PP(V^*) \ar[d]^p \\
Y \ar[rr]^g && S\times\PP(V^*)
}
$$
Note that $Z = Y\times_{(S\times\PP(V^*))}(X\times\PP(V^*))$ and the square
in the diagram is exact cartesian, because $p$ is flat. It follows that
$$
f^*g_* =
\alpha^*p^*g_* =
\alpha^*(\alpha\circ i)_*\phi^* =
\alpha^*\alpha_*i_*\phi^*.
$$
On the other hand, $\alpha$ is a divisorial embedding and $\CX_1$
is a zero locus of the line bundle $\CO_X(1)\boxtimes\CO_{\PP(V^*)}(1)$,
therefore $\alpha^*\alpha_*$ is a kernel functor and its kernel
$K_\alpha\in\D^b(\CX_1\times\CX_1)$ fits into exact triangle
$$
K_\alpha \to \Delta_*\CO_{\CX_1} \to
\Delta_*\alpha^*(\CO_X(-1)\otimes\CO_{\PP(V^*)}(-1))[2].
$$
Computing the convolution of this triangle with the kernel $\CO_Z$
of the functor $i_*\phi^*$ we obtain the claim.
\end{proof}

\begin{proofof}{of theorem~\ref{dbx1pe}}
Note that we have
$\CA_k(k)\boxtimes\D^b(\PP(V^*)) = f^*\D^b(S\times\PP(V^*))\otimes\CO_X(k)$,
and that
$\Phi_\CE(\D^b(Y)) = i_*\phi^*(\D^b(Y))$, so we need to prove that we have
the following semiorthogonal decomposition
\begin{equation}\label{dbcx1}
\hspace{-30pt}
\D^b(\CX_1) = \langle i_*\phi^*(\D^b(Y)),
f^*(\D^b(S\times\PP(V^*)))\otimes\CO_X(1),\dots,
f^*(\D^b(S\times\PP(V^*)))\otimes\CO_X(\ix-1)
\rangle.
\end{equation}
The semiorthogonality in question is verified quite easily.
Using lemma~\ref{sodx1} and proposition~\ref{phiceisff},
the question reduces to check that
$\Hom(f^*F\otimes\CO_X(k),i_*\phi^*G) = 0$ for all
$F\in\D^b(S\times\PP(V^*))$, $G\in\D^b(Y)$ and $1\le k\le \ix-1$.
But
\begin{multline*}
\Hom(f^*F\otimes\CO_X(k),i_*\phi^*G) =
\Hom(f^*F,(i_*\phi^*G)\otimes\CO_X(-k)) =
\Hom(f^*F,i_*(\phi^*G\otimes\CO_Z(-k))) = \\ =
\Hom(F,f_*i_*(\phi^*G\otimes\CO_Z(-k))) =
\Hom(F,g_*\phi_*(\phi^*G\otimes\CO_Z(-k))) =
\Hom(F,g_*(G\otimes\phi_*\CO_Z(-k)))
\end{multline*}
and $\phi_*\CO_Z(-k) = 0$ for $1\le k\le \ix-1$ since
$\phi$ is a projectivization of a rank $\ix$ vector bundle.

It remains to check that the RHS of~(\ref{dbcx1}) generates $\D^b(\CX_1)$.
First of all we will show that it contains $i_*\phi^*(\D^b(Y))\otimes\CO_X(k)$
for $0\le k\le \ix-1$. Indeed, for any $F\in\D^b(Y)$ the triangle
of lemma~\ref{fgip} twisted by $\CO_X(1)$ takes form
$$
f^*g_*F\otimes\CO_X(1) \to
i_*\phi^*F\otimes\CO_X(1) \to
i_*\phi^*(F\otimes g^*\CO_{\PP(V^*)}(-1))[2].
$$
Note that its first term is contained in $f^*(\D^b(S\times\PP(V^*)))\otimes\CO_X(1)$
and the last term is contained in $i_*\phi^*(\D^b(Y))$, therefore
the middle term is contained in the RHS of~(\ref{dbcx1}).
Twisting this triangle by $\CO_X(2)$ we get
$$
f^*g_*F\otimes\CO_X(2) \to
i_*\phi^*F\otimes\CO_X(2) \to
i_*\phi^*(F\otimes g^*\CO_{\PP(V^*)}(-1))\otimes\CO_X(1)[2].
$$
Note that its first term is contained in $f^*(\D^b(S\times\PP(V^*)))\otimes\CO_X(2)$
and the last term is contained in $i_*\phi^*(\D^b(Y))\otimes\CO_X(1)$, therefore
the middle term is contained in the RHS of~(\ref{dbcx1}).
Continuing in this way we deduce the claim by induction.

Now assume that $G$ is contained in the left orthogonal
to the RHS of~(\ref{dbcx1}). Then it follows from above that for any $0\le k\le \ix-1$
and any $F\in\D^b(Y)$ we have $\Hom(G,i_*\phi^*F\otimes\CO_X(k)) = 0$.
But by adjunction it equals to $\Hom(i^*G,\phi^*F\otimes\CO_Z(k))$
and since by \cite{O2} the subcategories $\phi^*(\D^b(Y))\otimes\CO_Z(k)$
with $0\le k\le \ix-1$ generate $\D^b(Z)$ it follows that $i^*G = 0$.
This means that $G$ is supported on $\CX_1\setminus Z$.
But $\CX_1\setminus Z$ is a $\PP^{\ix-2}$-bundle
over $(S\times\PP(V^*))\setminus Y$, hence the orthogonality of $G$ to the subcategory
$\lan
f^*(\D^b(S\times\PP(V^*)))\otimes\CO_X(1),\dots,
f^*(\D^b(S\times\PP(V^*)))\otimes\CO_X(\ix-1)
\ran$
implies that $G = 0$.
\end{proofof}


Now we are going to apply theorems~\ref{themain} and \ref{refthemain}
for this special case of Homological Projective Duality.

Let $F$ be another vector bundle on $S$, $\rank(F) = r$,
and let $\phi:F \to E^*$ be a morphism of vector bundles.
Consider the projectivizations $\PP_S(E)$ and $\PP_S(F)$.
Let $p:\PP_S(E) \to S$ and $q:\PP_S(F) \to S$ be the projections
and let $\CO_{\PP_S(E)/S}(1)$ and $\CO_{\PP_S(F)/S}(1)$ denote
the Grothendieck ample line bundles.
Note that $\phi$ induces a section of the vector bundle
$F^*\boxtimes\CO_{\PP_S(E)/S}(1)$ on $\PP_S(E)$, and a section
of the vector bundle $E^*\boxtimes\CO_{\PP_S(F)/S}(1)$ on $\PP_S(F)$.
Let $X_F \subset \PP_S(E)$ and $Y_F \subset \PP_S(F)$ denote their zero loci.

\begin{theorem}\label{ef}
If\/
$\codim_{\PP_S(E)} X_F = \rank F$,
$\codim_{\PP_S(F)} Y_F = \rank E$ and
$\dim X_F\times_S Y_F = \dim S - 1$
then there exist semiorthogonal decompositions
\begin{equation}\label{xyf}
\arraycolsep = 2pt
\hspace{-30pt}
\begin{array}{llll}
\D^b(X_F) &=&  \langle \D^b(Y_F),
p^*\D^b(S)\otimes\CO_{\PP_S(E)/S}(1),\dots,
p^*\D^b(S)\otimes\CO_{\PP_S(E)/S}(\ix-r)\rangle, &
\text{if $r < \ix$}\\
\D^b(Y_F) &=&  \langle
q^*\D^b(S)\otimes\CO_{\PP_S(F)/S}(\ix-r),\dots,
q^*\D^b(S)\otimes\CO_{\PP_S(F)/S}(-1),
\D^b(X_F) \rangle, \ &
\text{if $r > \ix$},
\end{array}
\end{equation}
and an equivalence $\D^b(X_F) \cong \D^b(Y_F)$ if $r = \ix$.
\end{theorem}
\begin{proof}
First of all, consider the case when the morphism $\phi:F \to E^*$
can be represented as a composition $F \to V^*\otimes \CO_S \to E^*$
of a monomorphism of vector bundles to a trivial vector bundle
followed by an epimorphism of vector bundles.
In this case the claim of theorem follows from theorem~\ref{refthemain}.

Indeed, take $r = \rank F$, $T = \Gr_S(r,V^*\otimes\CO_S)$, the relative Grassmannian,
and let $\CL \subset V^*\otimes\CO_T$ be the tautological subbundle of rank $r$.
Then $\CX_\CL$ and $\CY_\CL$ in the notations of theorem~\ref{refthemain}
are the universal families of linear sections of $X = \PP_S(E)$ and $Y = \PP_S(E^\perp)$
respectively. It is easy to see that the dimension assumptions of theorem~\ref{refthemain}
are satisfied, hence we have semiorthogonal decompositions
$$
\arraycolsep = 2pt
\begin{array}{llll}
\D^b(\CX_\CL) &=&  \langle \D^b(\CY_\CL),
p^*\D^b(T)\otimes\CO_{\PP_S(E)/S}(1),\dots,
p^*\D^b(T)\otimes\CO_{\PP_S(E)/S}(\ix-r)\rangle, &
\text{if $r < \ix$}\\
\D^b(\CY_\CL) &=&  \langle
q^*\D^b(T)\otimes\CO_{\PP_S(F)/S}(\ix-r),\dots,
q^*\D^b(T)\otimes\CO_{\PP_S(F)/S}(-1),
\D^b(\CX_\CL) \rangle, \ &
\text{if $r > \ix$},
\end{array}
$$
and an equivalence $\D^b(\CX_\CL) \cong \D^b(\CY_\CL)$, if $r = \ix$.

The embedding $F \to V^*\otimes\CO_S$ gives a section $\sigma:S \to T = \Gr_S(V^*\otimes\CO_S)$
such that $F \cong \sigma^*\CL$. Consider $\sigma$ as a base change.
Note that $\CX_\CL\times_T S = X_F$. On the other hand,
$\CY_\CL\times_T S$ is the zero locus of a section of the vector bundle
$\CO_{\PP_S(E^\perp)/S}(1)\otimes F^{\perp*}$ on $\PP_S(E^\perp)$.
But looking at the commutative diagram
$$
\xymatrix{
E^\perp\boxtimes\CO_S \ar[r]^{\phi^\perp} &
\CO_S\boxtimes F^{\perp*} \\
E^\perp\boxtimes\CO_S \ar[r] \ar@{=}[u] &
V^*\otimes\CO_S\boxtimes\CO_S \ar[r] \ar[u] &
E^*\boxtimes\CO_S \\
&
\CO_S\boxtimes F \ar@{=}[r] \ar[u] &
\CO_S\boxtimes F \ar[u]_\phi
}
$$
it is easy to deduce that it also can be represented as the zero locus of a section
of the vector bundle $E^*\otimes\CO_{\PP_S(F)/S}(1)$ on $\PP_S(F)$,
i.e.\ that $\CY_\CL\times_T S = Y_F$.

Finally note that the dimension assumptions of the theorem and lemma~\ref{fisf}
imply that this base change is faithful for a pair $(\CX_\CL,\CY_\CL)$.
Applying the faithful base change theorem~\ref{phitsod} we deduce the claim.

The general case follows from the above case by theorem~\ref{loc_sod}.
Indeed, all inclusion functors in the desired decompositions
are $S$-linear and for every point $s\in S$ there exists
an open neighborhood $U \subset S$ over which the morphism
$\phi_{|U}:F_{|U} \to E^*_{|U}$ can be represented as a composition
$F_{|U} \to V^*\otimes \CO_U \to E^*_{|U}$ of a monomorphism of vector bundles
followed by an epimorphism of vector bundles.
\end{proof}

Consider the case $r = \ix$.
Then the dimension assumptions of theorem~\ref{ef} can be rewritten as
$$
\dim X_F = \dim Y_F = \dim X_F\times_S Y_F = \dim S - 1.
$$
Note that $D := p(X_F) = q(Y_F) \subset S$ is the degeneration locus
of $\phi$, i.e. the zero locus of $\det \phi :\det F \to \det E^*$.
Since $\dim D = \dim S - 1$ it follows that $p:X_F \to D$ and
$q:Y_F \to D$ are birational. Therefore $q^{-1}\circ p$ is a birational
transformation $\xymatrix@1{X_F \ar@{-->}[r] & Y_F}$.
It is easy to check that
$\omega_{X_F} \cong \omega_{Y_F} \cong \omega_D \cong \omega_S\otimes\det E^*\otimes\det F^*$,
hence this transformation is a flop.

\begin{corollary}
If $\dim X_F = \dim Y_F = \dim X_F\times_S Y_F = \dim S - 1$ then
the kernel functor with kernel given by the structure sheaf
of the fiber product $X_F\times_D Y_F = X_F\times_S Y_F$ is an equivalence
of categories $\D^b(X_F) \cong \D^b(Y_F)$.
\end{corollary}

For example, let $S = \PP^4$, $E = \CO_S(-1) \oplus \CO_S(-1)$,
$F = \CO_S \oplus \CO_S$ and $\phi:F \to E^*$
given by the matrix
$\phi = \left(\begin{smallmatrix}x & y\\z & u\end{smallmatrix}\right)$,
where $(x:y:z:u:v)$ are the homogeneous coordinates on  $S$.
Then $D \subset S$ is the cone over $\PP^1\times\PP^1$
(given by equation $xu - yz = 0$), $X_F$ and $Y_F$ are small
resolutions of $D$ and $\xymatrix@1{X_F \ar@{-->}[r] & Y_F}$
is the standard flop.

\section{Examples}

In this section we will give several examples of
Homologically Projectively Dual varieties.

%
%


\subsection{A stupid example}

Let $X$ be any smooth algebraic variety with a projective map $f:X \to \PP(V)$.
Then taking $\CA_0 = \D^b(X)$ we get a Lefschetz decomposition with only one term.
We will call it the {\sf stupid Lefschetz decomposition}.

\begin{proposition}\label{stupid}\label{hpd_stupid}
The universal hyperplane section $\CX_1 \subset X\times\PP(V^*)$
with the projection $\CX_1 \to \PP(V^*)$ is
Homologically Projectively Dual to $X \to \PP(V)$
with respect to the stupid Lefschetz decomposition.
\end{proposition}
\begin{proof}
Just take $\CE$ to be the structure shift of the diagonal in $\CX_1\times\CX_1$.
\end{proof}

Alternatively, one can consider the stupid Lefschetz decomposition as a particular case
of the decomposition~(\ref{odec}) since any algebraic variety can be considered as
a projectivization of a line bundle over itself.

Let us describe the claim of theorem~\ref{themain} in this case.
Let $Y = \CX_1$ be the universal hyperplane section of $X$.
Let $L \subset V^*$ be a vector subspace.
Then $Y_L$ is the pencil of hyperplane sections of $X$ parameterized by
$\PP(L)$. It is fibered over $X$ with fiber equal to $\PP(L)$ over $X_L$,
and a hyperplane in $\PP(L)$ over $X \setminus X_L$.
Theorem~\ref{themain} implies that we have
a semiorthogonal decomposition
$$
\D^b(Y_L) = \lan \D^b(X)\otimes\CO_{\PP(L)}(1-\dim L),\dots,
\D^b(X)\otimes\CO_{\PP(L)}(-1),\D^b(X_L)\ran.
$$
E.g.\ for $\dim L = 1$ we have $Y_L = X_L$ and
for $L = 2$ we have $Y_L$ is the blowup of $X_L$ in $X$.
In the latter case the obtained semiorthogonal decomposition
coincides with the standard decomposition of the blowup.

\subsection{Curves}

Now, assume that $X$ is a smooth projective curve.

\begin{lemma}\label{ldcurve}
Derived category of a smooth connected projective curve $X$ admits a nontrivial
Lefschetz decomposition with respect to an effective line bundle $\CO_X(1)$,
only if $X \cong \PP^1$. In this case $\CO_X(1)$ is the positive generator of $\Pic X$,
and the decomposition takes form $\D^b(X) = \lan \CO_{\PP^1}(k),\CO_{\PP^1}(k+1)\ran$
for some $k$.
\end{lemma}
\begin{proof}
Assume that we have a nontrivial Lefschetz decomposition of $\D^b(X)$,
so that $\CA_1 \ne 0$. Let $F$ be a nontrivial object in $\CA_1$.
By definition of a Lefschetz decomposition we have
$$
\RHom(F(1),F) = 0.
$$
Since $X$ is a curve, every object in $\D^b(X)$
is a direct sum of its cohomology sheaves, and every sheaf on $X$
is a direct sum of a torsion sheaf and of a locally free sheaf.
If $G$ is a nontrivial torsion sheaf some shift of which is a direct summand of $F$,
then $G(1) \cong G$, hence $\Hom(G(1),G)\ne 0$, hence $\Hom(F(1),F)\ne 0$.
Therefore $F$ is a direct sum of shifts of locally free sheaves.

Since $\CA_1$ is closed under direct summands and triangulated,
there exists a locally free sheaf $F \in \CA_1$. Then
$$
\RHom(F(1),F) = \RG(X,F\otimes F^*(-1)).
$$
But $F\otimes F^*$ has $\CO_X$ as a direct summand, hence
the condition $\RHom(F(1),F) = 0$ implies that the line bundle
$\CO_X(-1)$ on $X$ has no cohomology. By Riemann--Roch this is possible
only if $\deg\CO_X(1) = 1 - g$, where $g$ is the genus of $X$.
So, if $g\ge 1$ then $\CO_X(1)$ cannot be effective.
Therefore for $g \ge 1$ we cannot have a nontrivial Lefschetz decomposition.

Now assume that $g = 0$, so $X \cong \PP^1$. Then the above arguments show
that $\CO_X(1)$ is the positive generator of the $\Pic X$. Moreover,
since any locally free sheaf on $\PP^1$ is a direct sum of line bundles,
it follows that $\CA_1 = \lan\CO_{\PP^1}(k)\ran$ for some $k \in \ZZ$.
Then
$$
\lan\CO_{\PP^1}(k)\ran = \CA_1 \subset \CA_0 \subset \CA_1(1)^\perp =
\lan\CO_{\PP^1}(k+1)\ran^\perp = \lan\CO_{\PP^1}(k)\ran,
$$
and we are done.
\end{proof}

The above lemma shows that the only way to get a Homological Projective Duality
for a curve of positive genus is to consider the stupid Lefschetz decomposition.
Then as we have shown in proposition~\ref{stupid} the Homologically Projectively
Dual variety is the universal hyperplane section. Note that in this case
the map $\CX_1 \to \PP(V^*)$ is a finite covering (of degree equal
to the degree of $X$ in $\PP(V)$) ramified over the classical
projectively dual hypersurface $X^\vee \subset \PP(V^*)$.

The case of $X = \PP^1$ is treated in the following lemma.

\begin{lemma}\label{hpd_p1}
If $X = \PP(W) \subset \PP(V)$, a linear subspace, $\dim W = \ix$,
considered with the Lefschetz decomposition
$\D^b(X) = \lan f^*\CO_{\PP(V)}(k),\dots,f^*\CO_{\PP(V)}(k+\ix-1)\ran$,
then $Y = \PP(W^\perp) \subset \PP(V^*)$ is Homologically Projectively Dual to $X$.
\end{lemma}
\begin{proof}
We can consider $X$ as a projectivization of a trivial vector bundle $W$
over $\Spec\kk$. Then the claim follows from corollary~\ref{hpd_pe}.
\end{proof}

\subsection{Hirzebruch surfaces}

Let $S = \PP^1$ and $E = \CO_S \oplus \CO_S(-d)$, so that
$X = \PP_S(E)$ is the Hirzebruch surface $F_d$.
Take $V^* = H^0(S,E^*) \cong \kk \oplus \kk^{d+1}$.
Then $f:X \to \PP(V) = \PP^{d+1}$ maps $X$ onto
a cone over a Veronese rational curve of degree $d$
(the exceptional section of $X$ is contracted
to the vertex of the cone).

In this case $E^\perp = \Ker(V^* \to E^*) \cong \CO_S(-1)^d$,
hence $Y = \PP_S(E^\perp) \cong \PP^1\times\PP^{d-1}$.
The map $g:Y \to \PP(V^*)$ is a $d$-fold covering
onto the hyperplane in $\PP^d \subset \PP(V^*)$,
corresponding to the vertex of the cone.

\subsection{Two-dimensional quadric}

Let $S = \PP^1$ and $E = \CO_S(-1) \oplus \CO_S(-1)$,
so that $X = \PP_S(E) \cong \PP^1\times\PP^1$
is the two dimensional quadric.
Take $V^* = H^0(S,E^*) \cong \kk^4$.
Then $f:X \to \PP(V) = \PP^3$ is the standard embedding.

In this case
$E^\perp = \Ker(V^* \to E^*) \cong \CO_S(-1)\oplus\CO_S(-1)$,
hence $Y$ is also isomorphic to $\PP^1\times\PP^1$ and
the map $g:Y \to \PP(V^*) = \PP^3$ identifies it
with the projectively dual quadric to $X$.

In a forthcoming paper~\cite{K4} we will describe Homological Projective Duality
for all quadrics.

\subsection{Springer--Grothendieck resolution}

Let $G$ be a semisimple algebraic group,
$S = G/B$ be the flag variety of $G$
(the set of all Borel subgroups in $G$),
$\fg$ isbethe Lie algebra of $G$,
$\fb \subset \fg\otimes\CO_S$,
(resp.\ $\fn \subset \fg\otimes\CO_S$)
be the vector subbundle with fiber over a point of $G/B$
given by the corresponding Borel subalgebra
(resp.\ nilpotent subalgebra) of $\fg$.
Take $E = \fn$ (by the way $\fn$ is isomorphic
to the cotangent bundle of $S$), so that
$X = \PP_S(E) = \PP_{G/B}(\fn) \cong \PP_{G/B}(T^*_{G/B})$,
and $V^* = \fg^* \cong H^0(G/B,\fn^*)$.
Then $f:X \to \PP(\fg)$ maps $X$ onto the projectivization
of the nilpotent cone in $\fg$ and is well known
as the (projectivized) Springer resolution of the nilpotent cone.

In this case
$E^\perp = \Ker(\fg^* \to \fn^*) \cong \Ker(\fg \to \fg/\fb) = \fb$
(we identify $\fg$ with $\fg^*$ by the Killing form),
hence $Y$ is isomorphic to $\PP_{G/B}(\fb)$ and
the map $g:Y \to \PP(\fg)$ is known
as the (projectivized) simultaneous Springer--Grothendieck resolution.
Its generic fiber consists of $|W|$ points where $W$ is
the Weyl group of~$G$.

%

\end{document}